\documentclass[a4paper]{article}
\makeatletter
\renewcommand\normalsize{\@setfontsize\normalsize{14.4pt}{17pt}\abovedisplayskip 14pt plus3pt minus7pt\abovedisplayshortskip 0pt plus4pt\belowdisplayshortskip 7pt plus4pt minus3pt\belowdisplayskip\abovedisplayskip}
\normalsize
\def\@listi{\leftmargin\leftmargini \parsep 6pt plus2.5pt minus1pt \topsep 12pt plus5pt minus7pt \itemsep 6pt plus2.5pt minus1pt}
\let\@listI\@listi
\@listi
\renewcommand\@pnumwidth{1.55em}
\renewcommand\@tocrmarg{2.55em}
\makeatother
\usepackage[T2A,T1]{fontenc}
\usepackage[english,russian]{babel}

\usepackage[a4paper,left=3cm,right=1cm,top=2cm,bottom=2.5cm,includefoot=false]{geometry}
\usepackage{setspace}
\usepackage{amsmath}
\usepackage{newtxmath}
\usepackage{graphicx}
\usepackage{float}
\usepackage{longtable,booktabs,array,multirow,calc}
\usepackage{caption}
\usepackage{cite}
\usepackage[hidelinks,unicode]{hyperref}
\usepackage{titlesec}
\usepackage{indentfirst}
\usepackage[final,expansion=true,protrusion=true]{microtype}
\SetExpansion{encoding=T2A,family=Tempora-TLF}{stretch=40,shrink=40,step=1}
\usepackage{etoolbox}
\usepackage{fancyhdr}
\usepackage{latexml}
\iflatexml
\usepackage{amsaddr}
\fi
\newcommand{\dualcite}[3]{\iflatexml#3\else#2\fi}
\def\bibnumhtml ref#1\bibnumend{\hyperlink{bibref#1}{#1}}
\newcommand{\bibnum}[1]{\iflatexml\expandafter\bibnumhtml#1\bibnumend\else\citen{#1}\fi}

\DeclareMathSizes{14.4}{14}{9.8}{7}
\iflatexml
\else
\let\standardfrac\frac
\newcommand{\mainfrac}[2]{
  \genfrac{}{}{}{}{\displaystyle{\let\frac\standardfrac #1}}{\displaystyle{\let\frac\standardfrac #2}}
}
\fi
\AtBeginEnvironment{equation}{
  \iflatexml\else\let\frac\mainfrac\fi
  \setlength{\abovedisplayskip}{0.09\baselineskip}
  \setlength{\belowdisplayskip}{0.78\baselineskip}
  \setlength{\abovedisplayshortskip}{0.09\baselineskip}
  \setlength{\belowdisplayshortskip}{0.78\baselineskip}
}
\iflatexml
\else
\titleformat{\section}{\normalfont\bfseries\centering\fontsize{14}{21}\selectfont}{}{0pt}{}
\titleformat{\subsection}{\normalfont\bfseries\centering\fontsize{14}{18.2}\selectfont}{}{0pt}{}
\titlespacing*{\section}{0pt}{1.0\baselineskip}{0.6\baselineskip}
\titlespacing*{\subsection}{0pt}{0.8\baselineskip}{0.5\baselineskip}
\fi
\DeclareCaptionFont{figlinespacing}{\setstretch{1.3}}
\DeclareCaptionFont{tablelinespacing}{\setstretch{1.3}}
\AtBeginEnvironment{figure}{\vspace{0.38\baselineskip}}
\newcommand{\fullwidthnumitem}[2]{
  \iflatexml
    \item\relax\par #2\par
  \else
    \par\noindent
    \begingroup
    \parshape=2 1.25cm \dimexpr\linewidth-1.25cm\relax 0pt \linewidth
    \makebox[0.75cm][l]{#1.}#2\par
    \endgroup
  \fi
}
\newcommand{\tabledataA}[2]{
\multirow{9}{#1}{1} & 0.0782 & 38.1 & 34.976 & 19.102 & 18.954 & 14.0 \\
& 0.1229 & 49.7 & 35.025 & 23.031 & 23.530 & 21.1 \\
& 0.1706 & 60.4 & 35.018 & 26.727 & 27.681 & 28.2 \\
& 0.2217 & 73.1 & 34.979 & 31.753 & 32.644 & 35.3 \\
& 0.2764 & 79.7 & 35.004 & 34.079 & 34.810 & 42.5 \\
& 0.3354 & 81.8 & 30.423 & 34.489 & 35.016 & 49.6 \\
& 0.3996 & 83.3 & 27.764 & 34.643 & 34.993 & 56.7 \\
& 0.4701 & 84.9 & 23.315 & 34.838 & 34.998 & 63.8 \\
& 0.5489 & 86.0 & 20.313 & 34.842 & 35.001 & 70.9 \\
#2
\multirow{9}{#1}{1.3} & 0.1073 & 37.0 & 35.021 & 18.339 & 18.397 & 14.0 \\
& 0.1687 & 50.7 & 35.016 & 23.185 & 23.904 & 21.1 \\
& 0.2345 & 64.3 & 34.993 & 28.374 & 29.263 & 28.2 \\
& 0.3049 & 78.9 & 33.030 & 34.215 & 35.003 & 35.3 \\
& 0.3807 & 80.7 & 29.777 & 34.450 & 34.985 & 42.5 \\
& 0.4625 & 82.4 & 25.154 & 34.713 & 35.017 & 49.6 \\
& 0.5519 & 83.8 & 22.948 & 34.836 & 34.998 & 56.7 \\
& 0.6507 & 84.9 & 20.436 & 34.836 & 34.999 & 63.8 \\
& 0.7620 & 86.0 & 20.125 & 34.835 & 34.999 & 70.9 \\
}
\newcommand{\tabledataB}[2]{
\multirow{9}{#1}{1.6} & 0.1250 & 39.1 & 34.963 & 19.253 & 19.213 & 14.0 \\
& 0.1961 & 55.7 & 34.975 & 25.322 & 26.011 & 21.1 \\
& 0.2719 & 75.1 & 34.979 & 32.577 & 33.887 & 28.2 \\
& 0.3526 & 79.5 & 28.567 & 33.913 & 35.019 & 35.3 \\
& 0.4389 & 81.2 & 25.438 & 34.131 & 34.998 & 42.5 \\
& 0.5314 & 82.8 & 20.875 & 34.363 & 35.007 & 49.6 \\
& 0.6313 & 84.3 & 20.030 & 34.566 & 34.989 & 56.7 \\
& 0.7402 & 85.8 & 20.346 & 34.801 & 35.018 & 63.8 \\
& 0.8606 & 86.8 & 20.654 & 34.825 & 34.988 & 70.9 \\
#2
\multirow{9}{#1}{1.9} & 0.1658 & 40.6 & 35.029 & 20.377 & 20.123 & 14.0 \\
& 0.2606 & 68.0 & 35.020 & 31.453 & 31.862 & 21.1 \\
& 0.3620 & 77.3 & 28.170 & 34.854 & 35.005 & 28.2 \\
& 0.4707 & 78.7 & 24.955 & 34.845 & 34.998 & 35.3 \\
& 0.5876 & 80.1 & 22.804 & 34.852 & 35.007 & 42.5 \\
& 0.7139 & 81.3 & 21.996 & 34.824 & 34.981 & 49.6 \\
& 0.8519 & 82.6 & 20.808 & 34.829 & 34.987 & 56.7 \\
& 1.0044 & 84.0 & 20.175 & 34.857 & 35.018 & 63.8 \\
& 1.1766 & 85.3 & 19.910 & 34.818 & 34.980 & 70.9 \\
}
\newcommand{\tabledataC}[2]{
\multirow{9}{#1}{1} & 0.0782 & 37.4 & 34.989 & 29.773 & 20.671 & 14.0 \\
& 0.1229 & 46.4 & 34.202 & 35.020 & 24.658 & 21.1 \\
& 0.1706 & 48.2 & 30.139 & 35.007 & 24.994 & 28.2 \\
& 0.2217 & 49.8 & 27.744 & 34.988 & 25.447 & 35.3 \\
& 0.2764 & 51.4 & 27.323 & 34.997 & 25.982 & 42.5 \\
& 0.3354 & 53.0 & 26.756 & 34.992 & 26.566 & 49.6 \\
& 0.3996 & 54.7 & 26.524 & 35.017 & 27.239 & 56.7 \\
& 0.4701 & 56.5 & 25.978 & 35.019 & 28.007 & 63.8 \\
& 0.5489 & 58.3 & 26.011 & 34.995 & 28.818 & 70.9 \\
#2
\multirow{9}{#1}{1.3} & 0.1073 & 36.3 & 35.003 & 28.704 & 19.922 & 14.0 \\
& 0.1687 & 46.8 & 34.075 & 34.989 & 24.812 & 21.1 \\
& 0.2345 & 48.8 & 28.682 & 35.009 & 25.376 & 28.2 \\
& 0.3049 & 50.4 & 27.490 & 34.975 & 25.829 & 35.3 \\
& 0.3807 & 52.1 & 27.267 & 34.980 & 26.366 & 42.5 \\
& 0.4625 & 53.8 & 27.648 & 35.020 & 26.998 & 49.6 \\
& 0.5519 & 55.5 & 27.613 & 34.976 & 27.671 & 56.7 \\
& 0.6507 & 57.3 & 27.629 & 34.983 & 28.439 & 63.8 \\
& 0.7620 & 59.2 & 27.069 & 34.991 & 29.312 & 70.9 \\
}
\newcommand{\tabledataD}[2]{
\multirow{9}{#1}{1.6} & 0.1250 & 38.0 & 34.994 & 29.850 & 20.697 & 14.0 \\
& 0.1961 & 47.0 & 31.606 & 34.973 & 24.617 & 21.1 \\
& 0.2719 & 48.9 & 28.213 & 34.971 & 25.040 & 28.2 \\
& 0.3526 & 50.5 & 28.812 & 35.001 & 25.507 & 35.3 \\
& 0.4389 & 52.1 & 28.381 & 35.010 & 26.051 & 42.5 \\
& 0.5314 & 53.7 & 28.769 & 35.026 & 26.667 & 49.6 \\
& 0.6313 & 55.3 & 28.561 & 34.987 & 27.337 & 56.7 \\
& 0.7402 & 57.0 & 28.503 & 35.002 & 28.116 & 63.8 \\
& 0.8606 & 58.7 & 27.904 & 34.986 & 28.914 & 70.9 \\
#2
\multirow{9}{#1}{1.9} & 0.1658 & 39.2 & 34.987 & 31.114 & 22.158 & 14.0 \\
& 0.2606 & 46.7 & 31.153 & 35.006 & 25.318 & 21.1 \\
& 0.3620 & 48.8 & 29.873 & 35.024 & 25.811 & 28.2 \\
& 0.4707 & 50.5 & 30.215 & 34.994 & 26.329 & 35.3 \\
& 0.5876 & 52.4 & 30.021 & 35.012 & 27.087 & 42.5 \\
& 0.7139 & 54.2 & 30.163 & 34.993 & 27.873 & 49.6 \\
& 0.8519 & 56.1 & 29.608 & 35.001 & 28.712 & 56.7 \\
& 1.0044 & 57.9 & 29.404 & 34.971 & 29.518 & 63.8 \\
& 1.1766 & 59.8 & 28.884 & 35.019 & 30.347 & 70.9 \\
}
\newcommand{\tabledataE}[2]{
\multirow{9}{#1}{1} & 0.0782 & 79.4 & 35.005 & 22.468 & 16.464 & 14.0 \\
& 0.1229 & 84.7 & 35.016 & 24.159 & 21.610 & 21.1 \\
& 0.1706 & 88.2 & 34.983 & 26.829 & 26.639 & 28.2 \\
& 0.2217 & 95.3 & 35.002 & 32.742 & 32.928 & 35.3 \\
& 0.2764 & 83.9 & 34.983 & 32.088 & 32.250 & 42.5 \\
& 0.3354 & 82.3 & 30.545 & 34.498 & 34.986 & 49.6 \\
& 0.3996 & 75.3 & 25.763 & 34.058 & 34.984 & 56.7 \\
& 0.4701 & 69.8 & 20.034 & 33.667 & 34.988 & 63.8 \\
& 0.5489 & 65.8 & 17.993 & 33.375 & 34.984 & 70.9 \\
#2
\multirow{9}{#1}{1.3} & 0.1073 & 84.4 & 35.015 & 24.302 & 17.321 & 14.0 \\
& 0.1687 & 90.7 & 34.990 & 26.257 & 23.077 & 21.1 \\
& 0.2345 & 98.0 & 34.996 & 29.984 & 29.581 & 28.2 \\
& 0.3049 & 101.4 & 33.922 & 34.786 & 34.991 & 35.3 \\
& 0.3807 & 91.2 & 32.463 & 34.827 & 35.011 & 42.5 \\
& 0.4625 & 82.9 & 25.273 & 34.721 & 34.982 & 49.6 \\
& 0.5519 & 75.8 & 21.217 & 34.256 & 34.991 & 56.7 \\
& 0.6507 & 70.1 & 18.398 & 33.782 & 35.021 & 63.8 \\
& 0.7620 & 65.8 & 17.752 & 33.320 & 34.977 & 70.9 \\
}
\newcommand{\tabledataF}[2]{
\multirow{9}{#1}{1.6} & 0.1250 & 90.2 & 34.986 & 25.965 & 18.608 & 14.0 \\
& 0.1961 & 99.0 & 34.996 & 28.382 & 25.044 & 21.1 \\
& 0.2719 & 109.9 & 34.999 & 33.264 & 32.758 & 28.2 \\
& 0.3526 & 102.9 & 29.787 & 34.802 & 35.012 & 35.3 \\
& 0.4389 & 92.5 & 27.864 & 34.813 & 34.997 & 42.5 \\
& 0.5314 & 83.4 & 21.008 & 34.413 & 35.017 & 49.6 \\
& 0.6313 & 76.2 & 19.070 & 33.976 & 35.021 & 56.7 \\
& 0.7402 & 70.4 & 18.741 & 33.557 & 35.020 & 63.8 \\
& 0.8606 & 66.0 & 17.124 & 33.167 & 35.010 & 70.9 \\
#2
\multirow{9}{#1}{1.9} & 0.1658 & 90.1 & 34.987 & 26.120 & 19.663 & 14.0 \\
& 0.2606 & 110.2 & 35.010 & 31.893 & 29.362 & 21.1 \\
& 0.3620 & 112.2 & 32.317 & 34.751 & 34.990 & 28.2 \\
& 0.4707 & 98.8 & 28.715 & 34.784 & 34.987 & 35.3 \\
& 0.5876 & 89.3 & 23.257 & 34.821 & 35.000 & 42.5 \\
& 0.7139 & 81.8 & 22.031 & 34.836 & 34.994 & 49.6 \\
& 0.8519 & 76.1 & 20.271 & 34.878 & 35.020 & 56.7 \\
& 1.0044 & 70.4 & 18.713 & 34.312 & 35.022 & 63.8 \\
& 1.1766 & 65.9 & 17.210 & 33.710 & 35.002 & 70.9 \\
}
\newcommand{\tabledataG}[2]{
\multirow{9}{#1}{1} & 0.0782 & 80.0 & 34.982 & 29.713 & 21.052 & 14.0 \\
& 0.1229 & 79.6 & 35.015 & 33.617 & 24.854 & 21.1 \\
& 0.1706 & 72.2 & 34.984 & 34.880 & 26.143 & 28.2 \\
& 0.2217 & 64.4 & 31.879 & 35.016 & 26.474 & 35.3 \\
& 0.2764 & 58.4 & 27.509 & 35.008 & 26.611 & 42.5 \\
& 0.3354 & 53.4 & 26.782 & 35.015 & 26.621 & 49.6 \\
& 0.3996 & 49.5 & 26.341 & 35.017 & 26.630 & 56.7 \\
& 0.4701 & 46.4 & 25.805 & 34.991 & 26.679 & 63.8 \\
& 0.5489 & 44.4 & 18.935 & 35.034 & 26.993 & 70.9 \\
#2
\multirow{9}{#1}{1.3} & 0.1073 & 83.6 & 34.998 & 31.264 & 21.663 & 14.0 \\
& 0.1687 & 82.8 & 34.258 & 34.984 & 25.636 & 21.1 \\
& 0.2345 & 73.0 & 33.727 & 35.013 & 26.509 & 28.2 \\
& 0.3049 & 65.2 & 32.546 & 35.012 & 26.891 & 35.3 \\
& 0.3807 & 59.3 & 27.783 & 34.990 & 27.072 & 42.5 \\
& 0.4625 & 54.1 & 27.643 & 34.971 & 27.001 & 49.6 \\
& 0.5519 & 50.2 & 27.018 & 35.023 & 27.044 & 56.7 \\
& 0.6507 & 46.9 & 26.713 & 34.967 & 27.006 & 63.8 \\
& 0.7620 & 44.7 & 19.118 & 35.003 & 27.211 & 70.9 \\
}
\newcommand{\tabledataH}[2]{
\multirow{9}{#1}{1.6} & 0.1250 & 90.4 & 35.010 & 33.755 & 23.692 & 14.0 \\
& 0.1961 & 79.6 & 35.014 & 33.405 & 24.770 & 21.1 \\
& 0.2719 & 72.3 & 34.981 & 34.468 & 26.073 & 28.2 \\
& 0.3526 & 65.4 & 34.146 & 35.002 & 26.668 & 35.3 \\
& 0.4389 & 59.2 & 28.995 & 35.013 & 26.740 & 42.5 \\
& 0.5314 & 54.0 & 28.775 & 34.986 & 26.684 & 49.6 \\
& 0.6313 & 50.0 & 27.792 & 34.983 & 26.716 & 56.7 \\
& 0.7402 & 46.8 & 27.202 & 35.007 & 26.790 & 63.8 \\
& 0.8606 & 44.5 & 19.587 & 35.028 & 26.966 & 70.9 \\
#2
\multirow{9}{#1}{1.9} & 0.1658 & 85.6 & 35.018 & 31.877 & 23.546 & 14.0 \\
& 0.2606 & 74.0 & 35.017 & 30.942 & 24.005 & 21.1 \\
& 0.3620 & 68.0 & 35.020 & 32.483 & 25.386 & 28.2 \\
& 0.4707 & 64.3 & 34.993 & 34.444 & 27.059 & 35.3 \\
& 0.5876 & 59.6 & 30.755 & 35.020 & 27.821 & 42.5 \\
& 0.7139 & 54.6 & 30.196 & 35.007 & 27.927 & 49.6 \\
& 0.8519 & 50.6 & 28.893 & 34.971 & 27.957 & 56.7 \\
& 1.0044 & 47.3 & 27.861 & 34.966 & 27.912 & 63.8 \\
& 1.1766 & 44.9 & 20.053 & 35.033 & 27.955 & 70.9 \\
}

\AtBeginDocument{
  \setlength{\intextsep}{0.42\baselineskip}
}

\title{\iflatexml MATHEMATICAL MODELING OF THE MECHANICAL BEHAVIOR OF THREE-LAYER PLATES WITH A TETRACHIRAL HONEYCOMB CORE\else MATHEMATICAL MODELING OF THE MECHANICAL BEHAVIOR\\OF THREE-LAYER PLATES WITH A TETRACHIRAL\\HONEYCOMB CORE\fi}
\iflatexml
\author{Alexey V. Mazaev}
\email{alexeymazaev@outlook.com}
\else
\author{Alexey V. Mazaev\\[5pt]{\fontsize{13}{19.5}\selectfont\href{mailto:alexeymazaev@outlook.com}{alexeymazaev@outlook.com}}}
\fi
\date{}

\makeatletter
\renewcommand{\@maketitle}{
  \begin{center}
    \normalfont\normalsize\@title\par
    \vspace{\baselineskip}
    \normalfont\normalsize\@author\par
  \end{center}
  \vspace{0.45\baselineskip}
}
\iflatexml
\else
\renewcommand{\tableofcontents}{
  \begin{center}\normalfont\normalsize СОДЕРЖАНИЕ\end{center}
  \vspace{0.5\baselineskip}
  {\sloppy\@starttoc{toc}}
}
\renewcommand*{\l@section}{\@dottedtocline{1}{0em}{0em}}
\renewcommand*{\l@subsection}[2]{\@dottedtocline{2}{0.5cm}{0em}{\hspace*{0.5cm}#1}{#2}}
\fi
\makeatother

\begin{document}
\selectlanguage{english}
\maketitle
\thispagestyle{empty}
\begingroup
\let\small\normalsize
\renewenvironment{center}{}{}
\renewenvironment{quotation}{}{}
\renewcommand{\abstractname}{\iflatexml Abstract\else\vspace{0.5em}\fi{}}
\begin{abstract}
\iflatexml\else\hspace*{1.25cm}Abstract. \fi This work examines the mechanical behavior of three-layer plates with a tetrachiral honeycomb core and solid face layers under static bending conditions. The influence of discretization, relative density, and thickness of the honeycomb core on the stress state of the composite plates is investigated under two boundary conditions: rigid clamping and supports allowing elastic rotation. The first numerical experiment setup involves a constant thickness of each composite layer while varying the core relative density. The second experiment setup maintains a constant volume of the honeycomb core solid body, which causes its thickness to change as the relative density varies. Mathematical modeling is performed using the finite element method within the framework of linear elasticity, employing both three-dimensional modeling in Comsol Multiphysics and custom algorithms for solving a plane problem to analyze the stress state of the tetrachiral honeycomb-based multilayer plates. The technical process of manufacturing the composites is described, followed by laboratory tests under three-point bending conditions. Next, the diagrams showing the dependence of maximum stresses in the composite plate layers on the relative density and thickness of the honeycomb core are discussed in the first and second setups of the numerical experiments, respectively. The results demonstrate good agreement between the numerical data from the three-dimensional and plane finite element models. Furthermore, the laboratory data from the three-point bending tests qualitatively align with the numerical analysis.
\end{abstract}
\endgroup

\clearpage
\selectlanguage{russian}
\renewcommand{\figurename}{Рисунок}
\renewcommand{\tablename}{Таблица}
\begin{center}
\iflatexml
МАТЕМАТИЧЕСКОЕ МОДЕЛИРОВАНИЕ МЕХАНИЧЕСКОГО ПОВЕДЕНИЯ ТРЕХСЛОЙНЫХ ПЛАСТИН С СОТОВЫМ ЗАПОЛНИТЕЛЕМ ТЕТРАКИРАЛЬНОГО ТИПА
\else
МАТЕМАТИЧЕСКОЕ МОДЕЛИРОВАНИЕ МЕХАНИЧЕСКОГО\\
ПОВЕДЕНИЯ ТРЕХСЛОЙНЫХ ПЛАСТИН С СОТОВЫМ\\
ЗАПОЛНИТЕЛЕМ ТЕТРАКИРАЛЬНОГО ТИПА
\fi

\vspace{\baselineskip}
Мазаев Алексей Вячеславович
\end{center}
\vspace{0.30\baselineskip}

Аннотация. Данная работа посвящена изучению механического поведения трехслойных пластин с тетракиральным сотовым заполнителем и сплошными внешними слоями в условиях статического изгиба. Исследуется влияние дискретизации, относительной плотности и толщины сотового заполнителя на напряженное состояние композитных пластин при жестком защемлении и опирании с упругим поворотом. Первая постановка численных экспериментов предусматривает постоянную толщину каждого слоя композитов с варьированием относительной плотности заполнителя, а вторая постановка отличается постоянным объемом твердого тела сотового заполнителя, что приводит к изменению его толщины при изменении относительной плотности. Для математического моделирования используется метод конечных элементов в рамках линейной теории упругости посредством трехмерного моделирования в системе Comsol Multiphysics и путем разработки алгоритмов решения задачи плоской деформации для анализа напряженного состояния многослойных пластин с тетракиральными сотами. Описывается технологический процесс изготовления исследуемых композитов с последующими лабораторными испытаниями в условиях трехточечного изгиба. Далее обсуждаются диаграммы зависимости максимальных напряжений в слоях композитных пластин от относительной плотности и толщины сотового заполнителя в первой и во второй постановках численных экспериментов соответственно. В результате показано хорошее согласие численных данных трехмерного и плоского конечно-элементного моделирования, а лабораторные результаты трехточечного изгиба качественно согласуются с результатами численного анализа.

\clearpage
\tableofcontents
\clearpage

\iflatexml\else\phantomsection\fi
\section{ВВЕДЕНИЕ}

Актуальность работы. При разработке материалов конструкционного назначения основное внимание уделяется выгодному сочетанию механических свойств и объемно-массовых характеристик. С развитием производства сплошные пластины все чаще заменяются современными композиционными материалами с улучшенными свойствами, одну из разновидностей которых составляют слоистые композиты с сотовым заполнителем. Наибольшее распространение получили трехслойные пластины со сплошными внешними слоями и сотовым заполнителем шестиугольного типа. Известно, что заполнитель обеспечивает совместную работу внешних слоев и оказывает ключевое влияние на механические свойства трехслойного композита. Наряду с развитием новых технологий, в том числе 3D-печати, набирают популярность соты с новыми геометриями, механические свойства которых позволяют получать слоистые композиты с улучшенными характеристиками. В данной работе исследуются трехслойные пластины со сплошными внешними слоями и сотовым заполнителем тетракирального типа. Рассматриваемая сотовая структура представляет собой связку цилиндров, которые упорядоченно расположены по схеме квадратной решетки, при этом каждый из цилиндров содержит четыре тангенциально прикрепленных ребра. Таким образом, приставка «тетра-» указывает, что цилиндры сот имеют четыре прикрепленных ребра, а название «киральные» означает, что геометрия элементарной ячейки не совмещается со своим зеркальным отражением, в отличие, например, от шестиугольных сот. В данной работе изучается влияние дискретизации (количества элементарных ячеек), относительной плотности (коэффициента заполнения) и толщины тетракирального сотового заполнителя на напряженное состояние трехслойных композитов при статическом изгибе для различных граничных условий. Математическое моделирование производится в рамках теории упругости методом конечных элементов путем трехмерного моделирования в системе Comsol Multiphysics, а также с помощью разработанных автором моделей для анализа напряженного состояния многослойных пластин с тетракиральными сотами посредством решения плоской задачи теории упругости.

Научная новизна работы состоит в следующем:

\iflatexml
\renewcommand{\labelenumi}{\arabic{enumi}.}
\begin{enumerate}
\fi
\fullwidthnumitem{1}{Предложены эффективные математические конечно-элементные модели слоистой пластины с тетракиральными сотовыми прослойками с использованием совместных и несовместных плоских конечных элементов для описания напряженно-деформированного состояния в условиях статического воздействия в рамках теории упругости;}
\fullwidthnumitem{2}{Разработаны алгоритмы численной реализации математического моделирования напряженно-деформированного состояния многослойных пластин с тетракиральными сотами в условиях статического изгиба посредством решения плоской задачи теории упругости методом конечных элементов с применением совместных и несовместных элементов;}
\fullwidthnumitem{3}{Выполнено численное моделирование напряженного состояния трехслойных пластин на основе первой постановки численных экспериментов при постоянной толщине слоев и варьируемой относительной плотности заполнителя методом конечных элементов посредством алгоритмов решения плоской задачи и путем трехмерного моделирования в Comsol Multiphysics с последующим сопоставлением численных результатов с данными лабораторных испытаний;}
\fullwidthnumitem{4}{Выполнено численное моделирование напряженного состояния трехслойных пластин на основе второй постановки численных экспериментов при постоянном объеме твердого тела сот и варьируемой толщине заполнителя методом конечных элементов посредством алгоритмов решения плоской задачи и путем трехмерного моделирования в Comsol Multiphysics;}
\fullwidthnumitem{5}{Исследовано влияние дискретизации, относительной плотности и толщины тетракирального сотового заполнителя на напряженное состояние трехслойных пластин в условиях статического изгиба в обеих постановках численных экспериментов при жестком защемлении и опирании с упругим поворотом.}
\iflatexml
\end{enumerate}
\fi

Достоверность результатов работы основывается на корректности предложенных математических моделей в рамках плоской задачи, хорошем согласии численных результатов, полученных посредством этих моделей, с результатами трехмерного конечно-элементного моделирования в системе Comsol Multiphysics, а также на качественном соответствии численных результатов данным лабораторных испытаний.

\clearpage
\iflatexml\else\phantomsection\fi
\section{ГЛАВА 1. ОБЗОР ЛИТЕРАТУРЫ ПО ТЕМЕ ИССЛЕДОВАНИЯ}

Композитные пластины с сотовым заполнителем давно зарекомендовали себя в качестве достойной альтернативы сплошным пластинам с точки зрения конструкционного применения \dualcite{ref9,ref11,ref22}{{[}\bibnum{ref9}, \bibnum{ref11}, \bibnum{ref22}{]}}{{[}\bibnum{ref9}, \bibnum{ref11}, \bibnum{ref22}{]}}. При этом ключевыми преимуществами слоистых композитов являются их высокая удельная жесткость и прочность. Наиболее популярны трехслойные композиты с двумя внешними слоями и сотовым заполнителем, обеспечивающим совместную работу несущих слоев \dualcite{ref9,ref11,ref22,ref43,ref55,ref56,ref57,ref60,ref64,ref66,ref67,ref70,ref76,ref77,ref78,ref79,ref88,ref89,ref90,ref91,ref92,ref93,ref94,ref99,ref100,ref104}{{[}\bibnum{ref9}, \bibnum{ref11}, \bibnum{ref22}, \bibnum{ref43}, \bibnum{ref55}--\bibnum{ref57}, \bibnum{ref60}, \bibnum{ref64}, \bibnum{ref66}, \bibnum{ref67}, \bibnum{ref70}, \bibnum{ref76}--\bibnum{ref79}, \bibnum{ref88}--\bibnum{ref94}, \bibnum{ref99}, \bibnum{ref100}, \bibnum{ref104}{]}}{{[}\bibnum{ref9}, \bibnum{ref11}, \bibnum{ref22}, \bibnum{ref43}, \bibnum{ref55}--\bibnum{ref57}, \bibnum{ref60}, \bibnum{ref64}, \bibnum{ref66}, \bibnum{ref67}, \bibnum{ref70}, \bibnum{ref76}--\bibnum{ref79}, \bibnum{ref88}--\bibnum{ref94}, \bibnum{ref99}, \bibnum{ref100}, \bibnum{ref104}{]}}. Как правило, характеристики сотового заполнителя, такие как геометрия элементарных ячеек, регулярность структуры, относительная плотность и высота заполнителя, определяют механические свойства композитной пластины. При этом наибольшее распространение в качестве заполнителя получили обычные шестиугольные соты \dualcite{ref55,ref60,ref64,ref66,ref67,ref79,ref90,ref91,ref92,ref93,ref94}{{[}\bibnum{ref55}, \bibnum{ref60}, \bibnum{ref64}, \bibnum{ref66}, \bibnum{ref67}, \bibnum{ref79}, \bibnum{ref90}--\bibnum{ref94}{]}}{{[}\bibnum{ref55}, \bibnum{ref60}, \bibnum{ref64}, \bibnum{ref66}, \bibnum{ref67}, \bibnum{ref79}, \bibnum{ref90}--\bibnum{ref94}{]}}, проявляющие положительный коэффициент Пуассона в плоскости \dualcite{ref9,ref11,ref58}{{[}\bibnum{ref9}, \bibnum{ref11}, \bibnum{ref58}{]}}{{[}\bibnum{ref9}, \bibnum{ref11}, \bibnum{ref58}{]}}. Однако наряду с развитием новых технологий, в том числе аддитивного производства \dualcite{ref63,ref83,ref98,ref102}{{[}\bibnum{ref63}, \bibnum{ref83}, \bibnum{ref98}, \bibnum{ref102}{]}}{{[}\bibnum{ref63}, \bibnum{ref83}, \bibnum{ref98}, \bibnum{ref102}{]}}, набирают популярность новые геометрии сот \dualcite{ref86,ref105}{{[}\bibnum{ref86}, \bibnum{ref105}{]}}{{[}\bibnum{ref86}, \bibnum{ref105}{]}}, среди которых немалую часть составляют соты с отрицательным коэффициентом Пуассона (ауксетические соты) \dualcite{ref15,ref39,ref40,ref41,ref51,ref62,ref69,ref75,ref80,ref82,ref95,ref105,ref106}{{[}\bibnum{ref15}, \bibnum{ref39}--\bibnum{ref41}, \bibnum{ref51}, \bibnum{ref62}, \bibnum{ref69}, \bibnum{ref75}, \bibnum{ref80}, \bibnum{ref82}, \bibnum{ref95}, \bibnum{ref105}, \bibnum{ref106}{]}}{{[}\bibnum{ref15}, \bibnum{ref39}--\bibnum{ref41}, \bibnum{ref51}, \bibnum{ref62}, \bibnum{ref69}, \bibnum{ref75}, \bibnum{ref80}, \bibnum{ref82}, \bibnum{ref95}, \bibnum{ref105}, \bibnum{ref106}{]}}.

\iflatexml\else\phantomsection\fi
\subsection{1.1. Трехслойные композиты с сотовым заполнителем в условиях изгиба}

Scarpa и Tomlinson \dualcite{ref90,ref91}{{[}\bibnum{ref90}, \bibnum{ref91}{]}}{{[}\bibnum{ref90}, \bibnum{ref91}{]}} рассмотрели повторно-входящие и шестиугольные соты с использованием теории ячеистого материала. Авторы также применили теорию многослойных пластин первого порядка для получения основных частот слоистых пластин с вышеуказанными сотами для случаев свободного опирания и цилиндрического изгиба. Авторы показали, что повторно-входящие соты, обладающие отрицательным коэффициентом Пуассона (ОКП), при определенных значениях геометрических параметров позволяют получить более высокие значения модуля сдвига из плоскости по сравнению с обычными шестиугольными сотами. В свою очередь это приводит к увеличению изгибной жесткости слоистых пластин. Lira, Scarpa и Rajasekaran \dualcite{ref71}{{[}\bibnum{ref71}{]}}{{[}\bibnum{ref71}{]}} посредством численного моделирования и экспериментов также показали, что ауксетическая сотовая структура повторно-входящего типа имеет повышенную удельную изгибную жесткость относительно шестиугольных сот. Аналогичный вывод следует из результатов лабораторных экспериментов Harland, Alshaer и Brooks с полимерным заполнителем \dualcite{ref60}{{[}\bibnum{ref60}{]}}{{[}\bibnum{ref60}{]}}. При этом ауксетические соты позволяют при сниженной массе получить одинаковые первые собственные частоты по сравнению с шестиугольными сотами \dualcite{ref71}{{[}\bibnum{ref71}{]}}{{[}\bibnum{ref71}{]}}.

Crupi, Epasto и Guglielmino \dualcite{ref55}{{[}\bibnum{ref55}{]}}{{[}\bibnum{ref55}{]}} теоретически и экспериментально исследовали влияние размера элементарных ячеек шестиугольного сотового заполнителя на механическое поведение сэндвич-композитов при трехточечном изгибе и ударной нагрузке. По результатам статического изгиба авторы выявили два режима разрушения, а после ударных испытаний с помощью томографического анализа показали, что разрушение сотового заполнителя происходит из-за смятия стенок сот вследствие потери устойчивости.

Sun и др. \dualcite{ref94}{{[}\bibnum{ref94}{]}}{{[}\bibnum{ref94}{]}} с помощью лабораторных экспериментов и конечно-элементного моделирования исследовали поведение трехслойных сэндвич-композитов при трехточечном изгибе и сжатии в плоскости основания. При этом изучалось влияние толщины внешних листов, размера шестиугольных ячеек и толщины их стенок, а также высоты сотового заполнителя на прочность, жесткость и способность композитов поглощать энергию. В результате авторы предложили аналитические выражения для определения пиковой нагрузки, поглощенной энергии и режима разрушения при трехточечном изгибе. Схожее исследование композитов с шестиугольным сотовым заполнителем с помощью конечно-элементного анализа провели Karakaya и Eksi \dualcite{ref64}{{[}\bibnum{ref64}{]}}{{[}\bibnum{ref64}{]}} в условиях трехточечного изгиба, но без анализа влияния толщины внешних листов.

Meran и Çetin \dualcite{ref79}{{[}\bibnum{ref79}{]}}{{[}\bibnum{ref79}{]}} методом конечных элементов исследовали напряженное состояние трехслойных сэндвич-композитов с шестиугольным сотовым заполнителем в условиях жесткого защемления с двух сторон при статическом изгибе, а также определили собственные частоты колебаний. В работе изучалось влияние размера элементарных ячеек сот на максимальные эквивалентные напряжения и прогиб композитных пластин при постоянной относительной плотности сотового заполнителя. На основе диаграмм распределения напряжений для большинства моделей можно сделать вывод, что увеличение размера элементарных ячеек сот может приводить к снижению максимальных эквивалентных напряжений в сотовом заполнителе. Kumar, Angra и Chanda \dualcite{ref67}{{[}\bibnum{ref67}{]}}{{[}\bibnum{ref67}{]}} с помощью конечно-элементного анализа также исследовали влияние размера элементарных ячеек шестиугольных сот на жесткость сэндвич-пластин в условиях изгиба.

В отличие от данного исследования и работы \dualcite{ref79}{{[}\bibnum{ref79}{]}}{{[}\bibnum{ref79}{]}}, в публикациях \dualcite{ref55,ref64,ref67,ref94}{{[}\bibnum{ref55}, \bibnum{ref64}, \bibnum{ref67}, \bibnum{ref94}{]}}{{[}\bibnum{ref55}, \bibnum{ref64}, \bibnum{ref67}, \bibnum{ref94}{]}} не рассматривается условие постоянного объема твердого тела сотовых прослоек. Таким образом, соты с разным размером элементарных ячеек имеют разную относительную плотность, что не позволяет оценить непосредственное влияние размера элементарных ячеек на механические свойства композитов.

Shah \dualcite{ref92}{{[}\bibnum{ref92}{]}}{{[}\bibnum{ref92}{]}} с помощью метода конечных элементов исследовал механические свойства ряда сотовых прослоек. При этом был проведен квазистатический анализ прочности трехслойных сэндвич-композитов с шестиугольными и треугольными сотами при трехточечном изгибе, в рамках которого рассмотрено влияние размера элементарных ячеек и относительной плотности сот на критическую нагрузку. Согласно результатам численных экспериментов, увеличение размера элементарных ячеек при одинаковой относительной плотности в основном приводило к снижению критической нагрузки. При этом использование шестиугольного сотового заполнителя привело к большей прочности относительно треугольного заполнителя.

Smardzewski и Prekrat \dualcite{ref89}{{[}\bibnum{ref89}{]}}{{[}\bibnum{ref89}{]}} численно исследовали влияние четырех геометрических типов сотовой прослойки на механические свойства трехслойных панелей при кручении и трехточечном изгибе. Использовались соты следующих типов: повторно-входящего, двойной стрелки, звездообразного и на базе цилиндров; при этом прослойки имели сопоставимую относительную плотность. По результатам конечно-элементного моделирования трехточечного изгиба авторы определили модули упругости и сдвига, а также прочность композитных пластин.

Araújo и др. \dualcite{ref42}{{[}\bibnum{ref42}{]}}{{[}\bibnum{ref42}{]}} с помощью конечно-элементного анализа и лабораторных экспериментов исследовали прочность и жесткость сотовых пластин в условиях трехточечного изгиба. Использовались следующие геометрические формы сот: обычная шестиугольная, типа лотоса и шестиугольная со скругленными углами. При этом сотовые пластины каждого типа моделировались с четырьмя значениями относительной плотности. Однако в рамках некоторых групп образцы значительно отличались по относительной плотности, а одна из групп имела меньшее количество элементарных ячеек, что затрудняет оценку непосредственного влияния типа ячейки на механические характеристики сот.

Essassi и др. \dualcite{ref56,ref57}{{[}\bibnum{ref56}, \bibnum{ref57}{]}}{{[}\bibnum{ref56}, \bibnum{ref57}{]}} экспериментально исследовали влияние относительной плотности на коэффициент Пуассона повторно-входящего сотового заполнителя, а также на механические свойства трехслойных сэндвич-композитов при статическом изгибе. В результате было показано, что при увеличении относительной плотности повторно-входящих сот увеличивается их коэффициент Пуассона, а также ожидаемо увеличивается изгибная и сдвиговая жесткость сэндвич-композита. Очевидно, что при увеличении относительной плотности за счет утолщения стенок элементарных ячеек сотовые структуры стремятся к сплошной среде. При этом если в качестве твердого тела используется классический материал, коэффициент Пуассона ауксетических сот будет переходить из области отрицательных значений в область положительных, как это графически показано в главе \hyperref[sec:chapter3]{3} данной работы.

В работах \dualcite{ref43,ref66,ref93}{{[}\bibnum{ref43}, \bibnum{ref66}, \bibnum{ref93}{]}}{{[}\bibnum{ref43}, \bibnum{ref66}, \bibnum{ref93}{]}} экспериментально исследовано влияние толщины сотового заполнителя на прочность, а также жесткость \dualcite{ref43,ref93}{{[}\bibnum{ref43}, \bibnum{ref93}{]}}{{[}\bibnum{ref43}, \bibnum{ref93}{]}} трехслойных пластин в условиях трехточечного изгиба. По результатам испытаний было показано, что увеличение толщины сотового заполнителя с пропорциональным увеличением объема твердого тела сот ожидаемо приводит к повышению жесткости и прочности композитов.

Brischetto и др. \dualcite{ref49}{{[}\bibnum{ref49}{]}}{{[}\bibnum{ref49}{]}} экспериментально исследовали трехслойные сэндвич-композиты с шестиугольным сотовым заполнителем при трехточечном изгибе. Образцы были изготовлены методом экструзионной 3D-печати из полимерных материалов. В результате были построены графики нагрузки-прогиба и напряжения-деформации для исследуемых образцов, а также определены модули упругости при изгибе.

Hou и др. \dualcite{ref61}{{[}\bibnum{ref61}{]}}{{[}\bibnum{ref61}{]}} исследовали в условиях трехточечного изгиба сэндвич-композиты с полиморфными сотовыми прослойками, геометрическая форма которых основана на повторно-входящем типе. Сотовые структуры изготавливались из полимера с помощью 3D-печати, а внешние пластины состояли из углепластика. Было показано, что градиентная геометрия сот, полученная изменением длин ребер ячеек и углов, образуемых ребрами, непосредственно влияет на локализацию начального разрушения сот.

Xiao и др. \dualcite{ref101}{{[}\bibnum{ref101}{]}}{{[}\bibnum{ref101}{]}} исследовали трехслойные сэндвич-композиты в условиях квазистатического изгиба. В качестве материала шестиугольных сот использовался алюминиевый сплав, а внешние листы состояли из углепластика. В результате были получены экспериментальные графики нагрузки-прогиба, исследованы режимы разрушения и способность к поглощению энергии композитов. Произведена оценка влияния толщины внешних листов, укладки в них углеволокна, а также скорости нагружения на прочность композитов. Полученные результаты прошли численную проверку с помощью метода конечных элементов.

Belouettar и др. \dualcite{ref48}{{[}\bibnum{ref48}{]}}{{[}\bibnum{ref48}{]}} исследовали алюминиевые и алюмополимерные сэндвич-композиты на статическую и усталостную прочность в условиях четырехточечного изгиба. Внешние слои выполнялись из алюминиевых листов, а шестиугольные сотовые прослойки изготавливались из алюминия и композита на основе арамидных волокон. Авторы сделали вывод о том, что при циклических испытаниях изгибная жесткость не является однозначным показателем хорошего состояния образца, поскольку активная область композита может накопить значительные повреждения с менее выраженной потерей жесткости.

Lister \dualcite{ref72}{{[}\bibnum{ref72}{]}}{{[}\bibnum{ref72}{]}} с помощью численного, аналитического и экспериментального подходов исследовал сэндвич-композиты с шестиугольным сотовым заполнителем в условиях трехточечного изгиба. В качестве материала внешних слоев использовался углепластик, а для сотовых прослоек -- арамидно-бумажный композит. Было показано, что ориентация ленты шестиугольного сотового заполнителя непосредственно влияет на механические свойства: композиты имеют наибольшую прочность и жесткость при продольной ориентации ленты, а при поперечной ориентации, наоборот, наименьшую. Снижение прочности вследствие изменения ориентации ленты сот от продольной к поперечной также вызвано увеличением антикластической кривизны, которая приводит к уменьшению площади контакта при изгибе и, следовательно, увеличению напряжений. Показано, что использование различной толщины внешних листов может служить дополнительным фактором увеличения жесткости сэндвич-композитов по сравнению с симметричными композитами при одинаковой массе.

\iflatexml\else\phantomsection\fi
\subsection{1.2. Особенности механического поведения тетракиральных сот в плоскости}

Alderson и др. \dualcite{ref39}{{[}\bibnum{ref39}{]}}{{[}\bibnum{ref39}{]}} определили значения модуля Юнга и коэффициента Пуассона у ряда сотовых структур кирального типа. В качестве методов использовали численный анализ с помощью конечно-элементного моделирования и лабораторные эксперименты на образцах из нейлона. В работе показано, что семейство киральных сот может демонстрировать ауксетическое поведение.

Chen и др. \dualcite{ref52,ref53}{{[}\bibnum{ref52}, \bibnum{ref53}{]}}{{[}\bibnum{ref52}, \bibnum{ref53}{]}} описали механическое поведение тетракиральной решетки с помощью плоской ортотропной микрополярной теории. В рамках предположения о жесткости кругов прямоугольной решетки были выведены аналитические выражения для 13 приведенных констант микрополярного материала. Также была предложена процедура численного усреднения, которая учитывает деформируемость кругов решетки. Показано, что деформация круга влияет на киральное, механическое и ауксетическое поведение решетки. Модуль Юнга и коэффициент Пуассона тетракиральных структур зависят от ориентации нагрузки, вследствие чего тетракиральные соты проявляют ауксетичность лишь в узком диапазоне направлений, в отличие от трикиральных структур.

Bacigalupo и Gambarotta \dualcite{ref45}{{[}\bibnum{ref45}{]}}{{[}\bibnum{ref45}{]}} рассмотрели нелокальное усреднение гексакиральной и тетракиральной сотовой структуры с помощью двух подходов. В первом случае структура смоделирована как балка-решетка с последующим усреднением в виде микрополярного континуума. Во втором подходе, предложенном авторами, рассмотрены периодические ячейки из деформируемых частей посредством градиентного усреднения. Были выведены аналитические выражения для приведенных констант сотовых структур с помощью микрополярного усреднения, а также показана зависимость модуля Юнга и коэффициента Пуассона тетракиральных структур от ориентации одноосного напряжения.

Bacigalupo, Badino и Gambarotta \dualcite{ref46}{{[}\bibnum{ref46}{]}}{{[}\bibnum{ref46}{]}} развили подход усреднения для многослойных балочно-решетчатых метаматериалов с чередующейся киральной микроструктурой. Авторами предложен метод динамического перехода от дискретной модели к континуальной, позволяющий описывать дисперсионные свойства многослойной решетки с учетом взаимного соединения отдельных слоев. Данная работа показывает, что континуальные модели киральных структур могут быть расширены за пределы квазистатического описания и использованы для анализа волновых эффектов, обусловленных микроструктурой.

Изотропные киральные микрополярные среды также называют гемитропными \dualcite{ref13,ref14,ref24,ref52,ref53}{{[}\bibnum{ref13}, \bibnum{ref14}, \bibnum{ref24}, \bibnum{ref52}, \bibnum{ref53}{]}}{{[}\bibnum{ref13}, \bibnum{ref14}, \bibnum{ref24}, \bibnum{ref52}, \bibnum{ref53}{]}}, а в ряде работ -- полуизотропными \dualcite{ref13,ref14,ref24}{{[}\bibnum{ref13}, \bibnum{ref14}, \bibnum{ref24}{]}}{{[}\bibnum{ref13}, \bibnum{ref14}, \bibnum{ref24}{]}}. Для описания тетракиральных решеток при этом используется ортотропная киральная микрополярная модель \dualcite{ref52,ref53}{{[}\bibnum{ref52}, \bibnum{ref53}{]}}{{[}\bibnum{ref52}, \bibnum{ref53}{]}}.

Согласно результатам работ \dualcite{ref45,ref52,ref53}{{[}\bibnum{ref45}, \bibnum{ref52}, \bibnum{ref53}{]}}{{[}\bibnum{ref45}, \bibnum{ref52}, \bibnum{ref53}{]}}, балочно-решетчатые модели тетракиральных сот с жесткими круглыми узлами демонстрируют нулевой коэффициент Пуассона при одноосной деформации вдоль основных направлений в плоскости, что не согласуется с экспериментальными результатами работы \dualcite{ref39}{{[}\bibnum{ref39}{]}}{{[}\bibnum{ref39}{]}}. При этом учет деформируемости кругов и других особенностей микроструктуры может приводить к ненулевым значениям.

Mousanezhad и др. \dualcite{ref80}{{[}\bibnum{ref80}{]}}{{[}\bibnum{ref80}{]}} провели теоретическое и численное исследование механического поведения сотовых структур кирального и иерархического типа, геометрические особенности которых встречаются в природе. Для изучения упругих констант использовался энергетический подход. В результате были выведены аналитические выражения для приведенного модуля Юнга, модуля сдвига и коэффициента Пуассона киральных, антикиральных и иерархических сотовых структур на базе квадрата и шестиугольника. Аналитические результаты прошли численную проверку с помощью метода конечных элементов. Показано, что иерархия и киральность оказывают значительное влияние на механические свойства структур.

Zhong и др. \dualcite{ref107}{{[}\bibnum{ref107}{]}}{{[}\bibnum{ref107}{]}} вывели аналитические выражения для приведенного модуля Юнга, модуля сдвига и коэффициента Пуассона тетракиральных сот. В работе использовался метод эллиптических интегралов в рамках теории больших прогибов тонких балок Эйлера-Бернулли. Теоретические результаты были подтверждены с помощью численных расчетов. При моделировании нелинейного механического поведения тетракиральных сот было выявлено, что коэффициент Пуассона имеет положительные значения при растяжении и отрицательные значения при сжатии. Такое поведение авторы объяснили влиянием взаимосвязи между сдвиговой и нормальной деформациями при одноосном нагружении. Помимо этого, при нелинейном моделировании коэффициент Пуассона проявил зависимость от величины деформации, а при линейном оказался равен нулю, что согласуется с аналитическим результатом работы \dualcite{ref80}{{[}\bibnum{ref80}{]}}{{[}\bibnum{ref80}{]}}.

Qi и др. \dualcite{ref85}{{[}\bibnum{ref85}{]}}{{[}\bibnum{ref85}{]}} исследовали характер деформирования тетракиральных ауксетических сот при квазистатическом и динамическом одноосном сжатии с помощью численных и теоретических методов, а также лабораторных экспериментов. При квазистатическом продольном сжатии тетракиральные соты демонстрировали неравномерное поперечное сжатие в режиме «Z», вследствие которого в середине экспериментального образца появлялась наклонная полоса деформированных элементарных ячеек, разделяющая образец на две части, а по краям наклонной полосы образовывались выпуклости с меньшей относительной плотностью. При динамическом продольном сжатии тетракиральные соты демонстрировали последовательное зигзагообразное уплотнение образца в режиме «I» с равномерным поперечным сжатием. Отсутствие выпуклостей при динамическом деформировании авторы объяснили запаздыванием волны напряжений в сотовой структуре по сравнению с ударной волной. При увеличении коэффициента \(\alpha_{h}\), представляющего собой отношение среднего радиуса цилиндров и длины тангенциально прикрепленных ребер, монотонно увеличивалась прочность при квазистатическом и динамическом сжатии, а также уменьшался эффект выпучивания. Однако у образцов с одинаковой относительной плотностью коэффициент \(\alpha_{h}\) практически не влиял на прочность при динамическом сжатии. Были предложены аналитические выражения для определения напряжения плато в двух постановках экспериментов. Помимо этого, авторы определили коэффициент Пуассона в плоскости у тетракиральных сот при квазистатическом и динамическом сжатии с использованием метода конечных элементов. ОКП определялся по перемещениям двух репрезентативных узлов без учета образующихся выпуклостей. Моделирование в условиях квазистатического сжатия при малых деформациях показало зависимость коэффициента Пуассона от величины деформации. При динамическом нагружении увеличение скорости удара повышало значение коэффициента Пуассона тетракиральных сот, а при уменьшении \(\alpha_{h}\) коэффициент монотонно уменьшался. Также получено аналитическое выражение для определения минимального значения приведенного коэффициента Пуассона при предварительной деформации уплотнения в условиях квазистатического сжатия.

Alomarah и др. \dualcite{ref41}{{[}\bibnum{ref41}{]}}{{[}\bibnum{ref41}{]}} исследовали при квазистатическом одноосном сжатии известные сотовые структуры: повторно-входящую, тетракиральную и антитетракиральную. В результате были определены значения ОКП, характер деформирования и эффективность поглощения энергии сот, а также изучена новая повторно-входящая киральная структура. В ходе исследования использовалось численное моделирование с помощью метода конечных элементов, а также проводились лабораторные испытания образцов из полиамида. При квазистатическом сжатии тетракиральные соты демонстрировали последовательное зигзагообразное уплотнение без выраженного режима деформации «Z» \dualcite{ref85}{{[}\bibnum{ref85}{]}}{{[}\bibnum{ref85}{]}}, а скорее в режиме деформации «I» с минимальным эффектом выпуклости. Учитывая малое значение коэффициента \(\alpha_{h}\) у образцов данной работы, полученные результаты слабо согласуются с результатами работы \dualcite{ref85}{{[}\bibnum{ref85}{]}}{{[}\bibnum{ref85}{]}}. В ходе численных и лабораторных экспериментов ОКП тетракиральных сот не продемонстрировал выраженной зависимости от величины деформации сжатия в сравнении с результатами работ \dualcite{ref85,ref107}{{[}\bibnum{ref85}, \bibnum{ref107}{]}}{{[}\bibnum{ref85}, \bibnum{ref107}{]}}. При вычислении коэффициента Пуассона вместо двух репрезентативных узлов \dualcite{ref85}{{[}\bibnum{ref85}{]}}{{[}\bibnum{ref85}{]}} авторы использовали десять репрезентативных точек в центральной части сот для определения средних перемещений в продольном и поперечном направлениях. В результате численных и лабораторных экспериментов тетракиральные соты продемонстрировали наибольшее значение модуля Юнга.

Schneider и др. \dualcite{ref95}{{[}\bibnum{ref95}{]}}{{[}\bibnum{ref95}{]}} с помощью лабораторных экспериментов и конечно-элементного моделирования исследовали жесткость, прочность, характер деформирования и поглощение энергии многослойных тетракиральных решеток. При одинаковой массе двухслойная конфигурация TC2 и четырехслойная TC4 продемонстрировали увеличение поглощенной энергии по сравнению с однослойной TC1. Авторы также показали, что последовательность и ориентация слоев определяют сохранение или утрату ауксетического поведения, тогда как переход от мезомасштабных к макромасштабным образцам преимущественно влияет на жесткость, прочность, характер разрушения и поглощение энергии. Наблюдаемый масштабный эффект авторы связывают с различиями в свойствах фотополимера, обусловленными особенностями процесса, включая различия в степени отверждения образцов разных размеров.

Lu, Tan и Tay \dualcite{ref74}{{[}\bibnum{ref74}{]}}{{[}\bibnum{ref74}{]}} численно и теоретически исследовали ауксетическое поведение тетракиральных сотовых структур. Наличие связи между сдвиговой и продольной деформациями в тетракиральных сотах при одноосном нагружении приводит к неоднозначности определения коэффициента Пуассона. Исходя из работы \dualcite{ref80}{{[}\bibnum{ref80}{]}}{{[}\bibnum{ref80}{]}}, известно, что традиционный способ вычисления коэффициента Пуассона с помощью матрицы податливости не показывает отрицательные значения. Однако в лабораторных экспериментах тетракиральные соты демонстрируют ауксетическое поведение \dualcite{ref39,ref41}{{[}\bibnum{ref39}, \bibnum{ref41}{]}}{{[}\bibnum{ref39}, \bibnum{ref41}{]}}. Lu, Tan и Tay использовали альтернативное описание ауксетичности на базе репрезентативного элемента объема с периодическими граничными условиями. Круги структуры рассматривались жесткими, а прикрепленные ребра деформировались в соответствии с теорией балки Эйлера. В результате были предложены численные подходы для определения приведенного коэффициента Пуассона либо с использованием матрицы жесткости, либо через поперечную и продольную деформацию при периодических граничных условиях. Авторами получено аналитическое выражение, связывающее геометрические параметры тетракиральных сот, для вычисления приведенного коэффициента Пуассона.

Lorato и др. \dualcite{ref73}{{[}\bibnum{ref73}{]}}{{[}\bibnum{ref73}{]}} описали механическое поведение трикиральных, тетракиральных и гексакиральных сотовых структур в направлении из плоскости. Были предложены аналитические выражения для определения поперечного модуля Юнга и границ Фойгта и Рейсса для жесткости при поперечном сдвиге. С помощью конечно-элементного анализа были верифицированы аналитические результаты и произведена оценка влияния толщины сот на жесткость при поперечном сдвиге. Теоретические модели прошли экспериментальную проверку на лабораторных образцах.

Рассмотренные выше работы демонстрируют, что аналитическое описание и численное определение коэффициента Пуассона тетракиральных сот не всегда показывают их ауксетичность, однако в ряде лабораторных экспериментов зарегистрирован ОКП тетракиральных сот. Исследование механического поведения киральных структур с помощью метода конечных элементов является удобным и проверенным подходом, однако данный метод чувствителен к граничным условиям, учету сдвиговой и нормальной деформаций, области усреднения и прочим факторам. В то же время ауксетичность может оказывать дополнительное влияние на жесткость и прочность сотовых структур. Опираясь на физические соображения и аналитическое выражение из работы \dualcite{ref85}{{[}\bibnum{ref85}{]}}{{[}\bibnum{ref85}{]}} для определения эффективного коэффициента Пуассона тетракиральных сот, можно заключить, что при увеличении относительной плотности сот коэффициент Пуассона переходит из области отрицательных значений в область положительных, если материал сот имеет классические свойства. Эта зависимость будет графически продемонстрирована далее в главе \hyperref[sec:chapter3]{3}.

\iflatexml\else\phantomsection\fi
\subsection{1.3. Матрица жесткости конечного элемента в рамках теории упругости в условиях плоской деформации}

Если в процессе деформации точки тела перемещаются параллельно одной координатной плоскости (например, \emph{xy}), то такое деформированное состояние называется плоской деформацией. В этом случае компоненты перемещения являются функциями координат \emph{x} и \emph{y}: \(u_{x} = u_{x}(x,y),\) \(u_{y} = u_{y}(x,y),\) \(u_{z} = 0\) \dualcite{ref1,ref6,ref7,ref19,ref20,ref21,ref25,ref28,ref29,ref30}{{[}\bibnum{ref1}, \bibnum{ref6}, \bibnum{ref7}, \bibnum{ref19}--\bibnum{ref21}, \bibnum{ref25}, \bibnum{ref28}--\bibnum{ref30}{]}}{{[}\bibnum{ref1}, \bibnum{ref6}, \bibnum{ref7}, \bibnum{ref19}--\bibnum{ref21}, \bibnum{ref25}, \bibnum{ref28}--\bibnum{ref30}{]}}. При этом относительная деформация в направлении из плоскости \(\varepsilon_{zz}\) и деформации сдвига \(\gamma_{xz}\) и \(\gamma_{yz}\) равны нулю. Касательные напряжения \(\tau_{xz}\) и \(\tau_{yz}\) также равны нулю. Таким образом, дифференциальные уравнения равновесия в частных производных по координатам \emph{x} и \emph{y} принимают вид:

\begin{equation}
\begin{matrix}
\frac{\partial\sigma_{xx}}{\partial x} + \frac{\partial\tau_{yx}}{\partial y} + R_{x} = 0, \\[0.55em]
\frac{\partial\tau_{xy}}{\partial x} + \frac{\partial\sigma_{yy}}{\partial y} + R_{y} = 0,
\end{matrix}
\end{equation}

\noindent где \(\sigma_{xx}\) и \(\sigma_{yy}\) -- нормальные напряжения вдоль координатных осей \emph{x} и \emph{y} соответственно, \(\tau_{xy} = \tau_{yx}\) -- касательные напряжения, \(R_{x}\) и \(R_{y}\) -- проекции объемных сил на координатные оси.

Для напряжений на наклонных площадках имеем:

\begin{equation}
\begin{matrix}
\sigma_{nx} = \sigma_{xx}n_{x} + \tau_{yx}n_{y}, \\
\sigma_{ny} = \tau_{xy}n_{x} + \sigma_{yy}n_{y},
\end{matrix}
\end{equation}

\noindent где \(n_{x}\) и \(n_{y}\) -- косинусы углов между внешней нормалью к наклонной площадке и координатными осями \emph{x} и \emph{y} соответственно, \(\sigma_{nx}\) и \(\sigma_{ny}\) -- проекции полного напряжения на наклонной площадке на координатные оси.

Статические граничные условия выражаются соотношениями:

\begin{equation}
\begin{matrix}
p_{x} = \sigma_{xx}n_{x} + \tau_{yx}n_{y}, \\
p_{y} = \tau_{xy}n_{x} + \sigma_{yy}n_{y},
\end{matrix}
\end{equation}

\noindent где \(p_{x}\) и \(p_{y}\) -- составляющие вектора поверхностной нагрузки.

Уравнения Коши имеют вид:

\begin{equation}
\begin{matrix}
\varepsilon_{xx} = \frac{\partial u_{x}}{\partial x}, \\[0.55em]
\varepsilon_{yy} = \frac{\partial u_{y}}{\partial y}, \\[0.55em]
\gamma_{xy} = \frac{\partial u_{y}}{\partial x} + \frac{\partial u_{x}}{\partial y},
\end{matrix}
\end{equation}

\noindent где \(\varepsilon_{xx}\) и \(\varepsilon_{yy}\) -- нормальные деформации вдоль координатных осей \emph{x} и \emph{y} соответственно, а \(\gamma_{xy}\) -- деформация сдвига.

Условие совместности деформаций определяется выражением:

\begin{equation}
\frac{\partial^{2}\varepsilon_{xx}}{\partial y^{2}} + \frac{\partial^{2}\varepsilon_{yy}}{\partial x^{2}} = \frac{\partial^{2}\gamma_{xy}}{\partial x\partial y}.
\end{equation}

Матричная запись закона Гука имеет вид:

\begin{equation}
\label{eq:6}
\sigma = \chi\varepsilon = \frac{E}{(1 + \mu)(1 - 2\mu)}\begin{pmatrix}
1 - \mu & \mu & 0 \\
\mu & 1 - \mu & 0 \\
0 & 0 & \frac{1 - 2\mu}{2}
\end{pmatrix}\left( \begin{aligned}
 & \varepsilon_{xx} \\
 & \varepsilon_{yy} \\
 & \gamma_{xy}
\end{aligned} \right),
\end{equation}

\noindent где \emph{E} -- модуль упругости, \emph{\(\mu\)} -- коэффициент Пуассона, \emph{\(\chi\)} -- матрица упругих постоянных для плоской деформации, \(\varepsilon\) -- вектор деформаций, \(\sigma\) -- вектор напряжений.

Рассмотрим упругое равновесие деформированного конечного элемента \dualcite{ref5,ref6,ref7,ref20,ref21,ref25,ref27,ref36}{{[}\bibnum{ref5}--\bibnum{ref7}, \bibnum{ref20}, \bibnum{ref21}, \bibnum{ref25}, \bibnum{ref27}, \bibnum{ref36}{]}}{{[}\bibnum{ref5}--\bibnum{ref7}, \bibnum{ref20}, \bibnum{ref21}, \bibnum{ref25}, \bibnum{ref27}, \bibnum{ref36}{]}}. Предположим, что его узловые перемещения \(v\) получили произвольные бесконечно малые приращения, определяемые вектором \(\delta v = \begin{pmatrix}
\delta v_{1}^{x} & \delta v_{1}^{y} & \delta v_{2}^{x} & \delta v_{2}^{y} & \ldots
\end{pmatrix}^{T},\) где \(v_{q}^{x}\) и \(v_{q}^{y}\) -- узловые перемещения по осям \emph{x} и \emph{y} соответственно, \emph{q} -- номер узла конечного элемента.

Перемещения \(\widetilde{u} = \begin{pmatrix}
{\widetilde{u}}_{x} & {\widetilde{u}}_{y}
\end{pmatrix}^{T}\) произвольной точки элемента связаны с узловыми перемещениями равенством \(\widetilde{u} = \alpha v,\) где \(\alpha = \begin{pmatrix}
\alpha_{x} & \alpha_{y}
\end{pmatrix}^{T},\) \(\alpha_{x}\) и \(\alpha_{y}\) -- матрицы-строки, элементами которых являются аппроксимирующие функции. Таким образом, компоненты вектора \(\widetilde{u}\) получат приращения \(\delta\widetilde{u} = \alpha\delta v,\) а согласно принципу возможных перемещений:

\begin{equation}
\delta U^{e} = \delta A^{e},
\end{equation}

\noindent где \(\delta A^{e}\) -- работа внешних сил на возможных перемещениях \(\delta\widetilde{u},\) \(\delta U^{e}\) -- вариация потенциальной энергии деформации конечного элемента.

Заменим напряжения на поверхности элемента эквивалентными им сосредоточенными узловыми силами \({P_{q}^{e}}^{x}\) и \({P_{q}^{e}}^{y},\) действующими в направлении узловых перемещений \(P^{e} = \begin{pmatrix}
{P_{1}^{e}}^{x} & {P_{1}^{e}}^{y} & {P_{2}^{e}}^{x} & {P_{2}^{e}}^{y} & \ldots
\end{pmatrix}^{T}.\) Объемные и поверхностные силы также заменим эквивалентными узловыми силами \({\widetilde{P}}^{e} = \begin{pmatrix}
\mbox{$\widetilde{P}_{1}^{e}$}^{x} & \mbox{$\widetilde{P}_{1}^{e}$}^{y} & \mbox{$\widetilde{P}_{2}^{e}$}^{x} & \mbox{$\widetilde{P}_{2}^{e}$}^{y} & \ldots
\end{pmatrix}^{T}.\)

Будем считать силы \(P^{e} + {\widetilde{P}}^{e}\) эквивалентными действительным нагрузкам при условии выполнения равенства:

\begin{equation}
\label{eq:8}
\delta A^{e} = (\delta v)^{T}(P^{e} + {\widetilde{P}}^{e}).
\end{equation}

Вектор деформаций \(\varepsilon\) связан с вектором узловых перемещений \(v\) соотношением \(\varepsilon = \beta v,\) где \(\beta\) -- матрица связи между узловыми перемещениями и деформациями:

\begin{equation}
\beta = L\alpha = \begin{pmatrix}
\frac{\partial}{\partial x} & 0 \\
0 & \frac{\partial}{\partial y} \\
\frac{\partial}{\partial y} & \frac{\partial}{\partial x}
\end{pmatrix}\begin{pmatrix}
\alpha_{x} \\
\alpha_{y}
\end{pmatrix},
\end{equation}

\noindent \emph{L} -- матричный дифференциальный оператор.

Таким образом, приращение узловых перемещений вызовет приращение деформаций:

\begin{equation}
\delta\varepsilon = \beta\delta v.
\end{equation}

Вариацию потенциальной энергии деформации можно определить интегрированием:

\begin{equation}
\label{eq:11}
\delta U^{e} = \int_{\tau^{e}}^{}{\delta\varepsilon^{T}\sigma d\tau},
\end{equation}

\noindent где \(\tau^{e}\) -- объем конечного элемента.

Произведя замены \(\delta\varepsilon = \beta\delta v\) и \(\sigma = \chi\beta v\) в выражении  (\ref{eq:11}), получим:

\begin{equation}
\label{eq:12}
\delta U^{e} = \int_{\tau^{e}}^{}{(\beta\delta v)^{T}\chi\beta vd\tau} = \int_{\tau^{e}}^{}{(\delta v)^{T}\beta^{T}\chi\beta vd\tau}.
\end{equation}

Векторы \(v\) и \(\delta v\) в  (\ref{eq:12}) можно вынести за знак интеграла в силу постоянства их элементов, таким образом:

\begin{equation}
\delta U^{e} = (\delta v)^{T}\left\lbrack \int_{\tau^{e}}^{}{\beta^{T}\chi\beta d\tau} \right\rbrack v.
\end{equation}

Введем обозначение:

\begin{equation}
k^{e} = \int_{\tau^{e}}^{}{\beta^{T}\chi\beta d\tau}
\end{equation}

\noindent и представим выражение для \(\delta U^{e}\) в виде:

\begin{equation}
\label{eq:15}
\delta U^{e} = (\delta v)^{T}k^{e}v.
\end{equation}

Приравнивая  (\ref{eq:8}) к  (\ref{eq:15}) в соответствии с принципом возможных перемещений, получим соотношение, которое связывает узловые силы с узловыми перемещениями в стандартной форме, а матрица \(k^{e}\) является матрицей жесткости конечного элемента:

\begin{equation}
k^{e}v = P^{e} + {\widetilde{P}}^{e}.
\end{equation}

\clearpage
\iflatexml
\section[ГЛАВА 2. МАТЕМАТИЧЕСКОЕ МОДЕЛИРОВАНИЕ МЕХАНИЧЕСКОГО ПОВЕДЕНИЯ СЛОИСТЫХ ПЛАСТИН С ТЕТРАКИРАЛЬНЫМИ СОТАМИ МЕТОДОМ КОНЕЧНЫХ ЭЛЕМЕНТОВ]{ГЛАВА 2. МАТЕМАТИЧЕСКОЕ МОДЕЛИРОВАНИЕ МЕХАНИЧЕСКОГО ПОВЕДЕНИЯ СЛОИСТЫХ ПЛАСТИН С ТЕТРАКИРАЛЬНЫМИ СОТАМИ МЕТОДОМ КОНЕЧНЫХ ЭЛЕМЕНТОВ}
\else
\phantomsection
\section[\texorpdfstring{\protect\mbox{\protect\scalebox{0.95}[1]{ГЛАВА 2. МАТЕМАТИЧЕСКОЕ МОДЕЛИРОВАНИЕ МЕХАНИЧЕСКОГО}} ПОВЕДЕНИЯ СЛОИСТЫХ ПЛАСТИН С ТЕТРАКИРАЛЬНЫМИ СОТАМИ МЕТОДОМ КОНЕЧНЫХ ЭЛЕМЕНТОВ}{ГЛАВА 2. МАТЕМАТИЧЕСКОЕ МОДЕЛИРОВАНИЕ МЕХАНИЧЕСКОГО ПОВЕДЕНИЯ СЛОИСТЫХ ПЛАСТИН С ТЕТРАКИРАЛЬНЫМИ СОТАМИ МЕТОДОМ КОНЕЧНЫХ ЭЛЕМЕНТОВ}]{ГЛАВА 2. МАТЕМАТИЧЕСКОЕ МОДЕЛИРОВАНИЕ МЕХАНИЧЕСКОГО ПОВЕДЕНИЯ СЛОИСТЫХ ПЛАСТИН С ТЕТРАКИРАЛЬНЫМИ СОТАМИ МЕТОДОМ КОНЕЧНЫХ ЭЛЕМЕНТОВ}
\fi

В данной главе выполнено численное моделирование напряженно-деформированного состояния трехслойных пластин с тетракиральным сотовым заполнителем при статическом изгибе путем трехмерного конечно-элементного моделирования в системе Comsol Multiphysics с использованием уравнений линейной теории упругости. Также выполнено построение математических конечно-элементных моделей слоистой пластины с тетракиральными сотовыми прослойками с применением совместных и несовместных плоских конечных элементов для описания напряженно-деформированного состояния в условиях статического воздействия в рамках теории упругости. На основе построенных моделей осуществлена разработка алгоритмов численного моделирования напряженно-деформированного состояния многослойных пластин с тетракиральными сотами при статическом изгибе.

\iflatexml\else\phantomsection\fi
\subsection{2.1. Моделирование напряженного состояния трехслойных композитов с тетракиральными сотами в условиях статического изгиба путем трехмерного конечно-элементного моделирования в системе Comsol Multiphysics}

\begin{figure}[H]
\centering
\includegraphics[keepaspectratio,width=6.29921in,height=3.24717in]{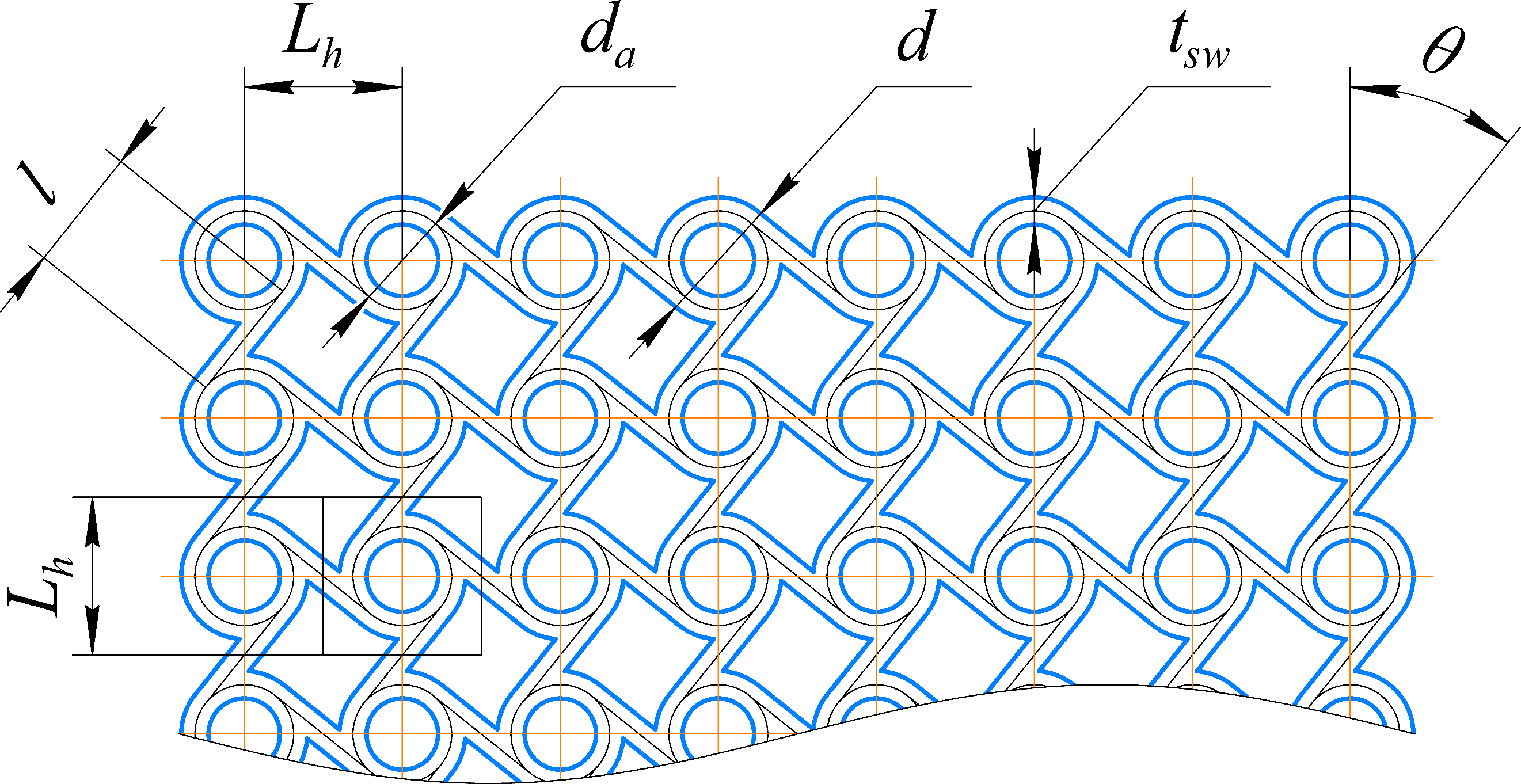}
\caption{Параметры сотовой структуры тетракирального типа \dualcite{ref76,ref77,ref78}{{[}\bibnum{ref76}--\bibnum{ref78}{]}}{{[}\bibnum{ref76}--\bibnum{ref78}{]}}}\label{fig:1}
\end{figure}

Рассмотрим тетракиральные сотовые заполнители с четырьмя значениями размера элементарных ячеек \(L_{h} = 1.6d_{a},\) \(d_{a} \in \left\{ 1, 1.3, 1.6, 1.9 \right\}\) мм \dualcite{ref76,ref77,ref78}{{[}\bibnum{ref76}--\bibnum{ref78}{]}}{{[}\bibnum{ref76}--\bibnum{ref78}{]}}. Рассматриваемые соты имеют постоянное отношение среднего радиуса цилиндров к длине тангенциально прикрепленных ребер \(r_{a}/l = \text{const}\) \((r_{a} = d_{a}/2),\) в связи с чем у попарно связанных цилиндров постоянен угол между средней линией ребер и отрезком, соединяющим центры цилиндров \(\theta = \arctan\left( 2r_{a}/l \right) = \text{const}\) (рисунок \ref{fig:1}).

\begin{figure}[H]
\centering
\includegraphics[keepaspectratio,width=6.69375in,height=2.68403in]{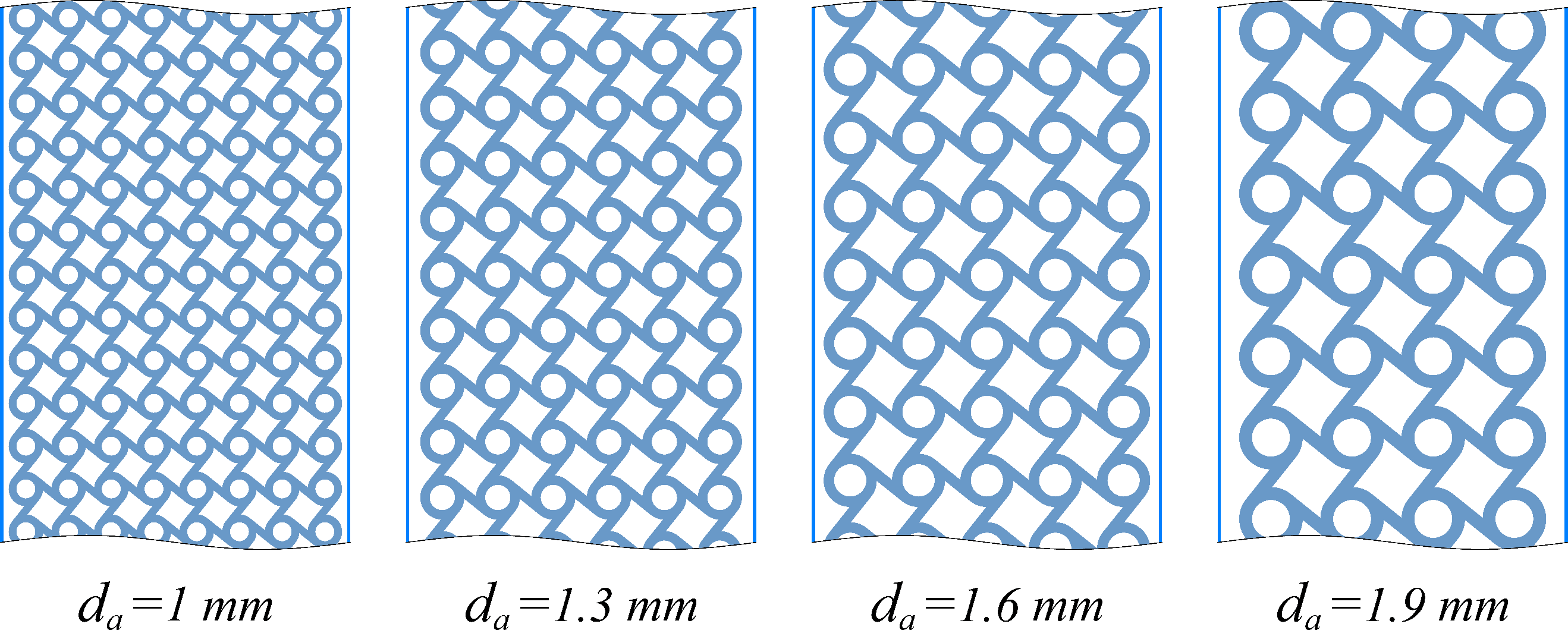}
\caption{Тетракиральные соты с разной дискретизацией и равной относительной плотностью \(\rho_{rel} = 35.3\)~\% \dualcite{ref78}{{[}\bibnum{ref78}{]}}{{[}\bibnum{ref78}{]}}}\label{fig:2}
\end{figure}

При этом соты имеют разную дискретизацию, а также равный диапазон изменения относительной плотности \(\rho_{rel}\) (\%), которая определяется как отношение площади твердого тела сот в плоскости к площади основания пластин (рисунок~\ref{fig:2}). Относительная плотность варьируется путем изменения толщины стенок сот \(t_{sw}.\)

При дополнительном условии \(\rho_{rel} = \text{const}\) у рассматриваемых сот также приблизительно постоянны отношения длины тангенциально прикрепленных ребер к внешнему радиусу цилиндров \(\alpha = l/r \approx \text{const}\) \((r = d/2)\) и толщины стенок сот к внешнему радиусу цилиндров \(\beta = t_{sw}/r \approx \text{const}\text{.}\)

При каждом из четырех значений \(L_{h}\) тетракиральные соты равномерно заполняют центральный слой композитов (рисунок~\ref{fig:2}).

\vspace{-0.575\baselineskip}

\begin{figure}[H]
\centering
\includegraphics[keepaspectratio,width=6.69291in,height=2.94618in]{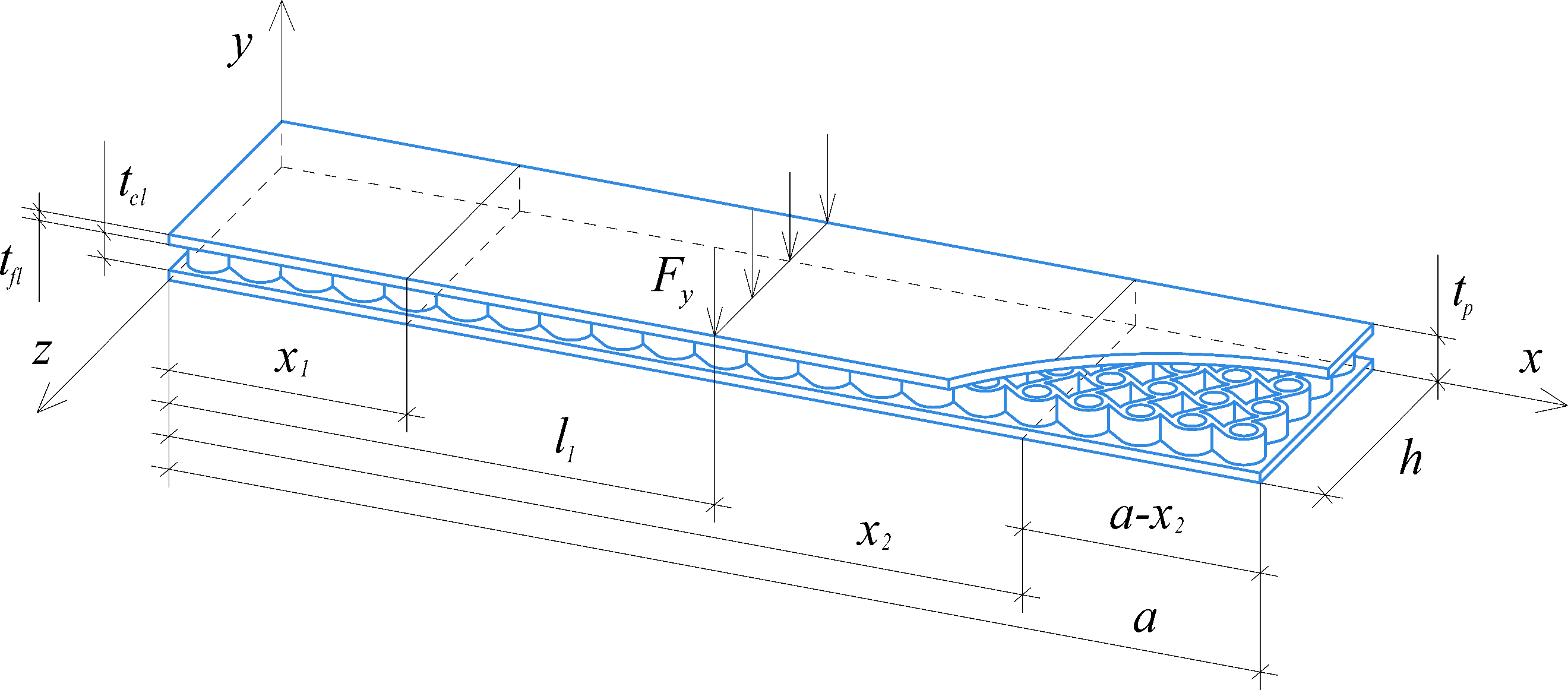}
\caption{Граничные условия, схема нагружения и параметры композитных пластин \dualcite{ref78}{{[}\bibnum{ref78}{]}}{{[}\bibnum{ref78}{]}}}\label{fig:3}
\end{figure}

Рассмотрим две постановки численных экспериментов \dualcite{ref77,ref78}{{[}\bibnum{ref77}, \bibnum{ref78}{]}}{{[}\bibnum{ref77}, \bibnum{ref78}{]}}. В первой постановке принимаются постоянными толщины слоев композитных пластин \(t_{fl}\) и \(t_{cl}\) (рисунок \ref{fig:3}) при изменении относительной плотности сот, а во второй постановке добавляется условие постоянного объема твердого тела сот, вследствие чего изменяется их толщина \(t_{cl}\) с изменением \(\rho_{rel}.\) Согласно принятым соотношениям \dualcite{ref1,ref28,ref34}{{[}\bibnum{ref1}, \bibnum{ref28}, \bibnum{ref34}{]}}{{[}\bibnum{ref1}, \bibnum{ref28}, \bibnum{ref34}{]}}, у тонких пластин отношение толщины к наименьшему размеру в плоскости \(t_{p}/h\) (рисунок \ref{fig:3}) должно быть менее \(1/5,\) чему соответствуют композиты в первой постановке экспериментов. Однако во второй постановке соотношение \(t_{p}/h\) у композитов при малых значениях \(\rho_{rel}\) характерно для относительно толстой пластины, а при больших значениях \(\rho_{rel}\) -- для тонкой пластины. При этом использовались следующие параметры пластин (рисунок \ref{fig:3}): \emph{a} -- длина, \emph{h} -- ширина, \(t_{fl}\) -- толщина внешних слоев, \(t_{cl}\) -- толщина сотового заполнителя, \(t_{p}\) -- общая толщина.

Нагружение пластины силой \(F_{y}\), приложенной к верхней грани по линии с координатой \(x = l_{1}\), эквивалентной распределенной нагрузке \(q_{y}\) вдоль ширины \emph{h} (рисунок \ref{fig:3}), является условием ее статического изгиба:

\begin{equation}
F_{y}/h = q_{y}\left( x = l_{1},\; y = t_{p},\; 0 \leq z \leq h \right).
\end{equation}

Будем рассматривать два варианта граничных условий.

При жестком защемлении пластины с двух сторон граничные условия имеют вид:

\begin{equation}
\label{eq:18}
\begin{matrix}
u_{x,y,z}\left( 0 \leq x \leq x_{1},\; y = 0,\; 0 \leq z \leq h \right) = 0, \\
\begin{matrix}
u_{x,y,z}\left( x_{2} \leq x \leq a,\; y = 0,\; 0 \leq z \leq h \right) = 0, \\
\begin{matrix}
u_{x,y,z}\left( 0 \leq x \leq x_{1},\; y = t_{p},\; 0 \leq z \leq h \right) = 0, \\
u_{x,y,z}\left( x_{2} \leq x \leq a,\; y = t_{p},\; 0 \leq z \leq h \right) = 0,
\end{matrix}
\end{matrix}
\end{matrix}
\end{equation}

\noindent где \(0 \leq x \leq x_{1}\) и \(x_{2} \leq x \leq a\) -- участки защемления.

При моделировании трехточечного изгиба пластины будем полагать, что она опирается нижней гранью по линиям с координатами \(x = x_{1}\) и \(x = x_{2}\) по всей ширине, допуская упругий поворот на опорах \dualcite{ref2,ref29}{{[}\bibnum{ref2}, \bibnum{ref29}{]}}{{[}\bibnum{ref2}, \bibnum{ref29}{]}}. В этом случае граничные условия имеют вид:

\begin{equation}
\begin{matrix}
u_{x,y,z}\left( x = x_{1},\; y = 0,\; 0 \leq z \leq h \right) = 0, \\
u_{x,y,z}\left( x = x_{2},\; y = 0,\; 0 \leq z \leq h \right) = 0.
\end{matrix}
\end{equation}

Таким образом, используется приближенная кинематическая схема, в которой повороты в окрестностях опор возникают вследствие упругого деформирования материала пластины.

Для ограничения прогиба композитных пластин в плоскости \emph{zy} во внешних слоях перемещения узлов вдоль координаты \emph{z} приравнивались нулю \(u_{z} = 0.\)

Известно, что метод конечных элементов является мощным инструментом решения задач математической физики в различных областях исследований \dualcite{ref26,ref31,ref32,ref33,ref35,ref38,ref47,ref50,ref59,ref65,ref68,ref87}{{[}\bibnum{ref26}, \bibnum{ref31}--\bibnum{ref33}, \bibnum{ref35}, \bibnum{ref38}, \bibnum{ref47}, \bibnum{ref50}, \bibnum{ref59}, \bibnum{ref65}, \bibnum{ref68}, \bibnum{ref87}{]}}{{[}\bibnum{ref26}, \bibnum{ref31}--\bibnum{ref33}, \bibnum{ref35}, \bibnum{ref38}, \bibnum{ref47}, \bibnum{ref50}, \bibnum{ref59}, \bibnum{ref65}, \bibnum{ref68}, \bibnum{ref87}{]}}.

В данной работе сначала было выполнено трехмерное конечно-элементное моделирование в модуле «Механика конструкций» системы «Comsol Multiphysics 5.6» \dualcite{ref96}{{[}\bibnum{ref96}{]}}{{[}\bibnum{ref96}{]}} с использованием уравнений линейной теории упругости.

\begin{figure}[H]
\centering
\includegraphics[keepaspectratio,width=6.69375in,height=3.35625in]{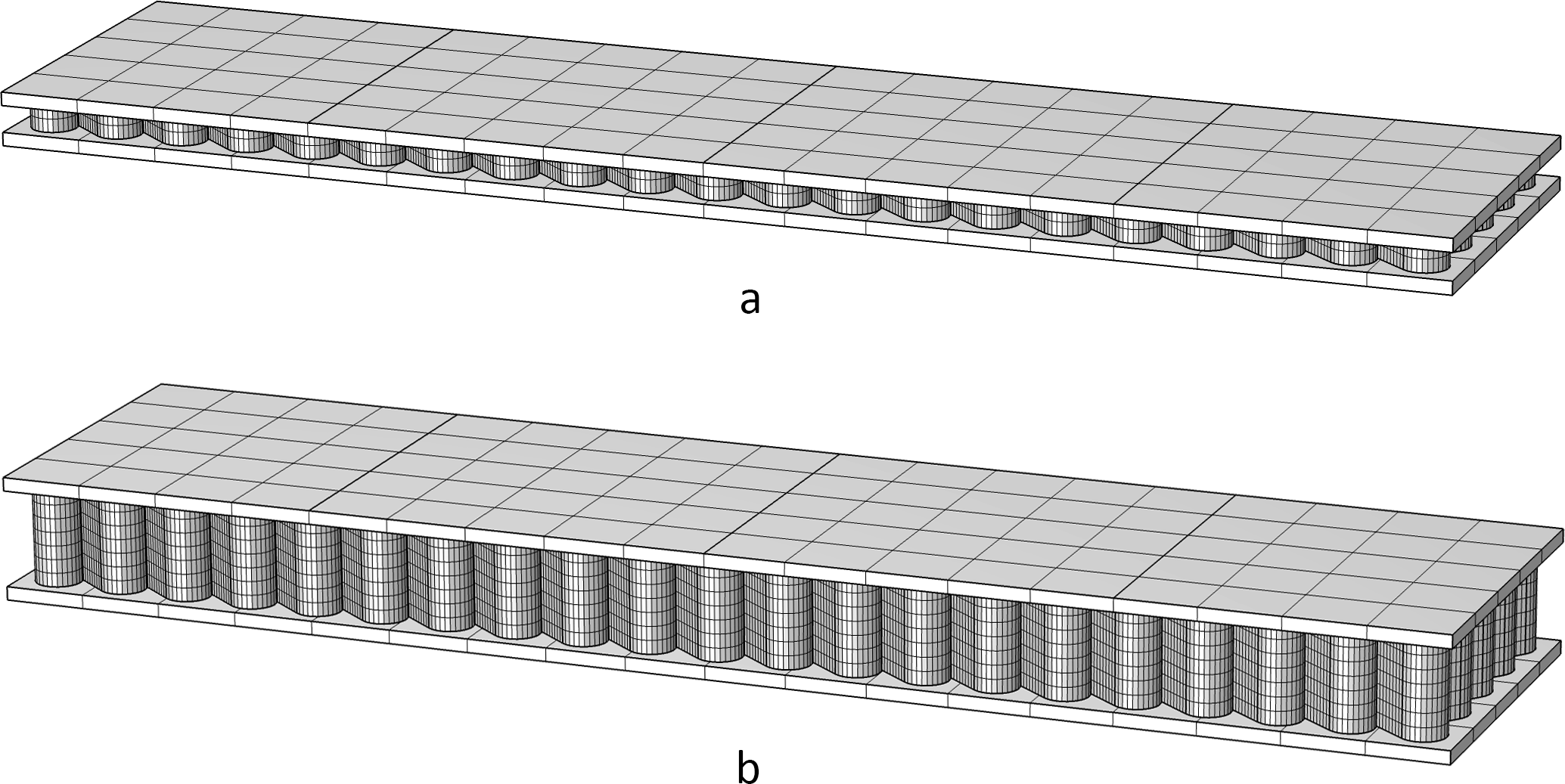}
\caption{Сетка конечных элементов тонкой (a) и относительно толстой (b) композитных пластин \dualcite{ref77}{{[}\bibnum{ref77}{]}}{{[}\bibnum{ref77}{]}}}\label{fig:4}
\end{figure}

В этом случае сетка конечных элементов композитных пластин построена отдельно для каждого слоя \dualcite{ref77}{{[}\bibnum{ref77}{]}}{{[}\bibnum{ref77}{]}}: четырехугольные призмы использованы для сплошных слоев, а треугольные призмы -- для сотовых прослоек (рисунок \ref{fig:4}). На границах сопряжения слоев принято условие непрерывности полевых переменных.

\begin{figure}[H]
\centering
\includegraphics[keepaspectratio,width=5.51181in,height=5.51181in]{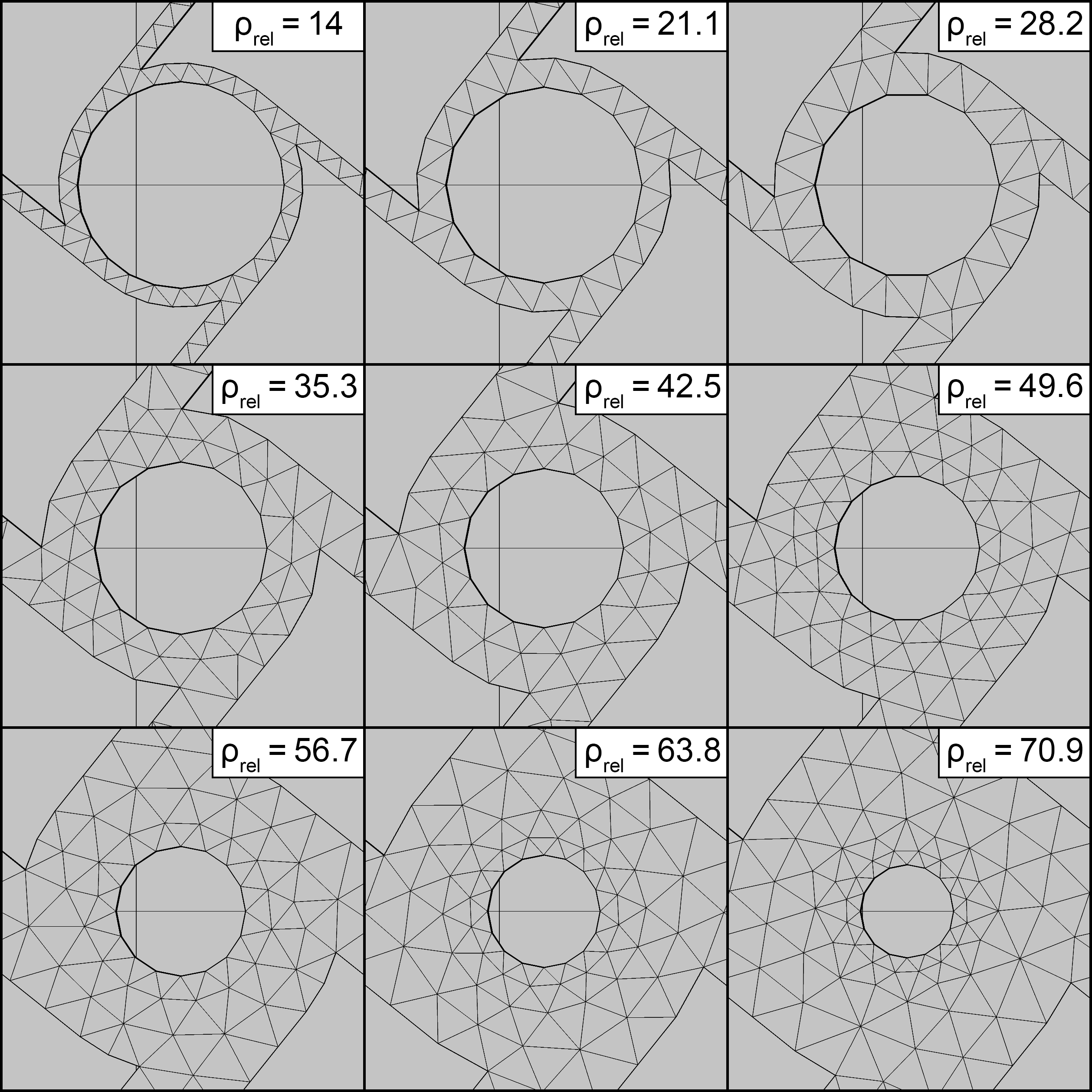}
\caption{Сетка конечных элементов элементарных ячеек тетракиральных структур при \(14 \leq \rho_{rel} \leq 70.9\) \% \dualcite{ref77}{{[}\bibnum{ref77}{]}}{{[}\bibnum{ref77}{]}}}\label{fig:5}
\end{figure}

При построении сетки конечных элементов применялись общие принципы. Тетракиральные соты в плоскости разбивались на минимально необходимое количество наиболее правильных треугольников при сохранении геометрической формы сот (рисунок \ref{fig:5}), а внешние слои композитов в плоскости разбивались на одинаковое количество наиболее правильных четырехугольников. При этом в первой постановке экспериментов тетракиральные соты содержат два слоя конечных элементов по толщине, а во второй постановке -- от 1 до 7 слоев согласно изменению толщины. Внешние слои композитов в обеих постановках содержат один слой конечных элементов по толщине. При использовании системы Comsol Multiphysics применялись конечные элементы серендипова семейства второго порядка. Таким образом, расчеты композитов должны иметь сопоставимую погрешность дискретизации, которая не препятствует анализу изменения напряженного состояния композитов.

\vspace{-0.575\baselineskip}

Эквивалентные напряжения согласно критерию предельного состояния Мизеса \dualcite{ref10,ref23}{{[}\bibnum{ref10}, \bibnum{ref23}{]}}{{[}\bibnum{ref10}, \bibnum{ref23}{]}} определяются соотношением:

\begin{equation}
\label{eq:20}
\sigma_{e} = \frac{1}{\sqrt{2}}\sqrt{\left( \sigma_{1} - \sigma_{2} \right)^{2} + \left( \sigma_{2} - \sigma_{3} \right)^{2} + \left( \sigma_{3} - \sigma_{1} \right)^{2}},
\end{equation}

\noindent где \(\sigma_{1}\), \(\sigma_{2}\) и \(\sigma_{3}\) -- главные напряжения.

\clearpage
\iflatexml\else\phantomsection\fi
\subsection{2.2. Моделирование напряженного состояния трехслойных композитов с тетракиральными сотами в условиях статического изгиба посредством алгоритмов решения плоской задачи методом конечных элементов}

Далее было выполнено моделирование сплошных пластин с помощью алгоритма решения плоской задачи методом конечных элементов в перемещениях \dualcite{ref5,ref6,ref7,ref8,ref16,ref20,ref21,ref25,ref27,ref76,ref77}{{[}\bibnum{ref5}--\bibnum{ref8}, \bibnum{ref16}, \bibnum{ref20}, \bibnum{ref21}, \bibnum{ref25}, \bibnum{ref27}, \bibnum{ref76}, \bibnum{ref77}{]}}{{[}\bibnum{ref5}--\bibnum{ref8}, \bibnum{ref16}, \bibnum{ref20}, \bibnum{ref21}, \bibnum{ref25}, \bibnum{ref27}, \bibnum{ref76}, \bibnum{ref77}{]}}. В качестве расчетной области принималась боковая грань пластины в плоскости \emph{xy} (рисунок \ref{fig:3}). В случае жесткого защемления вместо граничных условий  (\ref{eq:18}) использовались условия:

\begin{equation}
\begin{matrix}
u_{x,y,z}\left( x = x_{1},\; 0 \leq y \leq t_{p},\; 0 \leq z \leq h \right) = 0, \\
u_{x,y,z}\left( x = x_{2},\; 0 \leq y \leq t_{p},\; 0 \leq z \leq h \right) = 0.
\end{matrix}
\end{equation}

\begin{figure}[H]
\centering
\includegraphics[keepaspectratio,width=4.62598in,height=3.45559in]{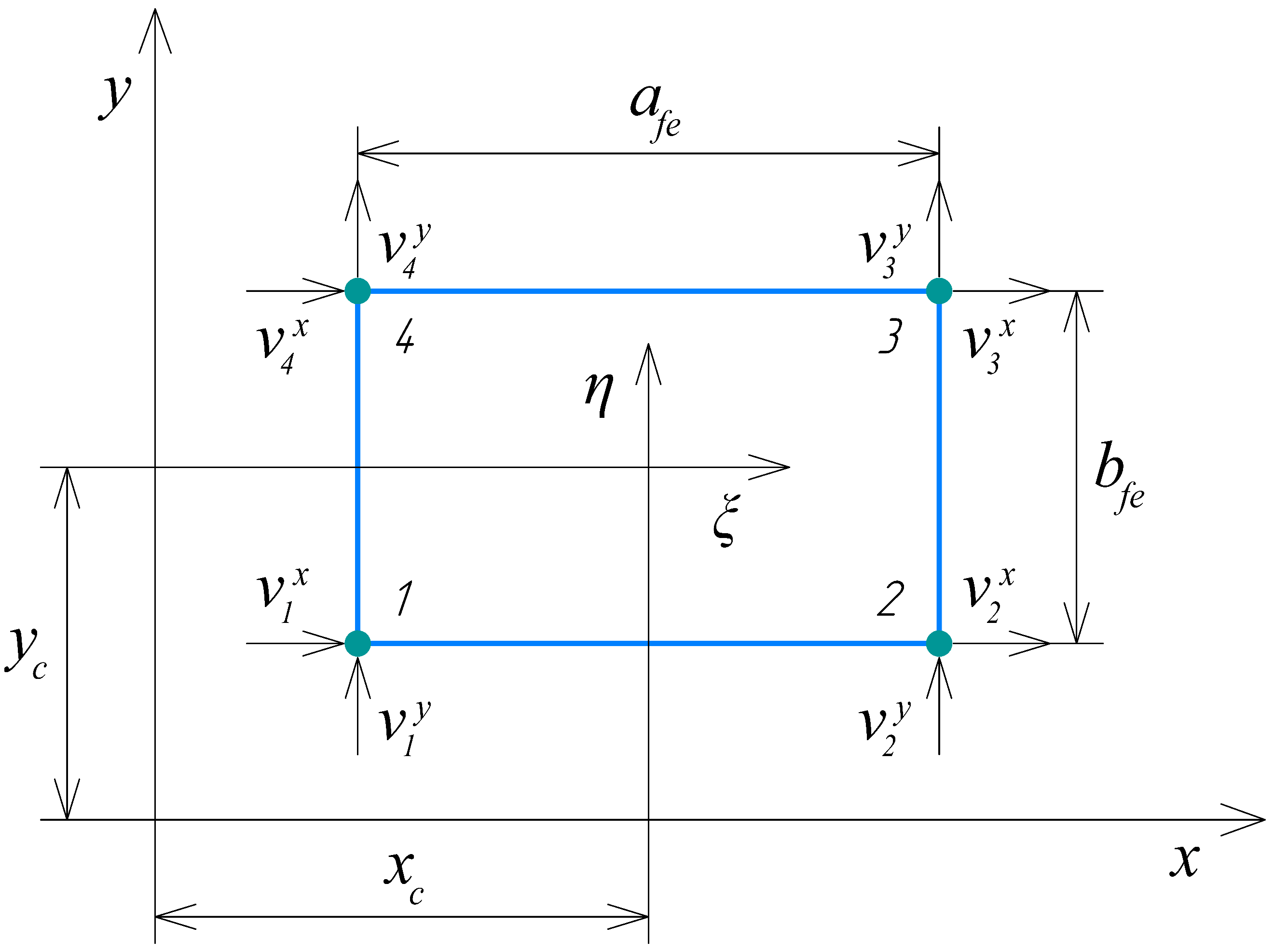}
\caption{Параметры прямоугольного конечного элемента \dualcite{ref78}{{[}\bibnum{ref78}{]}}{{[}\bibnum{ref78}{]}}}\label{fig:6}
\end{figure}

С этой целью расчетная область разбивалась на прямоугольные конечные элементы с четырьмя узлами в вершинах (рисунок \ref{fig:6}).

На рисунке \ref{fig:6} \(a_{fe}\) и \(b_{fe}\) -- размеры сторон прямоугольного элемента, \(x_{c}\) и \(y_{c}\) -- координаты его центра, \(\xi = 2\left( x - x_{c} \right)/a_{fe}\) и \(\eta = 2\left( y - y_{c} \right)/b_{fe}\) -- безразмерные координаты прямоугольного элемента, \(\xi_{1} = - 1,\) \(\eta_{1} = - 1,\) \(\xi_{2} = 1,\) \(\eta_{2} = - 1,\) \(\xi_{3} = 1,\) \(\eta_{3} = 1,\) \(\xi_{4} = - 1,\) \(\eta_{4} = 1,\) \(v^{x}\) и \(v^{y}\) -- компоненты перемещений узлов конечного элемента вдоль осей \emph{x} и \emph{y} соответственно.

В этом случае аналитическая форма матрицы жесткости конечного элемента \(k^{e}\) определяется путем интегрирования:

\begin{equation}
k_{r,s}^{e} = \int_{}^{}{\beta_{r}^{T}\chi\beta_{s}hdxdy =}\frac{a_{fe}b_{fe}h}{4}\int_{- 1}^{1}{\int_{- 1}^{1}{\beta_{r}^{T}\chi\beta_{s}d\xi d\eta}},
\end{equation}

\noindent где \(\beta_{r}\) (а также \(\beta_{s}\)) -- матрица связи между узловыми перемещениями и деформациями:

\begin{equation}
\beta_{r} = L\alpha_{r} = \begin{pmatrix}
\frac{\partial}{\partial x} & 0 \\
0 & \frac{\partial}{\partial y} \\
\frac{\partial}{\partial y} & \frac{\partial}{\partial x}
\end{pmatrix}\begin{pmatrix}
\psi_{r} & 0 \\
0 & \psi_{r}
\end{pmatrix} = \frac{1}{2}\begin{pmatrix}
b_{r} & 0 \\
0 & a_{r} \\
a_{r} & b_{r}
\end{pmatrix},
\end{equation}

\noindent \(\alpha_{r}\) -- матрица аппроксимирующих функций, \(\psi_{r} = \left( 1 + \xi_{r}\xi \right)\left( 1 + \eta_{r}\eta \right)/4,\) \(a_{r} = \eta_{r}\left( 1 + \xi_{r}\xi \right)/b_{fe},\) \(b_{r} = \xi_{r}\left( 1 + \eta_{r}\eta \right)/a_{fe},\) \emph{r} и \emph{s} -- номера блоков матриц \(\left( r,s \in \left\{ 1,2,\ldots,4 \right\} \right),\) \emph{h} -- размер конечного элемента по оси \emph{z}, который равен ширине пластины, \emph{\(\chi\)} -- матрица упругих постоянных для плоской деформации из выражения  (\ref{eq:6}).

Матрицу жесткости удобно разбить на два слагаемых:

\begin{equation}
k_{r,s}^{e} = k_{r,s}^{E} + k_{r,s}^{G},
\end{equation}

\noindent где \(k^{E}\) -- подматрица нормальных деформаций:

\begin{equation}
\label{eq:25}
k_{r,s}^{E} = \frac{a_{fe}b_{fe}h}{4}\int_{- 1}^{1}{\int_{- 1}^{1}{\beta_{r}^{T}\chi_{E}\beta_{s}d\xi d\eta}}
\end{equation}

\noindent и \(k^{G}\) -- подматрица деформаций сдвига:

\begin{equation}
\label{eq:26}
k_{r,s}^{G} = \frac{a_{fe}b_{fe}h}{4}\int_{- 1}^{1}{\int_{- 1}^{1}{\beta_{r}^{T}\chi_{G}\beta_{s}d\xi d\eta}},
\end{equation}

\noindent \(\chi_{E}\) и \(\chi_{G}\) -- матрицы упругих коэффициентов, которые определяются равенствами:

\begin{equation}
\chi_{E} = \frac{E}{(1 + \mu)(1 - 2\mu)}\begin{pmatrix}
1 - \mu & \mu & 0 \\
\mu & 1 - \mu & 0 \\
0 & 0 & 0
\end{pmatrix},
\end{equation}

\begin{equation}
\chi_{G} = G\begin{pmatrix}
0 & 0 & 0 \\
0 & 0 & 0 \\
0 & 0 & 1
\end{pmatrix},
\end{equation}

\noindent а \emph{G} -- модуль сдвига.

Вычисляя интегралы в  (\ref{eq:25}) и  (\ref{eq:26}), соответственно получим:

\begin{equation}
k_{r,s}^{E} = \frac{Eh}{4(1 + \mu)(1 - 2\mu)}\begin{bmatrix}
(1 - \mu)\gamma\xi_{r}\xi_{s}\left( 1 + \frac{\eta_{r}\eta_{s}}{3} \right) & \mu\xi_{r}\eta_{s} \\
\mu\eta_{r}\xi_{s} & (1 - \mu)\frac{\eta_{r}\eta_{s}}{\gamma}\left( 1 + \frac{\xi_{r}\xi_{s}}{3} \right)
\end{bmatrix},
\end{equation}

\begin{equation}
k_{r,s}^{G} = \frac{Gh}{4}\begin{bmatrix}
\frac{\eta_{r}\eta_{s}}{\gamma}\left( 1 + \frac{\xi_{r}\xi_{s}}{3} \right) & \eta_{r}\xi_{s} \\
\xi_{r}\eta_{s} & \gamma\xi_{r}\xi_{s}\left( 1 + \frac{\eta_{r}\eta_{s}}{3} \right)
\end{bmatrix},
\end{equation}

\noindent сумма которых определяет аналитическое выражение матрицы жесткости совместного конечного элемента, где \(\gamma = b_{fe}/a_{fe}\) -- безразмерный параметр.

\begin{figure}[H]
\centering
\includegraphics[keepaspectratio,width=5.90551in,height=2.95243in]{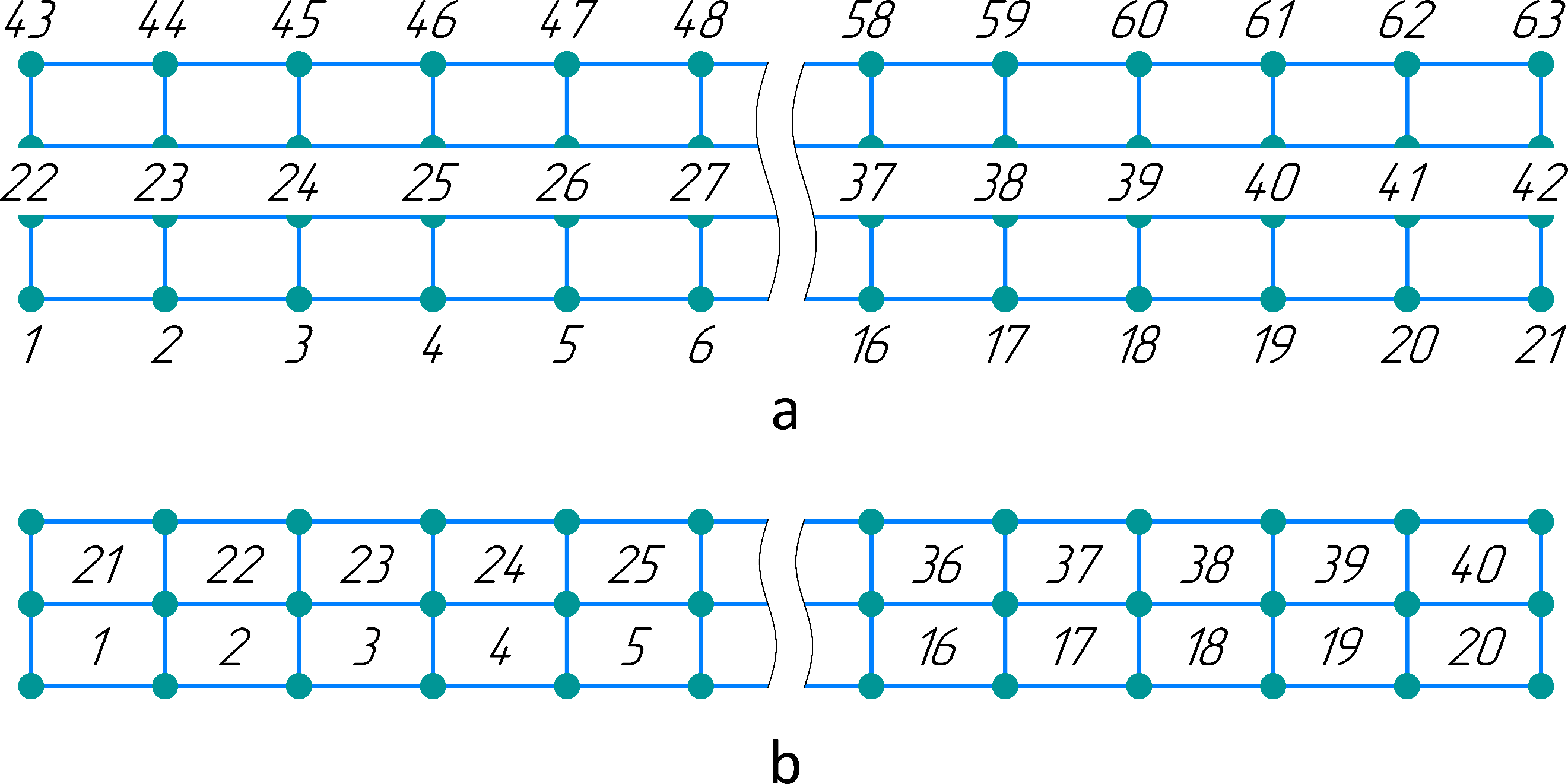}
\caption{Условная схема конечных элементов пластины: а) глобальная нумерация узлов, б) нумерация конечных элементов}\label{fig:7}
\end{figure}

\begin{figure}[H]
\centering
\vspace{0.38\baselineskip}
\includegraphics[keepaspectratio,width=5.70866in,height=1.10951in]{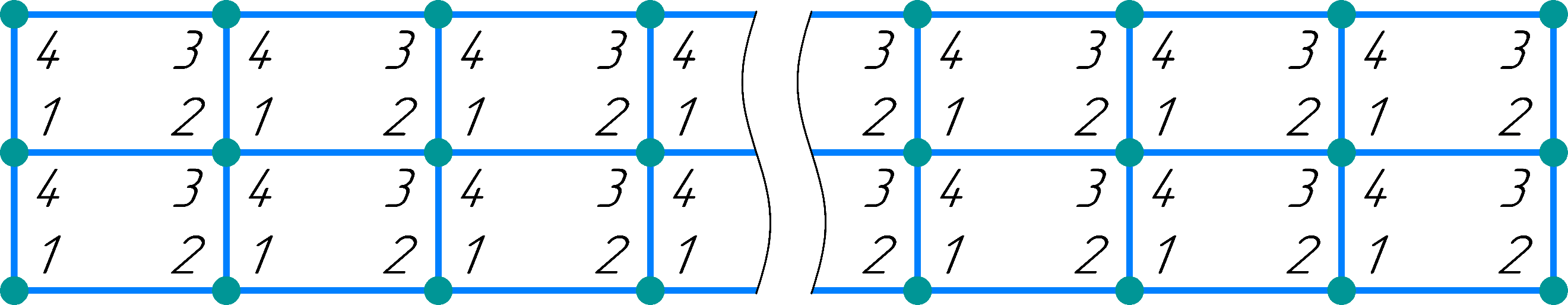}
\caption{Схема локальной нумерации узлов конечных элементов}\label{fig:8}
\end{figure}

Матрица соответствия локальных номеров узлов глобальным \emph{A} строится по принципу \(A_{m,i} = q\) \(\left( m \in \left\{ 1,2,\ldots,m_{f} \right\},\ i \in \left\{ 1,2,\ldots,i_{f} \right\},\ q \in \left\{ 1,2,\ldots,4 \right\} \right),\) где \emph{m} -- глобальный номер узла (рисунок \ref{fig:7}a), \(m_{f}\) -- количество глобальных узлов, \emph{i} -- номер конечного элемента (рисунок \ref{fig:7}b), \(i_{f}\) -- количество конечных элементов, \emph{q} -- локальный номер узла \emph{i}-го конечного элемента (рисунок \ref{fig:8}), если \(A_{m,i} \neq q,\) то \(A_{m,i} = 0.\)

Расширенная матрица жесткости \(k^{\exp}\) строится по принципу \(k^{\exp}(i)_{m,n} = k_{r,s}^{e},\) где \(r = A_{m,i},\) \(s = A_{n,i},\) \(m,n \in \left\{ 1,2,\ldots,m_{f} \right\}.\) Если \(r \vee s = 0,\) то \(k_{r,s}^{e} = \begin{pmatrix}
0 & 0 \\
0 & 0
\end{pmatrix}.\)

Тогда матрицу жесткости конечно-элементной модели \emph{K} можно определить суммированием расширенных матриц жесткости \(K = \sum_{i}^{}{k^{\exp}(i)}.\) Для учета внешнего закрепления узла конечно-элементной модели необходимо удалить строки \(i_{1} = 2m_{p} - 1,\) \(i_{2} = 2m_{p}\) и столбцы \(j_{1} = 2m_{p} - 1,\) \(j_{2} = 2m_{p}\) матрицы жесткости \emph{K}, где \(m_{p}\) -- номер закрепленного узла.

Перемещения узлов определяются выражением:

\begin{equation}
u_{a} = K_{a}^{- 1}P_{a},
\end{equation}

\noindent где \(K_{a}^{- 1}\) -- матрица, обратная матрице жесткости модели с учетом закрепленных узлов, \(P_{a}\) -- вектор узловых сил, \(P_{a} = \begin{pmatrix}
P_{a_{1}}^{x} & P_{a_{2}}^{y} & \begin{matrix}
\cdots & P_{a_{c - 1}}^{x} & P_{a_{c}}^{y}
\end{matrix}
\end{pmatrix}^{T},\) \(c = 2\left( m_{f} - p \right),\) \(P_{a_{c - 1}}^{x}\) и \(P_{a_{c}}^{y}\) -- узловые силы по осям \emph{x} и \emph{y} соответственно, \emph{p} -- количество закрепленных узлов.

\vspace{-0.575\baselineskip}

Полный вектор перемещений \(u = \begin{pmatrix}
u_{1}^{x} & u_{2}^{y} & \begin{matrix}
\cdots & u_{e - 1}^{x} & u_{e}^{y}
\end{matrix}
\end{pmatrix}^{T}\) \(\left( e = 2m_{f} \right)\) включает нулевые перемещения \(u_{k_{\varsigma}} = 0\) \(\left( \varsigma \in \left\{ 1,2 \right\} \right),\) где \(k_{1} = 2m_{p} - 1,\) \(k_{2} = 2m_{p},\) вектор \(u_{a} = \begin{pmatrix}
u_{a_{1}}^{x} & u_{a_{2}}^{y} & \begin{matrix}
\cdots & u_{a_{c - 1}}^{x} & u_{a_{c}}^{y}
\end{matrix}
\end{pmatrix}^{T}\) является подвектором \emph{u}, где \(u_{a_{c}}\) -- свободное перемещение.

Векторы перемещений по осям \emph{x} и \emph{y} определяются выражениями \(u_{m}^{x} = u_{m_{x}}\) и \(u_{m}^{y} = u_{m_{y}}\) соответственно, где \(m_{x} = 2m - 1,\) \(m_{y} = 2m.\)

Вектор перемещений \emph{v} узлов конечного элемента строится по принципу \(v(i) = \begin{pmatrix}
u_{m^{1}}^{x} & u_{m^{1}}^{y} & u_{m^{2}}^{x} & \begin{matrix}
u_{m^{2}}^{y} & u_{m^{3}}^{x} & u_{m^{3}}^{y} & \begin{matrix}
u_{m^{4}}^{x} & u_{m^{4}}^{y}
\end{matrix}
\end{matrix}
\end{pmatrix}^{T},\) где \(A_{m^{1},i} = 1,\) \(A_{m^{2},i} = 2,\) \(A_{m^{3},i} = 3,\) \(A_{m^{4},i} = 4,\) \(u_{m^{q}}^{x} = v_{q}^{x}\) и \(u_{m^{q}}^{y} = v_{q}^{y}\) -- узловые перемещения по осям \emph{x} и \emph{y} соответственно.

Вектор деформаций \(\varepsilon\) конечного элемента определяется из выражения:

\begin{equation}
\varepsilon(i,\xi,\eta) = \beta(\xi,\eta)v(i),
\end{equation}

\noindent где \(\beta(\xi,\eta)\) -- матрица связи между узловыми перемещениями и деформациями:

\begin{equation}
\beta(\xi,\eta) = \frac{1}{2}\begin{pmatrix}
b_{a}(1,\eta) & 0 & a_{a}(1,\xi) \\
0 & a_{a}(1,\xi) & b_{a}(1,\eta) \\
b_{a}(2,\eta) & 0 & a_{a}(2,\xi) \\
0 & a_{a}(2,\xi) & b_{a}(2,\eta) \\
b_{a}(3,\eta) & 0 & a_{a}(3,\xi) \\
0 & a_{a}(3,\xi) & b_{a}(3,\eta) \\
b_{a}(4,\eta) & 0 & a_{a}(4,\xi) \\
0 & a_{a}(4,\xi) & b_{a}(4,\eta)
\end{pmatrix}^{T},
\end{equation}

\noindent в которой \(a_{a}(q,\xi) = \eta_{q}\left( 1 + \xi_{q}\xi \right)/b_{fe},\) \(b_{a}(q,\eta) = \xi_{q}\left( 1 + \eta_{q}\eta \right)/a_{fe},\) \(q \in \left\{ 1,2,\ldots,4 \right\}.\)

Вектор узловых напряжений \(\sigma\) в вершинах \emph{i}-го конечного элемента определяется из выражения:

\begin{equation}
\sigma(i,\xi,\eta) = \chi\varepsilon(i,\xi,\eta),
\end{equation}

\noindent где \(\xi = \xi_{q}\) и \(\eta = \eta_{q},\) а \emph{\(\chi\)} -- матрица упругости из выражения  (\ref{eq:6}).

Эквивалентные напряжения \(\sigma_{e}\) в узлах конечного элемента определяются в соответствии с критерием Мизеса  (\ref{eq:20}). Следует отметить, что в алгоритмах решения плоской задачи в работах \dualcite{ref76,ref77,ref78}{{[}\bibnum{ref76}--\bibnum{ref78}{]}}{{[}\bibnum{ref76}--\bibnum{ref78}{]}} использовался «инженерный» подход, основанный на пренебрежении касательными напряжениями вследствие их малости в рамках постановок, поэтому в качестве главных принимались нормальные напряжения. Кроме того, с целью упрощения анализа и ввиду несущественной погрешности, в указанных алгоритмах использовался критерий Мизеса для плоского напряженного состояния. Данная работа, напротив, приведена к наибольшему соответствию теоретическим положениям.

Главные напряжения \(\sigma_{1}\), \(\sigma_{2}\) и \(\sigma_{3}\) выражаются через соотношения:

\iflatexml
\begin{equation}
\sigma_{1,2} = \frac{\sigma_{x} + \sigma_{y}}{2} \pm \sqrt{\left( \frac{\sigma_{x} - \sigma_{y}}{2} \right)^{2} + {\tau_{xy}}^{2}},\label{eq:35}
\end{equation}
\begin{equation}
\sigma_{3} = \sigma_{z} = \mu\left( \sigma_{x} + \sigma_{y} \right),
\end{equation}
\else
\begingroup
\setlength{\abovedisplayskip}{0.09\baselineskip}
\setlength{\belowdisplayskip}{0.43\baselineskip}
\setlength{\abovedisplayshortskip}{0.09\baselineskip}
\setlength{\belowdisplayshortskip}{0.43\baselineskip}
\begin{gather}
\sigma_{1,2} = \frac{\sigma_{x} + \sigma_{y}}{2} \pm \sqrt{\left( \frac{\sigma_{x} - \sigma_{y}}{2} \right)^{2} + {\tau_{xy}}^{2}},\label{eq:35}\\[1.00\baselineskip]
\sigma_{3} = \sigma_{z} = \mu\left( \sigma_{x} + \sigma_{y} \right),
\end{gather}
\endgroup
\fi

\noindent где \(\sigma_{x}\), \(\sigma_{y}\) и \(\sigma_{z}\) -- нормальные напряжения по осям \emph{x}, \emph{y} и \emph{z} соответственно, \(\tau_{xy}\) -- касательное напряжение.

Также с целью верификации конечно-элементного моделирования сплошных пластин был применен дополнительный алгоритм решения плоской задачи, выражения которого получены с использованием матрицы аппроксимирующих функций, подробно описанной в \dualcite{ref21}{{[}\bibnum{ref21}{]}}{{[}\bibnum{ref21}{]}}:

\begin{equation}
\alpha_{r}^{inc} = \begin{pmatrix}
\frac{\left( 1 + \xi_{r}\xi \right)\left( 1 + \eta_{r}\eta \right)}{4} & - \frac{\xi_{r}\eta_{r}}{8}\left( \frac{\mu}{\gamma}\xi^{2} + \gamma\eta^{2} - \frac{\mu}{\gamma} - \gamma \right) \\
 - \frac{\xi_{r}\eta_{r}}{8}\left( \frac{1}{\gamma}\xi^{2} + \mu\gamma\eta^{2} - \frac{1}{\gamma} - \mu\gamma \right) & \frac{\left( 1 + \xi_{r}\xi \right)\left( 1 + \eta_{r}\eta \right)}{4}
\end{pmatrix}.
\end{equation}

Повторяя предыдущие выкладки, получим аналитические выражения для подматриц нормальных и сдвиговых деформаций матрицы жесткости несовместного конечного элемента:

\begin{equation}
\begin{matrix}
k_{r,s}^{inc.E} = \frac{Eh}{4(1 + \mu)(1 - 2\mu)} \times \\
 \times \begin{bmatrix}
\gamma\xi_{r}\xi_{s}\left( 1 - \mu + \frac{1 - \mu - \mu^{2} - \mu^{3}}{3}\eta_{r}\eta_{s} \right) & \mu\xi_{r}\eta_{s} \\
\mu\eta_{r}\xi_{s} & \frac{\eta_{r}\eta_{s}}{\gamma}\left( 1 - \mu + \frac{1 - \mu - \mu^{2} - \mu^{3}}{3}\xi_{r}\xi_{s} \right)
\end{bmatrix},
\end{matrix}
\end{equation}

\begin{equation}
k_{r,s}^{inc.G} = \frac{Gh}{4}\begin{bmatrix}
\frac{\eta_{r}\eta_{s}}{\gamma} & \eta_{r}\xi_{s} \\
\xi_{r}\eta_{s} & \gamma\xi_{r}\xi_{s}
\end{bmatrix},
\end{equation}

\noindent а также выражение для матрицы связи между узловыми перемещениями и деформациями:

\begin{equation}
\beta^{inc}(\xi,\eta) = \frac{1}{2}\begin{pmatrix}
b_{a}(1,\eta) & e_{a}(1,\eta) & g_{a}(1) \\
c_{a}(1,\xi) & a_{a}(1,\xi) & f_{a}(1) \\
b_{a}(2,\eta) & e_{a}(2,\eta) & g_{a}(2) \\
c_{a}(2,\xi) & a_{a}(2,\xi) & f_{a}(2) \\
b_{a}(3,\eta) & e_{a}(3,\eta) & g_{a}(3) \\
c_{a}(3,\xi) & a_{a}(3,\xi) & f_{a}(3) \\
b_{a}(4,\eta) & e_{a}(4,\eta) & g_{a}(4) \\
c_{a}(4,\xi) & a_{a}(4,\xi) & f_{a}(4)
\end{pmatrix}^{T},
\end{equation}

\noindent где \(c_{a}(q,\xi) = - \mu\xi_{q}\eta_{q}\xi/b_{fe},\) \(e_{a}(q,\eta) = - \mu\xi_{q}\eta_{q}\eta/a_{fe},\) \(f_{a}(q) = \xi_{q}/a_{fe},\) \(g_{a}(q) = \eta_{q}/b_{fe}.\)

При расчете сплошных пластин с помощью вышеприведенных алгоритмов пластины разбивались по толщине на два слоя прямоугольных конечных элементов с приблизительно равными сторонами.

При получении алгоритма конечно-элементного моделирования композитных пластин \dualcite{ref17,ref18,ref78}{{[}\bibnum{ref17}, \bibnum{ref18}, \bibnum{ref78}{]}}{{[}\bibnum{ref17}, \bibnum{ref18}, \bibnum{ref78}{]}} используется матрица упругости трансверсально-изотропного материала для плоской деформации \dualcite{ref6,ref7,ref12}{{[}\bibnum{ref6}, \bibnum{ref7}, \bibnum{ref12}{]}}{{[}\bibnum{ref6}, \bibnum{ref7}, \bibnum{ref12}{]}}:

\begin{equation}
\label{eq:41}
\begin{matrix}
\chi = \frac{E_{2}}{\left( 1 + \mu_{1} \right)\left( 1 - \mu_{1} - 2n_{1}\mu_{2}^{2} \right)} \times \\
 \times \begin{pmatrix}
n_{1}\left( 1 - n_{1}\mu_{2}^{2} \right) & n_{1}\mu_{2}\left( 1 + \mu_{1} \right) & 0 \\
n_{1}\mu_{2}\left( 1 + \mu_{1} \right) & \left( 1 - \mu_{1}^{2} \right) & 0 \\
0 & 0 & m_{1}\left( 1 + \mu_{1} \right)\left( 1 - \mu_{1} - 2n_{1}\mu_{2}^{2} \right)
\end{pmatrix},
\end{matrix}
\end{equation}

\noindent \(n_{1} = E_{1}/E_{2},\) \(m_{1} = G_{2}/E_{2},\) \(E_{1}\) и \(\mu_{1}\) -- модуль упругости и коэффициент Пуассона соответственно, описывающие свойства материала в плоскости слоев, при этом \(G_{1} = E_{1}/[2\left( 1 + \mu_{1} \right)],\) а \(E_{2},\) \(\mu_{2}\) и \(G_{2}\) -- параметры материала из плоскости.

При этом матрицы упругих коэффициентов \(\chi_{E}\) и \(\chi_{G}\) определяются равенствами:

\begin{equation}
\chi_{E} = \frac{E_{2}}{\left( 1 + \mu_{1} \right)\left( 1 - \mu_{1} - 2n_{1}\mu_{2}^{2} \right)}\begin{pmatrix}
n_{1}\left( 1 - n_{1}\mu_{2}^{2} \right) & n_{1}\mu_{2}\left( 1 + \mu_{1} \right) & 0 \\
n_{1}\mu_{2}\left( 1 + \mu_{1} \right) & \left( 1 - \mu_{1}^{2} \right) & 0 \\
0 & 0 & 0
\end{pmatrix},
\end{equation}

\begin{equation}
\chi_{G} = G_{2}\begin{pmatrix}
0 & 0 & 0 \\
0 & 0 & 0 \\
0 & 0 & 1
\end{pmatrix}.
\end{equation}

При подстановке их в  (\ref{eq:25}) и  (\ref{eq:26}) соответственно получим:

\begin{equation}
\label{eq:44}
\begin{matrix}
k_{r,s}^{E} = \frac{E_{2}h}{4\left( 1 + \mu_{1} \right)\left( 1 - \mu_{1} - 2n_{1}\mu_{2}^{2} \right)} \times \\
 \times \begin{bmatrix}
\left( n_{1} - n_{1}^{2}\mu_{2}^{2} \right)\gamma\xi_{r}\xi_{s}\left( 1 + \frac{\eta_{r}\eta_{s}}{3} \right) & n_{1}\mu_{2}\left( 1 + \mu_{1} \right)\xi_{r}\eta_{s} \\
n_{1}\mu_{2}\left( 1 + \mu_{1} \right)\eta_{r}\xi_{s} & \left( 1 - \mu_{1}^{2} \right)\frac{\eta_{r}\eta_{s}}{\gamma}\left( 1 + \frac{\xi_{r}\xi_{s}}{3} \right)
\end{bmatrix},
\end{matrix}
\end{equation}

\begin{equation}
\label{eq:45}
k_{r,s}^{G} = \frac{G_{2}h}{4}\begin{bmatrix}
\frac{\eta_{r}\eta_{s}}{\gamma}\left( 1 + \frac{\xi_{r}\xi_{s}}{3} \right) & \eta_{r}\xi_{s} \\
\xi_{r}\eta_{s} & \gamma\xi_{r}\xi_{s}\left( 1 + \frac{\eta_{r}\eta_{s}}{3} \right)
\end{bmatrix}.
\end{equation}

При описании свойств тетракиральных сот используются аналитические выражения приведенных параметров, для справедливости которых в рамках модели требуется параллельность основного направления сот продольному сечению пластин. Приведенный модуль упругости в плоскости тетракиральных сот описывается выражением \dualcite{ref80}{{[}\bibnum{ref80}{]}}{{[}\bibnum{ref80}{]}}:

\begin{equation}
E_{1}^{tet} = \frac{E_{s}(t_{sw}/l)}{\cos^{2}\theta + \sin^{2}\theta/(t_{sw}/l)^{2}},
\end{equation}

\noindent а в перпендикулярном направлении к плоскости выражением \dualcite{ref58,ref73}{{[}\bibnum{ref58}, \bibnum{ref73}{]}}{{[}\bibnum{ref58}, \bibnum{ref73}{]}}:

\begin{equation}
E_{2}^{tet} = E_{s}\rho_{rel}^{tet},
\end{equation}

\noindent где \(\rho_{rel}^{tet}\) -- относительная плотность тетракиральных сот \dualcite{ref73}{{[}\bibnum{ref73}{]}}{{[}\bibnum{ref73}{]}}:

\begin{equation}
\rho_{rel}^{tet} = \frac{\beta\left\lbrack 2\alpha + \pi(2 - \beta) \right\rbrack - 2\left\lbrack \phi - (1 - \beta)\sin\phi \right\rbrack}{4\left\lbrack \left( 1 - \beta/2 \right)^{2} + \alpha^{2}/4 \right\rbrack},
\end{equation}

\noindent \(\phi = \text{arccos}(1 - \beta),\) \(E_{s}\) -- модуль упругости твердого тела сот.

В качестве приведенного модуля сдвига тетракиральных сот в направлении из плоскости воспользуемся верхней границей, определяемой выражением:

\begin{equation}
G_{2}^{tet.u} = G_{s}\beta\frac{\alpha\left( \cos^{2}\theta + \cos^{2}\left( \pi/2 - \theta \right) \right) + \pi}{4 + \alpha^{2}},
\end{equation}

\noindent где \(G_{s}\) -- модуль сдвига твердого тела.

Коэффициент Пуассона тетракиральных сот при деформировании в плоскости принимался равным нулю, \(\mu_{1} = 0\), согласно малодеформационной энергетической модели тетракиральной структуры \dualcite{ref80}{{[}\bibnum{ref80}{]}}{{[}\bibnum{ref80}{]}}. Такое допущение позволяет не использовать аналитические выражения, предложенные в работах \dualcite{ref74,ref85}{{[}\bibnum{ref74}, \bibnum{ref85}{]}}{{[}\bibnum{ref74}, \bibnum{ref85}{]}}, область применимости которых ограничена принятыми геометрическими и кинематическими предпосылками. При высокой относительной плотности сот, когда толщина стенок становится сопоставимой со свободной длиной тангенциальных ребер, эти предпосылки нарушаются, вследствие чего вычисленные значения эффективного коэффициента Пуассона могут утрачивать физический смысл. При нагружении в направлении, перпендикулярном плоскости сот, коэффициент Пуассона принимался равным коэффициенту Пуассона твердого тела \(\mu_{2} = \mu_{s}\) \dualcite{ref58}{{[}\bibnum{ref58}{]}}{{[}\bibnum{ref58}{]}}.

Для моделирования слоистого композита вводятся переменные \(i^{x} \in \left\{ 1,2,\ldots,i_{f}^{x} \right\},\) \(i^{y} \in \left\{ 1,2,\ldots,i_{f}^{y} \right\},\) где \(i^{x}\) -- порядковый номер конечного элемента по оси \emph{x}, \(i_{f}^{x}\) -- количество конечных элементов по оси \emph{x}, \(i^{y}\) -- порядковый номер конечного элемента по оси \emph{y} (слоя композита), \(i_{f}^{y}\) -- количество слоев. Матрица \emph{A} разбивается на подматрицы \(A^{pl}(i^{y}) = A_{m,i^{pl}},\) где \(i^{pl} \in \left\{ i_{1}^{pl},i_{1}^{pl} + 1,\ldots,i^{y}i_{f}^{x} \right\},\) \(i_{1}^{pl} = \left( i^{y} - 1 \right)i_{f}^{x} + 1.\) Конечные элементы отдельного слоя содержат общие значения параметров, вследствие чего уточняются выражения  (\ref{eq:44}) и  (\ref{eq:45}):

\begin{equation}
\begin{matrix}
k^{E}(i^{y})_{r,s} = \frac{E_{2}(i^{y})h}{4\left( 1 + \mu_{1}(i^{y}) \right)\left( 1 - \mu_{1}(i^{y}) - 2n_{1}(i^{y})\mu_{2}(i^{y})^{2} \right)} \times \\
 \times \begin{bmatrix}
\left( n_{1}(i^{y}) - n_{1}(i^{y})^{2}\mu_{2}(i^{y})^{2} \right)\gamma(i^{y})\xi_{r}\xi_{s}\left( 1 + \frac{\eta_{r}\eta_{s}}{3} \right) & n_{1}(i^{y})\mu_{2}(i^{y})\left( 1 + \mu_{1}(i^{y}) \right)\xi_{r}\eta_{s} \\
n_{1}(i^{y})\mu_{2}(i^{y})\left( 1 + \mu_{1}(i^{y}) \right)\eta_{r}\xi_{s} & \left( 1 - \mu_{1}(i^{y})^{2} \right)\frac{\eta_{r}\eta_{s}}{\gamma(i^{y})}\left( 1 + \frac{\xi_{r}\xi_{s}}{3} \right)
\end{bmatrix},
\end{matrix}
\end{equation}

\begin{equation}
k^{G}(i^{y})_{r,s} = \frac{G_{2}(i^{y})h}{4}\begin{bmatrix}
\frac{\eta_{r}\eta_{s}}{\gamma(i^{y})}\left( 1 + \frac{\xi_{r}\xi_{s}}{3} \right) & \eta_{r}\xi_{s} \\
\xi_{r}\eta_{s} & \gamma(i^{y})\xi_{r}\xi_{s}\left( 1 + \frac{\eta_{r}\eta_{s}}{3} \right)
\end{bmatrix},
\end{equation}

\noindent\(\gamma(i^{y}) = b_{fe}(i^{y})/a_{fe}.\)\par
Расширенная матрица жесткости конечного элемента \(k^{\exp}\) строится по принципу: \(k^{\exp}(i^{x},i^{y})_{m,n} = k_{r,s}^{e},\) где \(r = A^{pl}(i^{y})_{m,i^{x}},\) \(s = A^{pl}(i^{y})_{n,i^{x}},\) \(\left( m,n \in \left\{ 1,2,\ldots,m_{f} \right\} \right).\) Если \(r \vee s = 0,\) то \(k_{r,s}^{e} = \begin{pmatrix}
0 & 0 \\
0 & 0
\end{pmatrix}.\)

Матрица жесткости конечно-элементной модели \emph{K} определяется суммированием: \(K = \sum_{i^{y}}^{}{K^{pl}(i^{y})},\) \(K^{pl}(i^{y}) = \sum_{i^{x}}^{}{k^{\exp}(i^{x},i^{y})}.\)

Вектор перемещений \emph{v} узлов конечного элемента строится по принципу: \(v(i^{x},i^{y}) = \begin{pmatrix}
u_{m^{1}}^{x} & u_{m^{1}}^{y} & u_{m^{2}}^{x} & \begin{matrix}
u_{m^{2}}^{y} & u_{m^{3}}^{x} & u_{m^{3}}^{y} & \begin{matrix}
u_{m^{4}}^{x} & u_{m^{4}}^{y}
\end{matrix}
\end{matrix}
\end{pmatrix}^{T},\) где \(A^{pl}(i^{y})_{m^{1},i^{x}} = 1,\) \(A^{pl}(i^{y})_{m^{2},i^{x}} = 2,\) \(A^{pl}(i^{y})_{m^{3},i^{x}} = 3,\) \(A^{pl}(i^{y})_{m^{4},i^{x}} = 4,\) \(u_{m^{q}}^{x} = v_{q}^{x}\) и \(u_{m^{q}}^{y} = v_{q}^{y}.\)

Вектор деформаций \(\varepsilon\) конечного элемента определяется из выражения:

\begin{equation}
\varepsilon(i^{x},i^{y},\xi,\eta) = \beta(i^{y},\xi,\eta)v(i^{x},i^{y}),
\end{equation}

\noindent где \(\beta(i^{y},\xi,\eta)\) -- матрица связи между узловыми перемещениями и деформациями:

\begin{equation}
\beta(i^{y},\xi,\eta) = \frac{1}{2}\begin{pmatrix}
b_{a}(1,\eta) & 0 & a_{a}(i^{y},1,\xi) \\
0 & a_{a}(i^{y},1,\xi) & b_{a}(1,\eta) \\
b_{a}(2,\eta) & 0 & a_{a}(i^{y},2,\xi) \\
0 & a_{a}(i^{y},2,\xi) & b_{a}(2,\eta) \\
b_{a}(3,\eta) & 0 & a_{a}(i^{y},3,\xi) \\
0 & a_{a}(i^{y},3,\xi) & b_{a}(3,\eta) \\
b_{a}(4,\eta) & 0 & a_{a}(i^{y},4,\xi) \\
0 & a_{a}(i^{y},4,\xi) & b_{a}(4,\eta)
\end{pmatrix}^{T},
\end{equation}

\noindent где \(a_{a}(i^{y},q,\xi) = \eta_{q}\left( 1 + \xi_{q}\xi \right)/b_{fe}(i^{y}),\) \(b_{a}(q,\eta) = \xi_{q}\left( 1 + \eta_{q}\eta \right)/a_{fe}.\)

Вектор узловых напряжений \(\sigma\) в вершинах конечного элемента определяется из соотношения:

\begin{equation}
\sigma(i^{x},i^{y},\xi,\eta) = \chi(i^{y})\varepsilon(i^{x},i^{y},\xi,\eta),
\end{equation}

\noindent где \(\chi(i^{y})\) -- матрица упругих постоянных  (\ref{eq:41}), \(\xi = \xi_{q},\) \(\eta = \eta_{q}.\)

Главные напряжения \(\sigma_{1}\) и \(\sigma_{2}\) в плоскости расчета \emph{xy} (рисунок \ref{fig:3}) также выражаются через соотношение  (\ref{eq:35}), а главное напряжение из плоскости для трансверсально-изотропного материала \(\sigma_{3}^{ti}\) определяется выражением:

\begin{equation}
\sigma_{3}^{ti} = \sigma_{z}^{ti} = \mu_{1}\sigma_{x} + \mu_{2}\sigma_{y}\frac{E_{1}}{E_{2}}.
\end{equation}

С целью верификации конечно-элементного моделирования композитных пластин был применен дополнительный алгоритм решения плоской задачи, в котором сотовый заполнитель моделируется аналогично с использованием вышеприведенных выражений, а сплошные слои моделируются несовместными конечными элементами \dualcite{ref18,ref21,ref77,ref78}{{[}\bibnum{ref18}, \bibnum{ref21}, \bibnum{ref77}, \bibnum{ref78}{]}}{{[}\bibnum{ref18}, \bibnum{ref21}, \bibnum{ref77}, \bibnum{ref78}{]}}. При этом подматрицы нормальных и сдвиговых деформаций матрицы жесткости несовместного конечного элемента определяются выражениями:

\begingroup
\setlength{\arraycolsep}{2pt}
\begin{equation}
\begin{matrix}
k^{inc.E}(i^{y})_{r,s} = \frac{E(i^{y})h}{4\left( 1 + \mu(i^{y}) \right)\left( 1 - 2\mu(i^{y}) \right)} \times \\
 \times \begin{bmatrix}
\begin{gathered}
\gamma(i^{y})\xi_{r}\xi_{s}\times\\[-2pt]
{}\times\left(\begin{gathered}
1 - \mu(i^{y}) + \\[-2pt]
{}+\frac{1 - \mu(i^{y}) - \mu(i^{y})^{2} - \mu(i^{y})^{3}}{3}\eta_{r}\eta_{s}
\end{gathered}\right)
\end{gathered}
& \mu(i^{y})\xi_{r}\eta_{s} \\
\mu(i^{y})\eta_{r}\xi_{s}
& \begin{gathered}
\frac{\eta_{r}\eta_{s}}{\gamma(i^{y})}\times\\[-2pt]
{}\times\left(\begin{gathered}
1 - \mu(i^{y}) + \\[-2pt]
{}+\frac{1 - \mu(i^{y}) - \mu(i^{y})^{2} - \mu(i^{y})^{3}}{3}\xi_{r}\xi_{s}
\end{gathered}\right)
\end{gathered}
\end{bmatrix},
\end{matrix}
\end{equation}
\endgroup

\begin{equation}
k^{inc.G}(i^{y})_{r,s} = \frac{G(i^{y})h}{4}\begin{bmatrix}
\frac{\eta_{r}\eta_{s}}{\gamma(i^{y})} & \eta_{r}\xi_{s} \\
\xi_{r}\eta_{s} & \gamma(i^{y})\xi_{r}\xi_{s}
\end{bmatrix},
\end{equation}

\noindent а матрица связи между узловыми перемещениями и деформациями описывается выражением:

\begin{equation}
\beta^{inc}(i^{y},\xi,\eta) = \frac{1}{2}\begin{pmatrix}
b_{a}(1,\eta) & e_{a}(i^{y},1,\eta) & g_{a}(i^{y},1) \\
c_{a}(i^{y},1,\xi) & a_{a}(i^{y},1,\xi) & f_{a}(1) \\
b_{a}(2,\eta) & e_{a}(i^{y},2,\eta) & g_{a}(i^{y},2) \\
c_{a}(i^{y},2,\xi) & a_{a}(i^{y},2,\xi) & f_{a}(2) \\
b_{a}(3,\eta) & e_{a}(i^{y},3,\eta) & g_{a}(i^{y},3) \\
c_{a}(i^{y},3,\xi) & a_{a}(i^{y},3,\xi) & f_{a}(3) \\
b_{a}(4,\eta) & e_{a}(i^{y},4,\eta) & g_{a}(i^{y},4) \\
c_{a}(i^{y},4,\xi) & a_{a}(i^{y},4,\xi) & f_{a}(4)
\end{pmatrix}^{T},
\end{equation}

\noindent где \(c_{a}(i^{y},q,\xi) = - \mu(i^{y})\xi_{q}\eta_{q}\xi/b_{fe}(i^{y}),\) \(e_{a}(i^{y},q,\eta) = - \mu(i^{y})\xi_{q}\eta_{q}\eta/a_{fe},\) \(g_{a}(i^{y},q) = \eta_{q}/b_{fe}(i^{y}).\)

При анализе композитных пластин с помощью вышеприведенных алгоритмов решения плоской задачи слои пластин разбивались вдоль оси \emph{x} (рисунок \ref{fig:3}) на 36 конечных элементов, а внешние слои содержали вдоль оси \emph{y} 1 слой элементов. При этом в первой постановке экспериментов заполнитель разбивался вдоль оси \emph{y} на 1 слой элементов, а во второй постановке при \(0.7 \leq t_{cl} \leq 1.4\) мм и \(1.7 \leq t_{cl} \leq 3.6\) мм -- на 1 и 2 слоя конечных элементов соответственно.

\clearpage
\iflatexml
\section{ГЛАВА 3. РЕЗУЛЬТАТЫ МАТЕМАТИЧЕСКОГО МОДЕЛИРОВАНИЯ И ЛАБОРАТОРНЫХ ИСПЫТАНИЙ ТРЕХСЛОЙНЫХ ПЛАСТИН В УСЛОВИЯХ ИЗГИБА}
\else
\phantomsection
\setcounter{section}{2}\refstepcounter{section}
\section{ГЛАВА 3. РЕЗУЛЬТАТЫ МАТЕМАТИЧЕСКОГО МОДЕЛИРОВАНИЯ И ЛАБОРАТОРНЫХ ИСПЫТАНИЙ ТРЕХСЛОЙНЫХ ПЛАСТИН В УСЛОВИЯХ ИЗГИБА}
\fi
\label{sec:chapter3}

\begin{figure}[H]
\centering
\includegraphics[keepaspectratio,width=6.49606in,height=3.8577in]{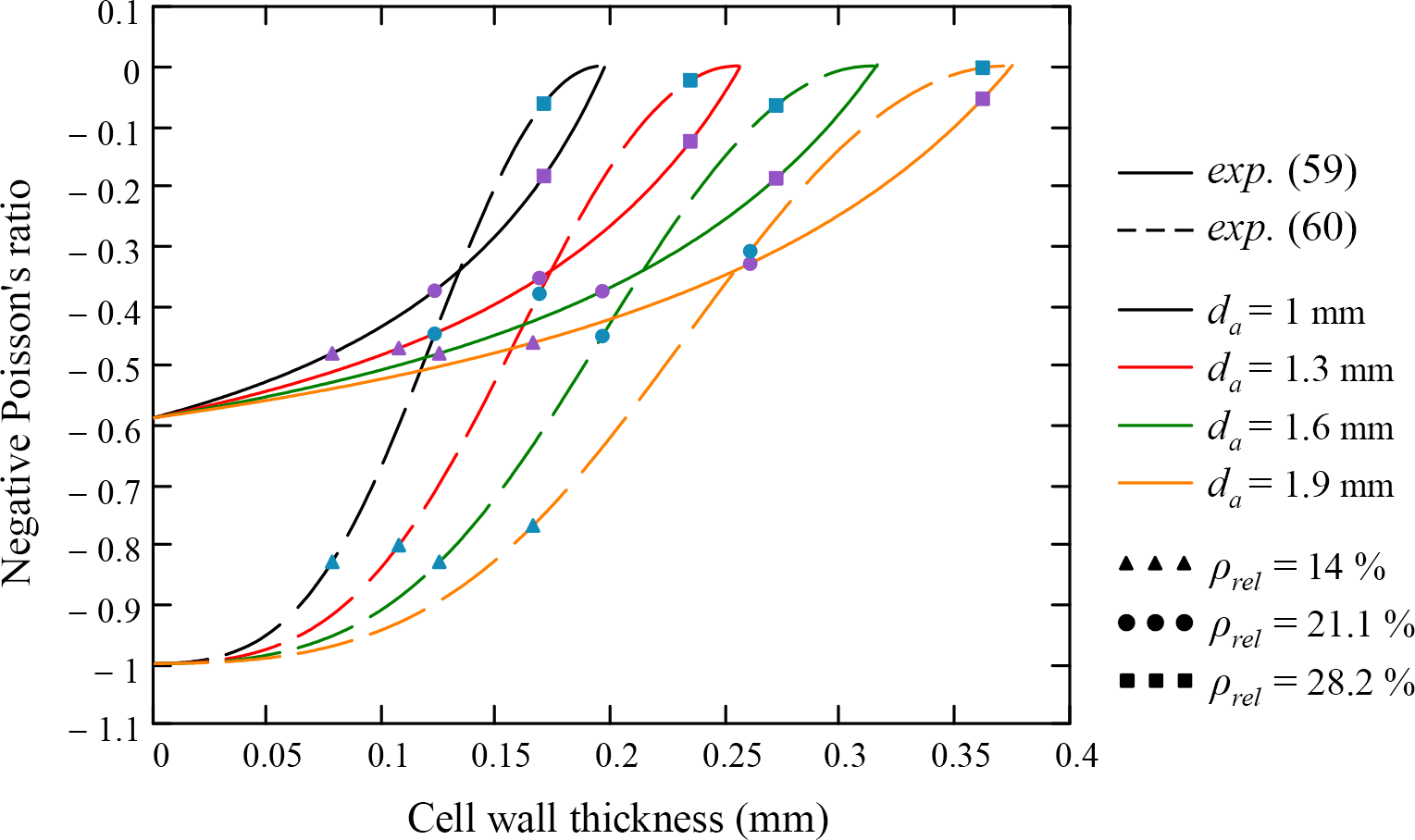}
\caption{Диаграмма изменения ОКП от толщины стенок элементарных ячеек тетракиральных сот при различной дискретизации \dualcite{ref77}{{[}\bibnum{ref77}{]}}{{[}\bibnum{ref77}{]}}}\label{fig:9}
\end{figure}

Построена диаграмма зависимости ОКП исследуемых тетракиральных сот от толщины стенок элементарных ячеек (рисунок \ref{fig:9}) \dualcite{ref77}{{[}\bibnum{ref77}{]}}{{[}\bibnum{ref77}{]}} по следующим аналитическим выражениям, которые предложены в работах \dualcite{ref85}{{[}\bibnum{ref85}{]}}{{[}\bibnum{ref85}{]}} и \dualcite{ref74}{{[}\bibnum{ref74}{]}}{{[}\bibnum{ref74}{]}} соответственно:

\begin{equation}
\label{eq:59}
\mu = - \frac{\sin\theta\left\lbrack \alpha_{h} - \alpha_{h}(\pi - 2\theta)(\alpha_{h} + \beta_{h}) \right\rbrack}{2\alpha_{h}(\alpha_{h} - \alpha_{h}\sin\theta - \beta_{h}\sin\theta)},
\end{equation}

\noindent \(\alpha_{h} = r_{a}/l,\) \(\beta_{h} = t_{sw}/l,\)

\begin{equation}
\mu = - \frac{(a_{h} - b_{h})^{2}\sin^{2}\theta\cos^{2}\theta}{2a_{h}b_{h}(\sin^{4}\theta + \cos^{4}\theta) + (a_{h} + b_{h})^{2}\sin^{2}\theta\cos^{2}\theta},
\end{equation}

\noindent где \(a_{h} = {l_{e}}^{3}/(24E_{s}I),\) \(b_{h} = l/(2t_{h}t_{sw}E_{s}),\) \(l_{e} = l - 2\sqrt{2rt_{sw} - {t_{sw}}^{2}},\) \(I = t_{h}{t_{sw}}^{3}/12,\) \(t_{h}\) -- толщина сотовой структуры.

У рассматриваемых сот при малых значениях относительной плотности (14, 21.1 и 28.2 \%) наблюдается ОКП. При этом у структур с разной дискретизацией, но равной \(\rho_{rel},\) значения ОКП находятся в узком диапазоне. Согласно выражению  (\ref{eq:59}), при увеличении \(\rho_{rel}\) эффективный коэффициент Пуассона переходит из области отрицательных значений в область положительных значений.

\begin{figure}[H]
\centering
\includegraphics[keepaspectratio,width=6.69291in,height=3.89725in]{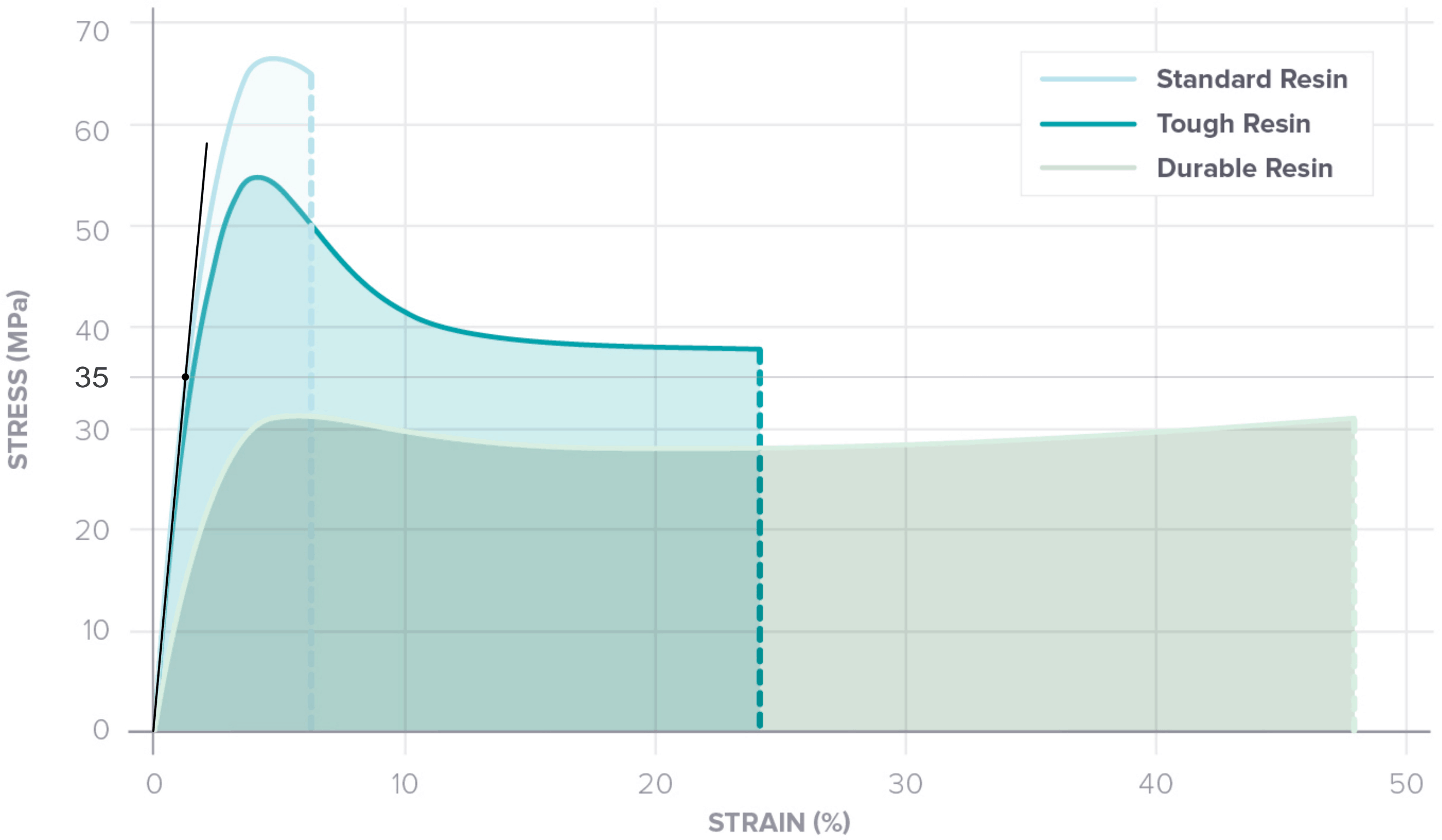}
\caption{Определение предела пропорциональности смолы Formlabs Clear (Standard Resin) \dualcite{ref78,ref97}{{[}\bibnum{ref78}, \bibnum{ref97}{]}}{{[}\bibnum{ref78}, \bibnum{ref97}{]}}}\label{fig:10}
\end{figure}

В первой постановке экспериментов для композитных пластин принимаются следующие значения параметров (рисунок \ref{fig:3}): \(a = 54\) мм, \(h = 13\) мм, \(t_{p} = 2\) мм, \(t_{fl} = 0.5\) мм, \(t_{cl} = 1\) мм, \(l_{1} = a/2,\) \(x_{1} = 12\) мм, \(x_{2} = 42\) мм, а вторая постановка отличается определением толщины сот: \(t_{cl} = t_{cl}(V_{cl},\rho_{rel}),\) \(14 \leq \rho_{rel} \leq 70.9\) \%, \(V_{cl} = 351\) мм\textsuperscript{3} -- объем твердого тела сот \dualcite{ref76,ref77}{{[}\bibnum{ref76}, \bibnum{ref77}{]}}{{[}\bibnum{ref76}, \bibnum{ref77}{]}}. В качестве материала, из которого изготовлен композит, используется полимер Formlabs Clear Resin \dualcite{ref81}{{[}\bibnum{ref81}{]}}{{[}\bibnum{ref81}{]}}, имеющий следующие свойства: модуль упругости при растяжении \(E = 2.8\) ГПа, \(\mu = 0.35,\) плотность \(\rho = 1200\) кг/м\textsuperscript{3}, предел пропорциональности \(\sigma_{pl} \approx 35\) МПа, который определялся по графику (рисунок \ref{fig:10}), представленному компанией Formlabs \dualcite{ref97}{{[}\bibnum{ref97}{]}}{{[}\bibnum{ref97}{]}}.

Исследуемые композиты имеют идеальное сцепление на границе раздела слоев, что позволяет рассматривать их как монолитные конструкции из полимера при возможности послойного анализа напряжений. При вычислении предельного состояния полимеров хорошее согласие с экспериментальными данными может обеспечить критерий Мизеса \dualcite{ref4,ref10}{{[}\bibnum{ref4}, \bibnum{ref10}{]}}{{[}\bibnum{ref4}, \bibnum{ref10}{]}}. Вследствие удобства использования и удовлетворительных результатов данный критерий широко применяется при анализе напряженного состояния полимерных конструкций \dualcite{ref37,ref42,ref54,ref60,ref84,ref103}{{[}\bibnum{ref37}, \bibnum{ref42}, \bibnum{ref54}, \bibnum{ref60}, \bibnum{ref84}, \bibnum{ref103}{]}}{{[}\bibnum{ref37}, \bibnum{ref42}, \bibnum{ref54}, \bibnum{ref60}, \bibnum{ref84}, \bibnum{ref103}{]}}. При этом было принято, что используемый материал Formlabs Clear является изотропным с незначительным различием в сопротивлении растяжению и сжатию.

\begin{figure}[H]
\centering
\includegraphics[keepaspectratio,width=6.69214in,height=2.60764in]{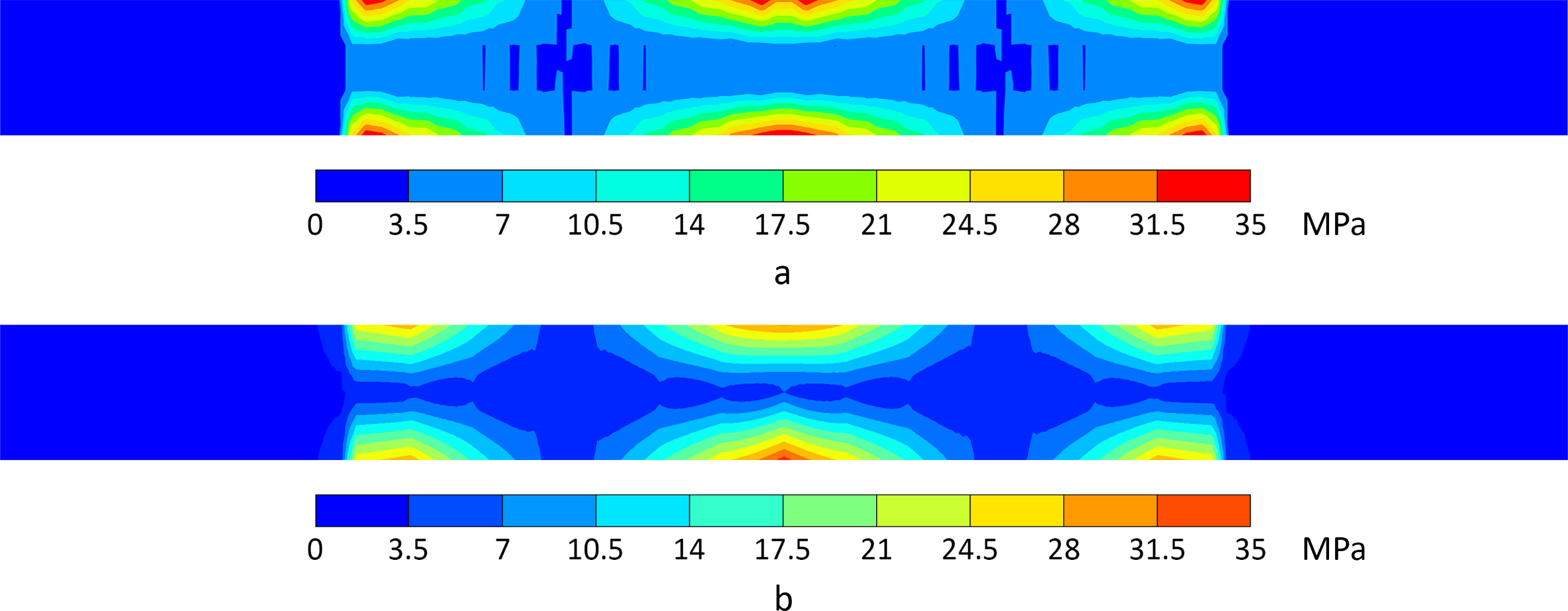}
\caption{Диаграммы распределения напряжений в сплошной пластине в условиях изгиба с жестким защемлением при \(\sigma_{\max} = \sigma_{pl},\) полученные с применением двух подходов: а) посредством алгоритма решения плоской задачи с использованием несовместных конечных элементов; б) используя систему Comsol Multiphysics}\label{fig:11}
\end{figure}

Получены диаграммы распределения напряжений в сплошной пластине с размерами \(54 \times 13 \times 2\) мм в условиях плоского деформированного состояния с помощью модуля «Механика конструкций» системы Comsol Multiphysics и алгоритма решения плоской задачи с использованием несовместных конечных элементов. При моделировании в системе Comsol Multiphysics пластина разбивалась по толщине на один слой прямоугольных конечных элементов, а при использовании алгоритма решения плоской задачи -- на два. По результатам расчета пластины с применением алгоритма в условиях жесткого защемления получен график изолиний с распределением напряжений (рисунок \ref{fig:11}a). В системе Comsol Multiphysics построена аналогичная диаграмма распределения напряжений (рисунок \ref{fig:11}b). Полученное значение \(F_{y}\) для пластины с жестким защемлением при \(\sigma_{\max} = \sigma_{pl}\) с применением алгоритма составляет 103.3~Н, а с применением Comsol Multiphysics -- 91.7 Н.

\begin{figure}[H]
\centering
\includegraphics[keepaspectratio,width=6.69214in,height=2.60764in]{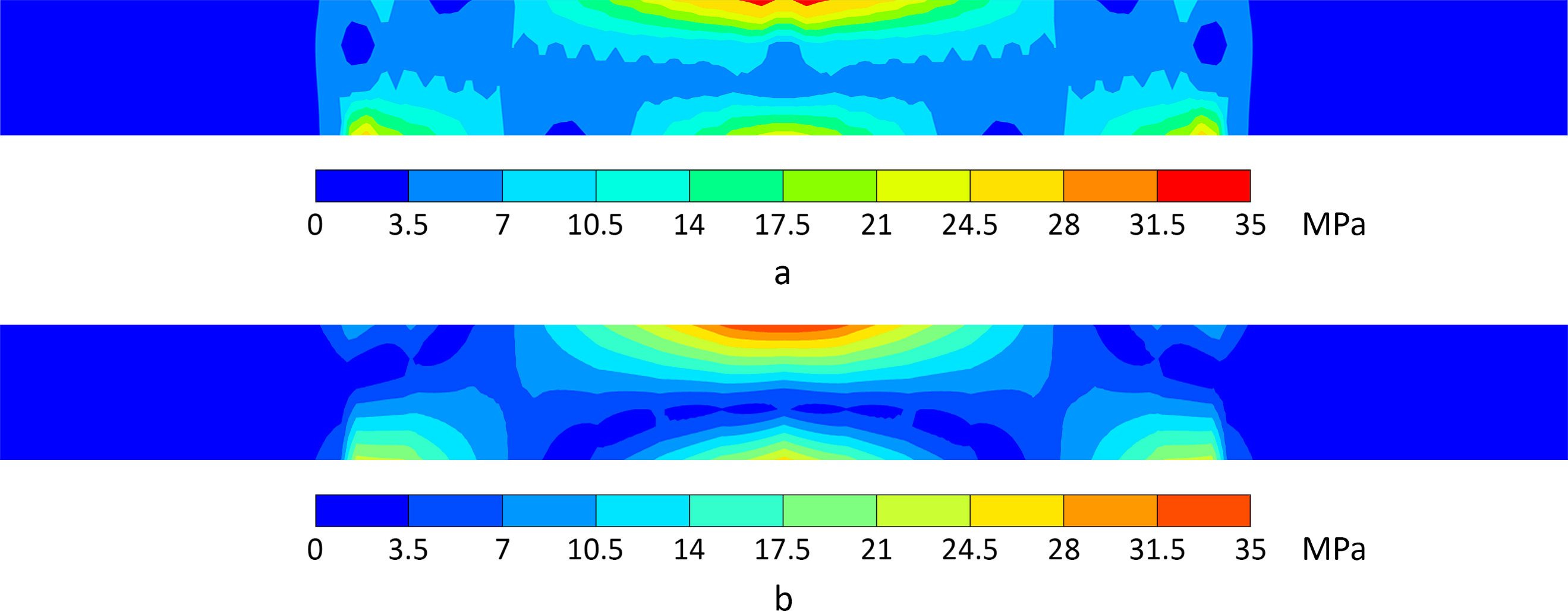}
\caption{Диаграммы распределения напряжений в сплошной пластине в условиях опирания с упругим поворотом при \(\sigma_{\max} = \sigma_{pl},\) полученные с применением двух подходов: а) посредством алгоритма решения плоской задачи с использованием несовместных конечных элементов; б) используя систему Comsol Multiphysics}\label{fig:12}
\end{figure}

Также представлены диаграммы распределения напряжений с применением алгоритма (рисунок \ref{fig:12}a) и системы Comsol Multiphysics (рисунок \ref{fig:12}b) в условиях опирания с упругим поворотом. Полученное значение \(F_{y}\) для пластины в условиях опирания с упругим поворотом при \(\sigma_{\max} = \sigma_{pl}\) с применением алгоритма составляет 67.1 Н, а с применением Comsol Multiphysics -- 67.8 Н. Решения, полученные с применением двух подходов, хорошо согласуются.

\begin{figure}[H]
\centering
\includegraphics[keepaspectratio,width=6.10236in,height=3.80229in]{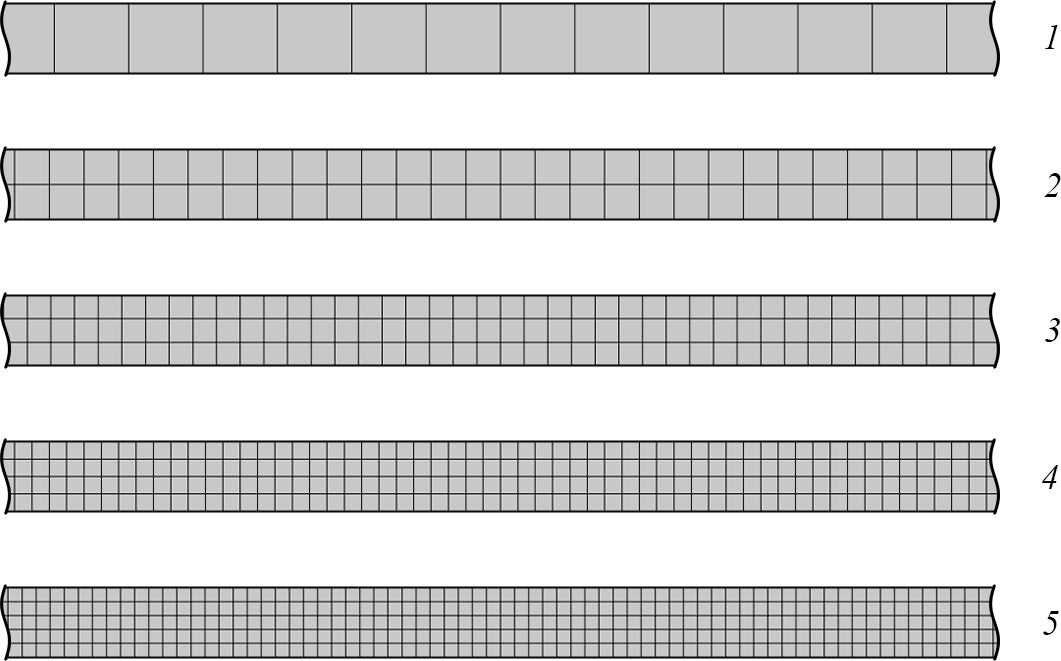}
\caption{Сгущение сетки конечных элементов сплошной пластины}\label{fig:13}
\end{figure}

\begin{figure}[H]
\centering
\includegraphics[keepaspectratio,width=5.71316in,height=8.08661in]{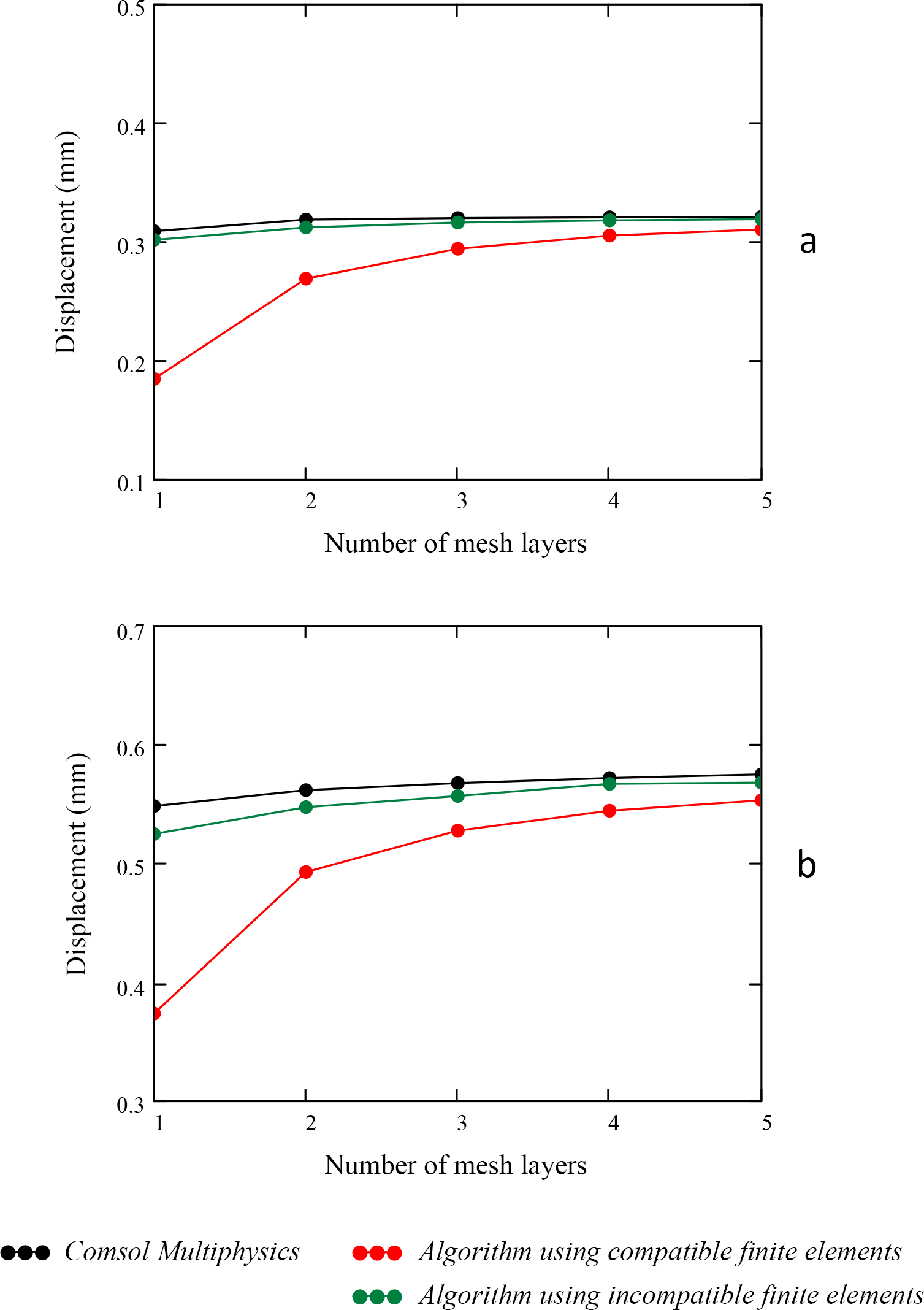}
\caption{Диаграмма зависимости \(u_{\max}\) в сплошной пластине от густоты сетки конечных элементов при \(F_{y} = 60\) Н в условиях жесткого защемления (a) и опирания с упругим поворотом (b)}\label{fig:14}
\end{figure}

\begin{figure}[H]
\centering
\includegraphics[keepaspectratio,width=5.70263in,height=8.08661in]{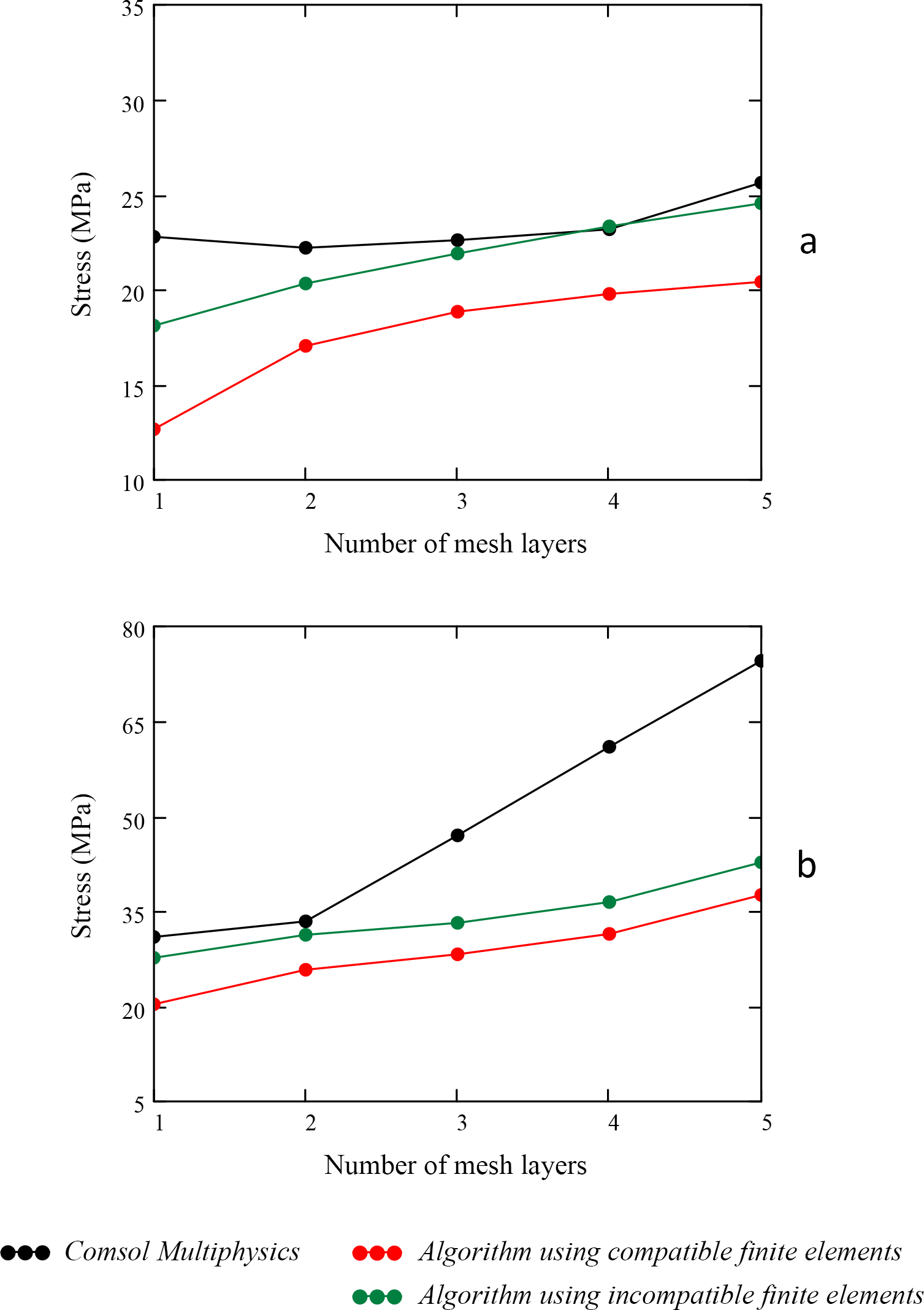}
\caption{Диаграмма зависимости \(\sigma_{\max}\) в сплошной пластине от густоты сетки конечных элементов при \(F_{y} = 60\) Н в условиях жесткого защемления (a) и опирания с упругим поворотом (b)}\label{fig:15}
\end{figure}

Произведена оценка влияния густоты сетки конечных элементов сплошной пластины на изменение в ней максимальных напряжений при \(F_{y} = 60\) Н в условиях изгиба с жестким защемлением и при опирании с упругим поворотом. Для решения задачи плоской деформации использовался модуль «Механика конструкций» системы Comsol Multiphysics, а также разработанные выше алгоритмы для анализа напряженно-деформированного состояния сплошных пластин с использованием совместных и несовместных элементов. При этом сетка конечных элементов содержала от 1 до 5 слоев прямоугольных элементов с приблизительно равными сторонами (рисунок \ref{fig:13}). Диаграммы зависимости \(u_{\max}\) и \(\sigma_{\max}\) в сплошной пластине от густоты сетки конечных элементов представлены на рисунках \ref{fig:14} и \ref{fig:15} соответственно.

Следует отметить, что в случае решения данной задачи в трехмерной постановке с использованием четырехугольных призм (рисунок \ref{fig:4}) с приблизительно равными ребрами при условии аналогичного сгущения сетки в продольном сечении пластины (рисунок \ref{fig:13}) и ограничения перемещений узлов \(u_{z} = 0\) (рисунок \ref{fig:3}) численные результаты будут соответствовать результатам решения плоской задачи.

\vspace{-0.575\baselineskip}

\begin{figure}[H]
\centering
\includegraphics[keepaspectratio,width=4.80315in,height=6.54893in]{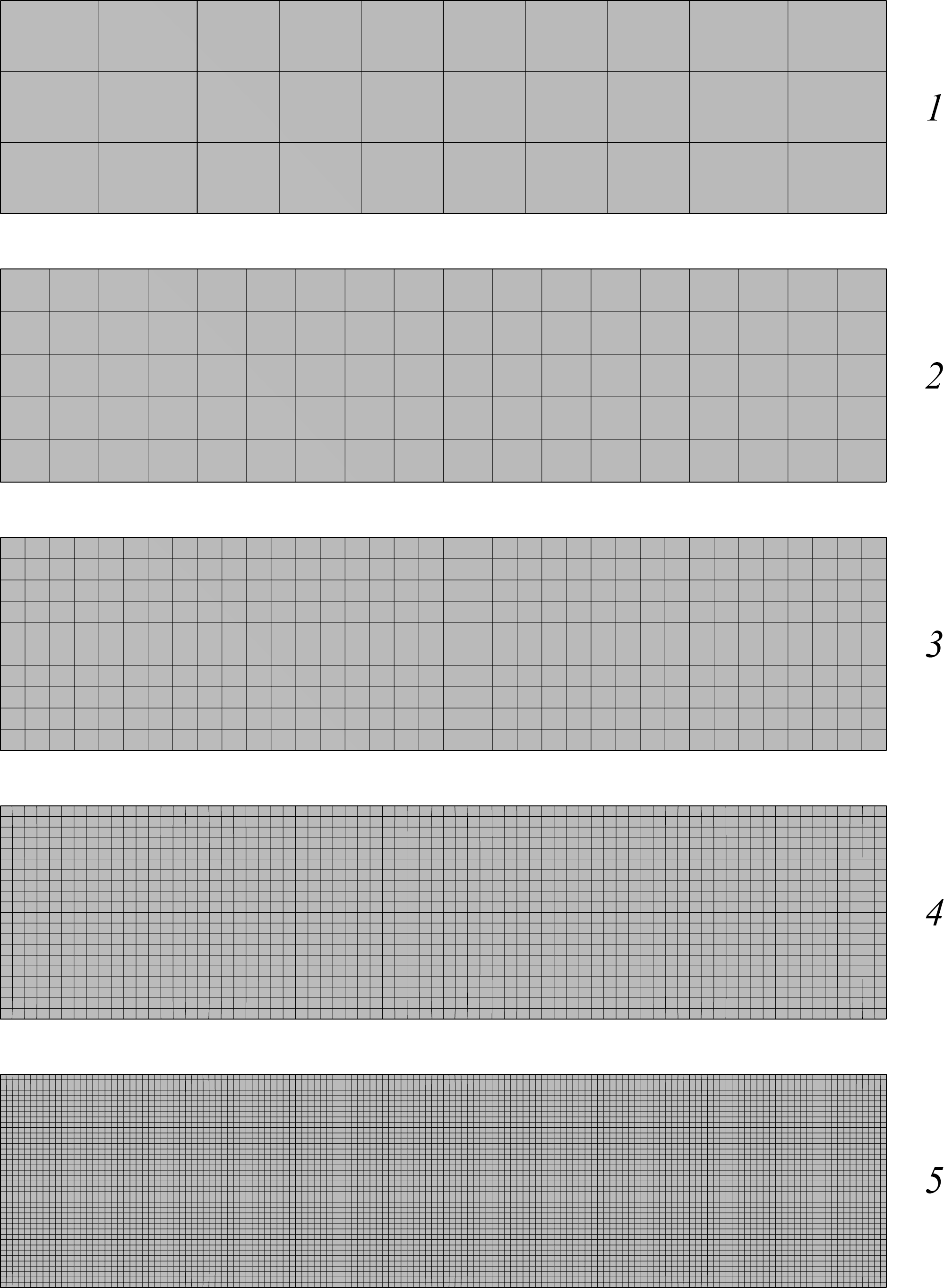}
\caption{Сгущение сетки конечных элементов сплошных слоев композитной пластины}\label{fig:16}
\end{figure}

\begin{figure}[H]
\centering
\includegraphics[keepaspectratio,width=4.67388in,height=7.87402in]{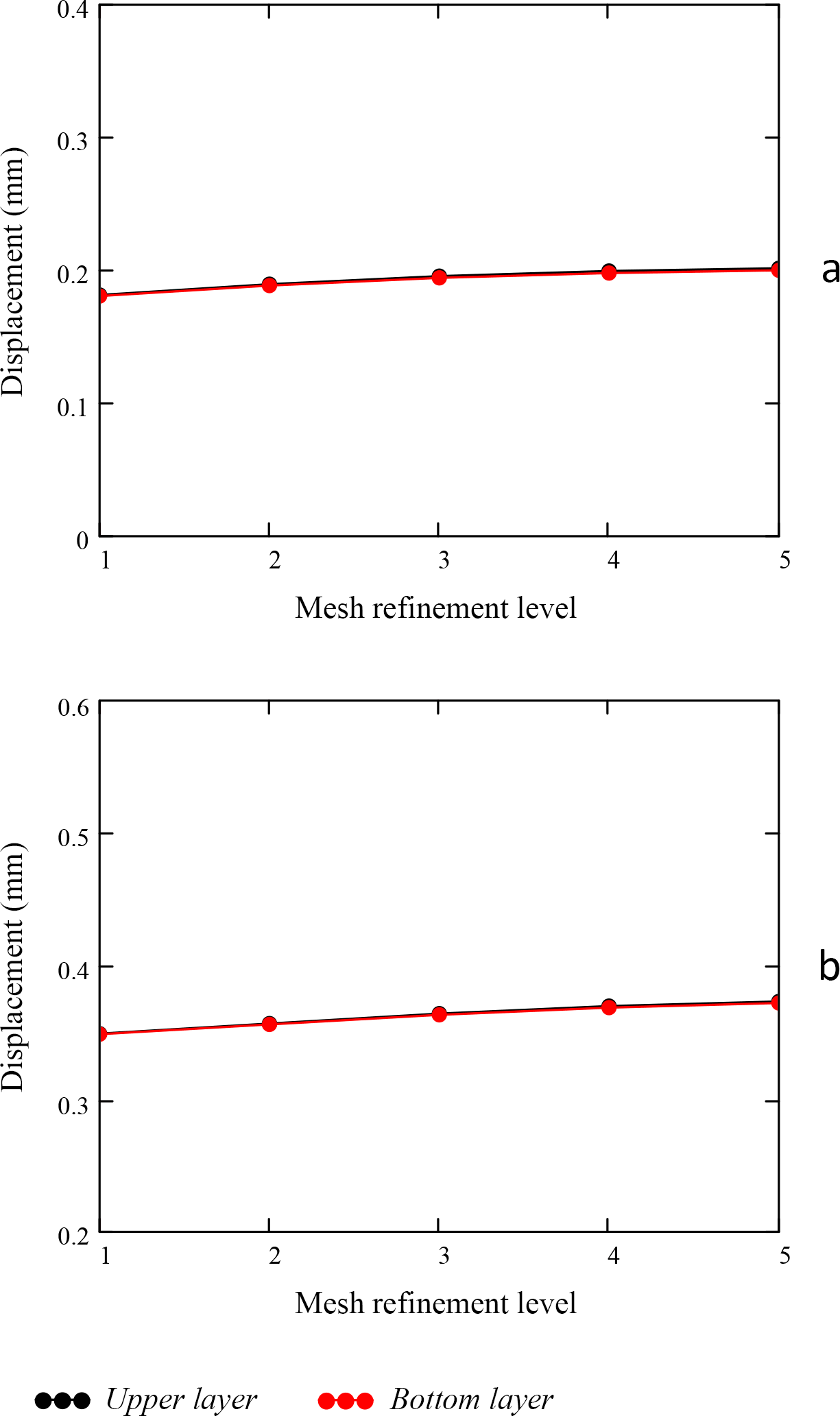}
\caption{Диаграмма зависимости \(u_{\max}\) в сплошных слоях композита от густоты сетки конечных элементов при~\(d_{a} = 1.6\) мм, \(\rho_{rel} = \text{42.5}\) \%, \(F_{y} = 30\) Н в условиях жесткого защемления (a) и опирания с упругим поворотом (b)}\label{fig:17}
\end{figure}

\begin{figure}[H]
\centering
\includegraphics[keepaspectratio,width=4.63896in,height=7.87402in]{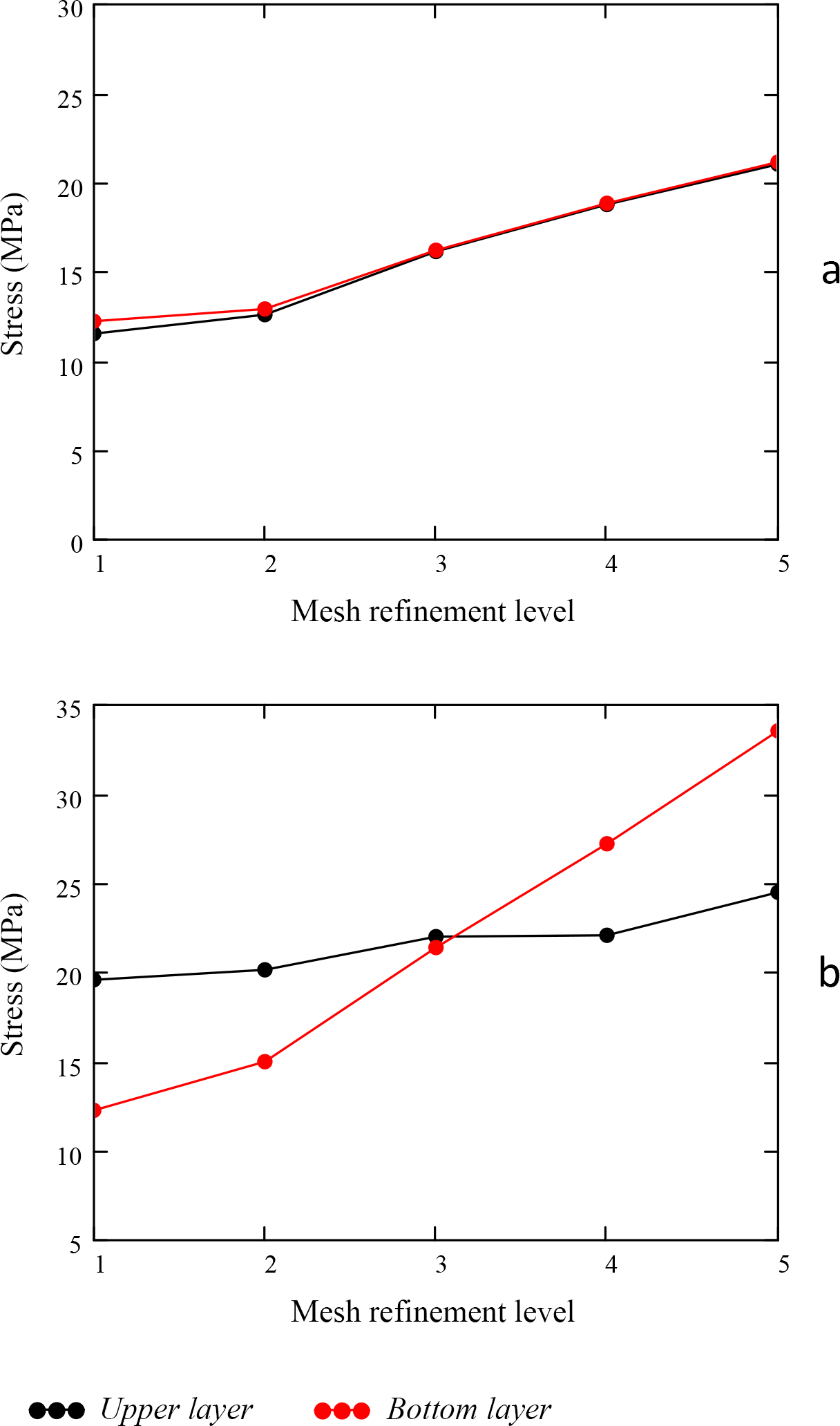}
\caption{Диаграмма зависимости \(\sigma_{\max}\) в сплошных слоях композита от густоты сетки конечных элементов при~\(d_{a} = 1.6\) мм, \(\rho_{rel} = \text{42.5}\) \%, \(F_{y} = 30\) Н в условиях жесткого защемления (a) и опирания с упругим поворотом (b)}\label{fig:18}
\end{figure}

\clearpage
Также в трехмерной постановке произведена оценка влияния густоты сетки конечных элементов внешних слоев композита (рисунок \ref{fig:16}) на изменение в них максимальных перемещений (рисунок \ref{fig:17}) и напряжений (рисунок \ref{fig:18}) в первой постановке экспериментов при \(d_{a} = 1.6\) мм, \(\rho_{rel} = \text{42.5}\) \%, \(F_{y} = 30\) Н в условиях изгиба с жестким защемлением и при опирании с упругим поворотом.

\vspace{-0.575\baselineskip}
\vspace{-0.575\baselineskip}

\clearpage
В области приложения нагрузки, а также в местах закрепления или опирания возникают локальные концентрации напряжений, величина которых зависит от способа идеализации граничных условий и степени сгущения конечно-элементной сетки. Поэтому определяемые в этих областях максимальные эквивалентные напряжения не рассматриваются как строгие сеточно-инвариантные локальные характеристики напряженного состояния. В настоящей работе они используются преимущественно в качестве сравнительного показателя для анализа качественных закономерностей изменения напряженного состояния при варьировании геометрических параметров.

При моделировании композитных пластин трехмерными конечными элементами для внешних слоев использовалась сетка 2 (рисунок \ref{fig:16}).

\iflatexml\else\phantomsection\fi
\subsection{3.1. Результаты численного моделирования композитных пластин в первой постановке численных экспериментов при постоянной толщине слоев и варьируемой относительной плотности заполнителя}

\begin{figure}[H]
\centering
\includegraphics[keepaspectratio,width=5.0828in,height=7.37782in]{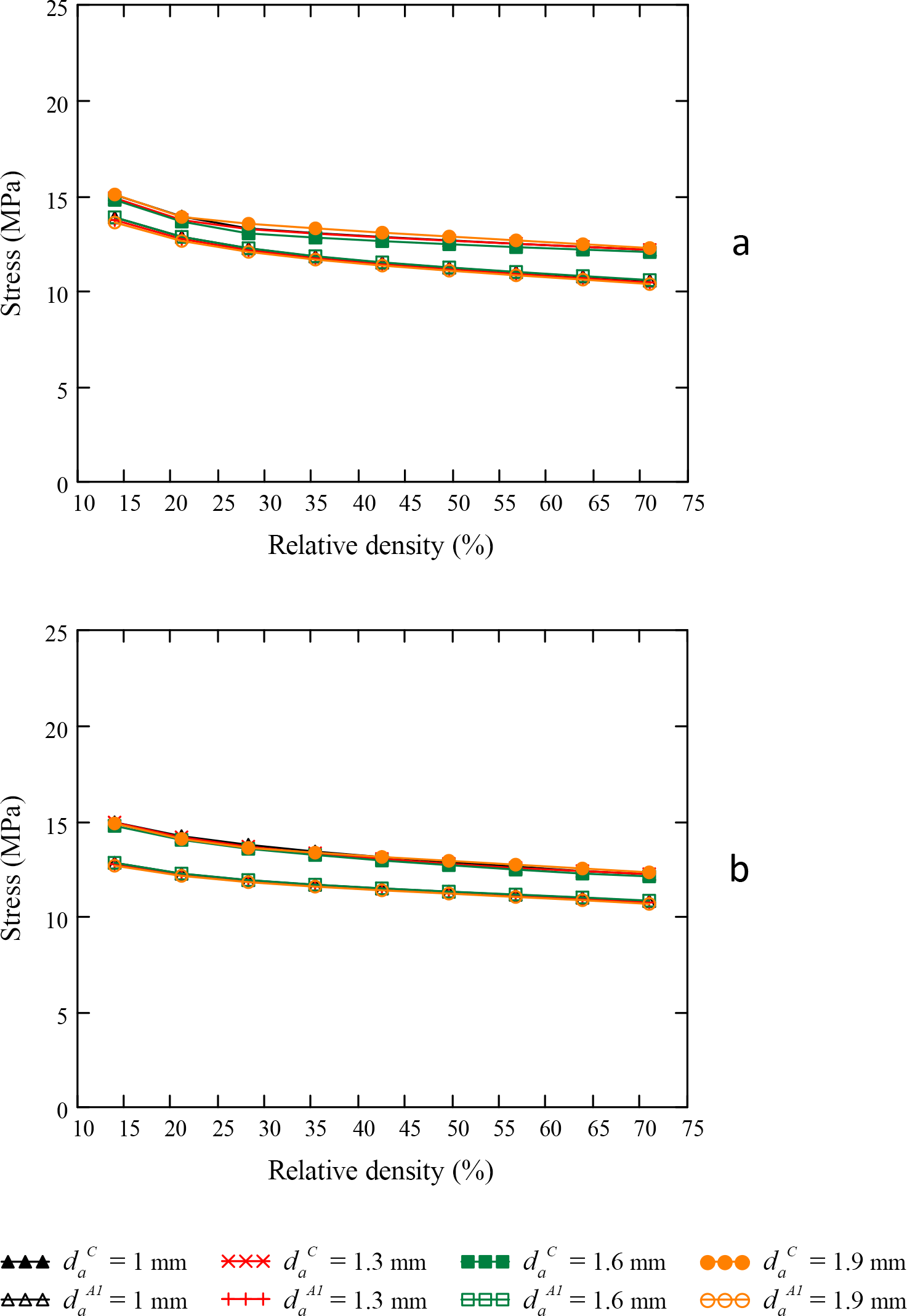}
\caption{Диаграммы зависимости \(\sigma_{\max}\) в верхнем (a) и нижнем (b) сплошном слое композитных пластин от \(\rho_{rel}\) при \(F_{y} = 30\) Н в первой постановке численных экспериментов при жестком защемлении, полученные с применением Comsol Multiphysics \dualcite{ref77,ref78}{{[}\bibnum{ref77}, \bibnum{ref78}{]}}{{[}\bibnum{ref77}, \bibnum{ref78}{]}} и Алгоритма 1}\label{fig:19}
\end{figure}

\begin{figure}[H]
\centering
\includegraphics[keepaspectratio,width=5.47153in,height=7.94541in]{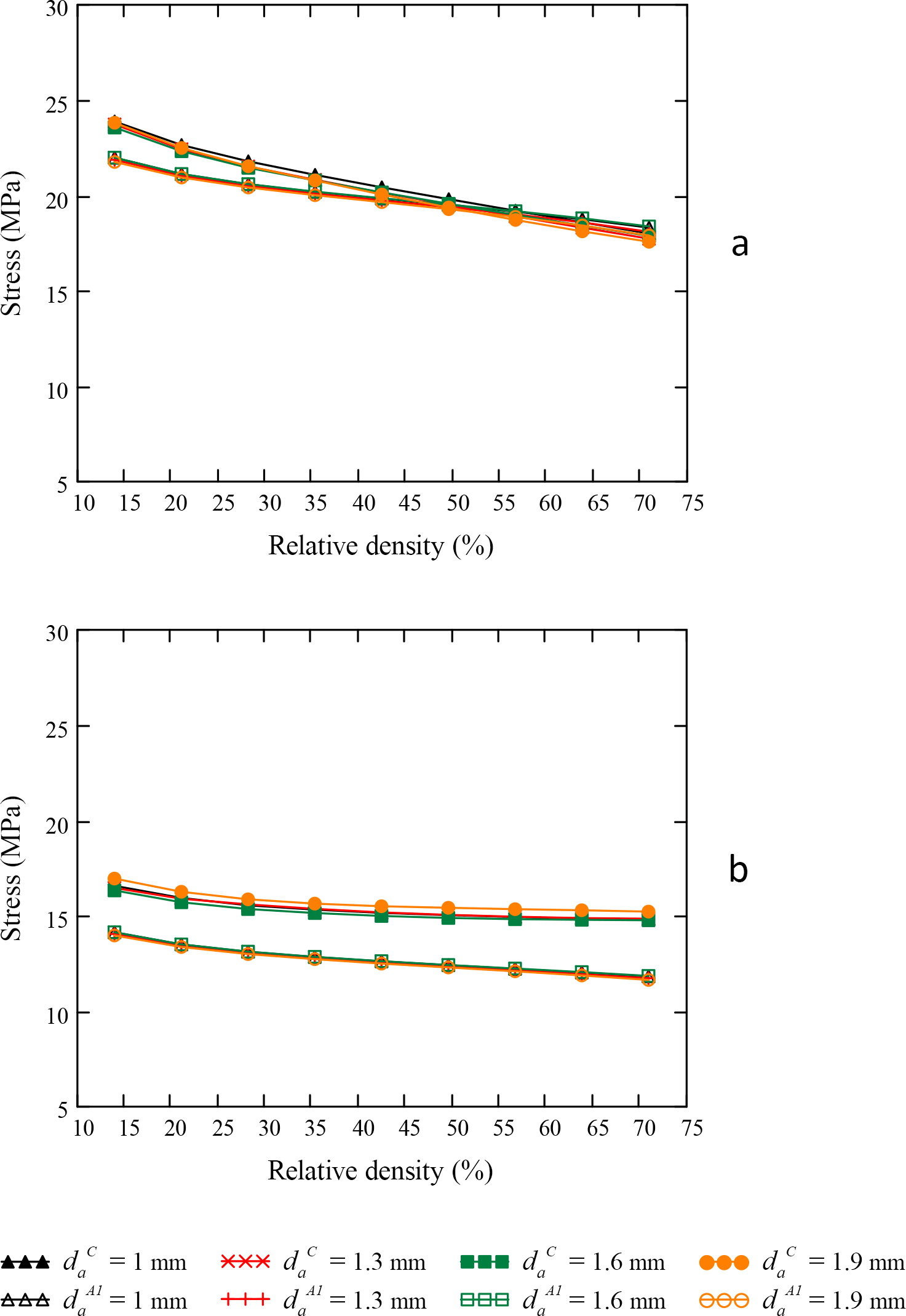}
\caption{Диаграммы зависимости \(\sigma_{\max}\) в верхнем (a) и нижнем (b) сплошном слое композитных пластин от \(\rho_{rel}\) при \(F_{y} = 30\) Н в первой постановке численных экспериментов при опирании с упругим поворотом, полученные с применением Comsol Multiphysics \dualcite{ref78}{{[}\bibnum{ref78}{]}}{{[}\bibnum{ref78}{]}} и Алгоритма 1}\label{fig:20}
\end{figure}

\begin{figure}[H]
\centering
\includegraphics[keepaspectratio,width=5.46643in,height=7.9343in]{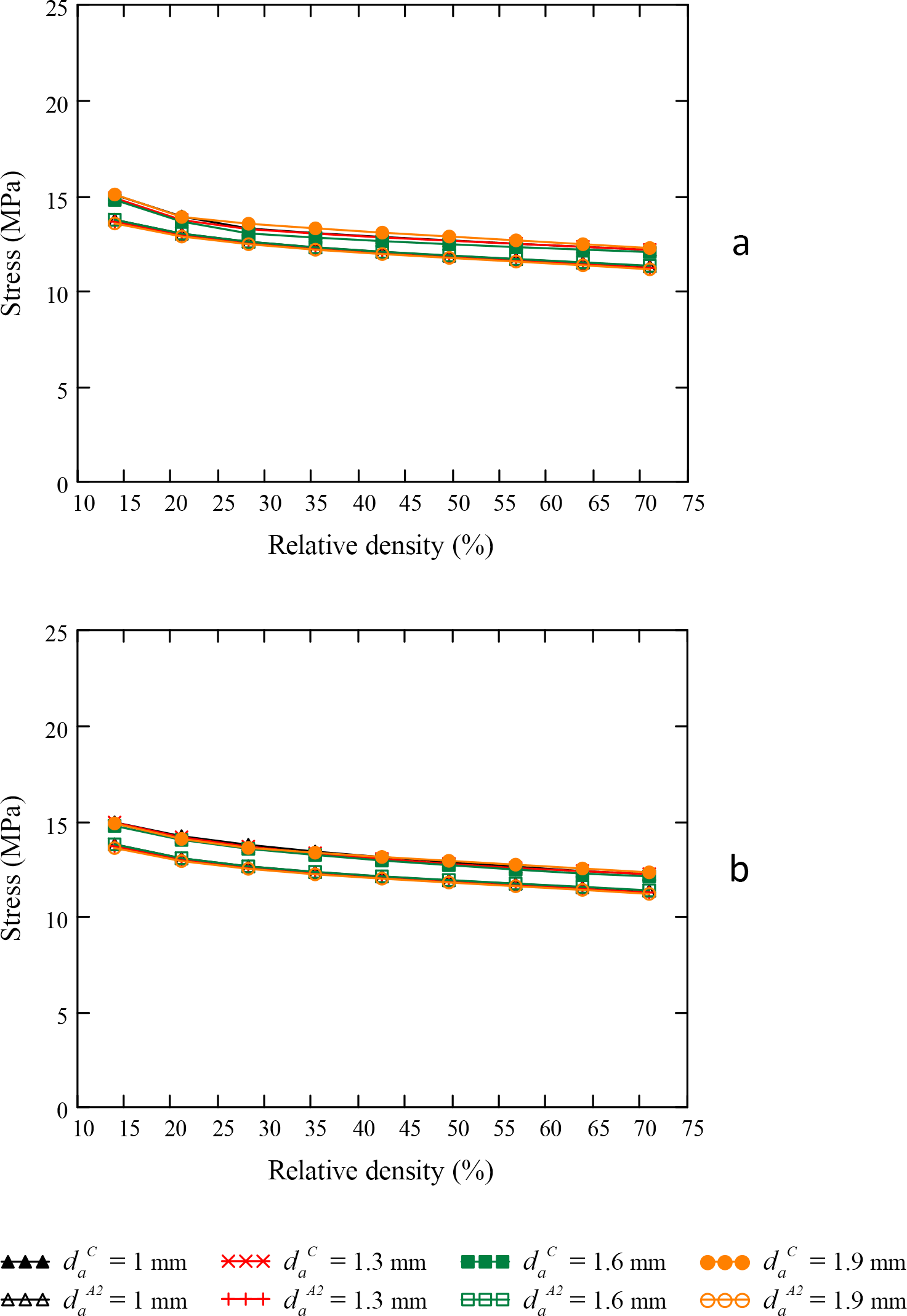}
\caption{Диаграммы зависимости \(\sigma_{\max}\) в верхнем (a) и нижнем (b) сплошном слое композитных пластин от \(\rho_{rel}\) при \(F_{y} = 30\) Н в первой постановке численных экспериментов при жестком защемлении, полученные с применением Comsol Multiphysics \dualcite{ref77,ref78}{{[}\bibnum{ref77}, \bibnum{ref78}{]}}{{[}\bibnum{ref77}, \bibnum{ref78}{]}} и Алгоритма 2}\label{fig:21}
\end{figure}

\begin{figure}[H]
\centering
\includegraphics[keepaspectratio,width=5.46388in,height=7.9343in]{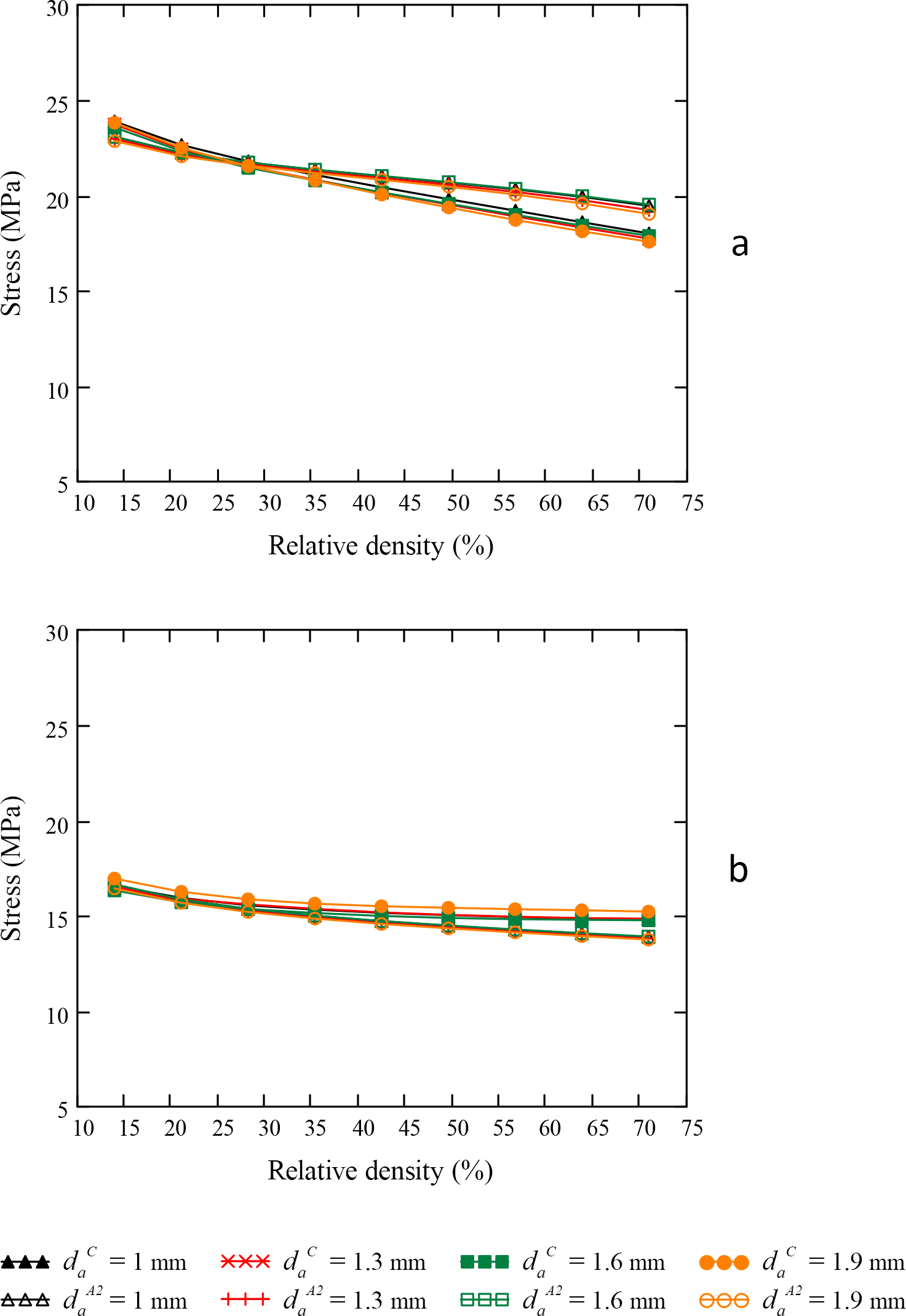}
\caption{Диаграммы зависимости \(\sigma_{\max}\) в верхнем (a) и нижнем (b) сплошном слое композитных пластин от \(\rho_{rel}\) при \(F_{y} = 30\) Н в первой постановке численных экспериментов при опирании с упругим поворотом, полученные с применением Comsol Multiphysics \dualcite{ref78}{{[}\bibnum{ref78}{]}}{{[}\bibnum{ref78}{]}} и Алгоритма 2}\label{fig:22}
\end{figure}

Диаграммы зависимости \(\sigma_{\max}\) во внешних слоях композитных пластин от \(\rho_{rel}\) заполнителя при \(F_{y} = 30\) Н в первой постановке численных экспериментов при жестком защемлении и опирании с упругим поворотом представлены на рисунках \ref{fig:19} и \ref{fig:20} соответственно. Данные диаграммы, помимо результатов моделирования в Comsol Multiphysics \(\left( {d_{a}}^{C} \right),\) содержат результаты анализа композитов с помощью алгоритма решения плоской задачи, в котором сплошные слои моделируются с применением совместных конечных элементов, далее будем называть его «Алгоритм 1» \(\left( {d_{a}}^{A1} \right).\) На рисунках \ref{fig:21} и \ref{fig:22} представлены диаграммы с аналогичной зависимостью в тех же граничных условиях, однако с использованием дополнительного алгоритма решения плоской задачи, в котором сплошные слои моделируются с применением несовместных конечных элементов, далее будем называть его «Алгоритм~2» \(\left( {d_{a}}^{A2} \right).\) Графики, полученные с применением Comsol Multiphysics и Алгоритма~2, хорошо согласуются (рисунки \ref{fig:21}, \ref{fig:22}), тогда как результаты Алгоритма 1 приводят к некоторой погрешности относительно данных Comsol Multiphysics (рисунки \ref{fig:19}, \ref{fig:20}).

\vspace{-0.575\baselineskip}
\vspace{-0.575\baselineskip}
\vspace{-0.575\baselineskip}

\begin{figure}[H]
\centering
\includegraphics[keepaspectratio,width=5.51181in,height=7.77019in]{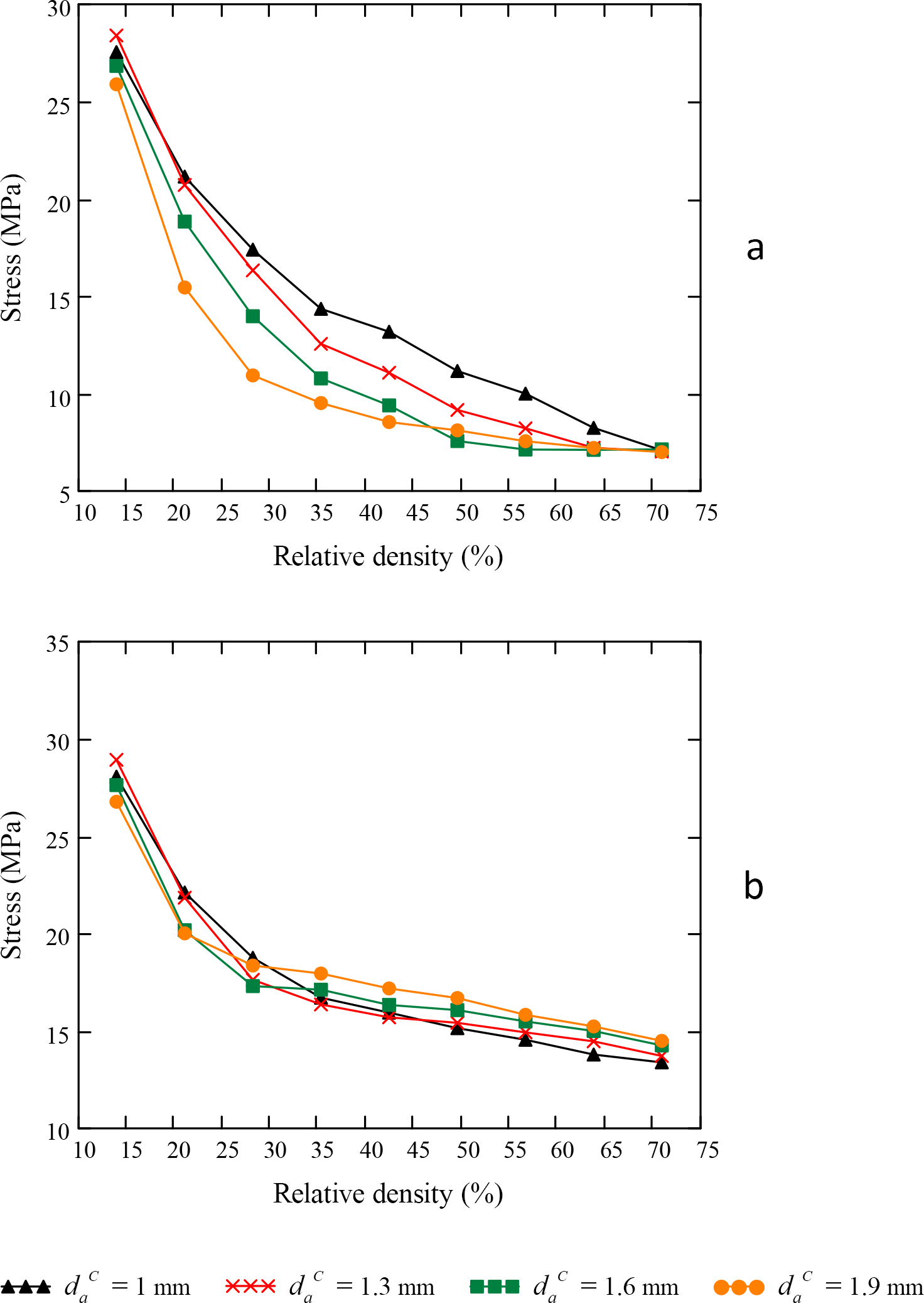}
\caption{Диаграммы зависимости \(\sigma_{\max}\) в сотовом заполнителе композитных пластин от \(\rho_{rel}\) при \(F_{y} = 30\)~Н в первой постановке численных экспериментов при жестком защемлении (a) \dualcite{ref77,ref78}{{[}\bibnum{ref77}, \bibnum{ref78}{]}}{{[}\bibnum{ref77}, \bibnum{ref78}{]}} и опирании с упругим поворотом~(b) \dualcite{ref78}{{[}\bibnum{ref78}{]}}{{[}\bibnum{ref78}{]}}}\label{fig:23}
\end{figure}

\clearpage
На рисунке \ref{fig:23} представлены диаграммы зависимости \(\sigma_{\max}\) в сотовом заполнителе композитных пластин от \(\rho_{rel}\) при \(F_{y} = 30\) Н в первой постановке численных экспериментов при жестком защемлении (рисунок \ref{fig:23}a) и опирании с упругим поворотом (рисунок \ref{fig:23}b). Диаграммы изменения напряжений в сотовых прослойках содержат результаты только трехмерного конечно-элементного моделирования в Comsol Multiphysics без сопоставления с результатами алгоритмов решения плоской задачи, поскольку в алгоритмах соты представляются как непрерывная среда с усредненными свойствами. В условиях жесткого защемления композитов при увеличении \(\rho_{rel}\) заполнителя от 14 до 71 \% разница между значениями \(\sigma_{\max}\) в тетракиральных сотах с различной дискретизацией сначала интенсивно увеличивается, а затем уменьшается до минимального значения (рисунок \ref{fig:23}a). Однако при опирании с упругим поворотом на всем диапазоне изменения \(\rho_{rel}\) максимальные напряжения в тетракиральных сотах с различной дискретизацией имеют одинаковый разброс в пределах погрешности метода конечных элементов (рисунок \ref{fig:23}b). При этом дискретизация сот не влияет на значения \(\sigma_{\max}\) во внешних слоях композитных пластин (рисунки \ref{fig:19}\iflatexml--\else--\fi\ref{fig:22}).

\begin{figure}[H]
\centering
\includegraphics[keepaspectratio,width=5.51181in,height=7.77019in]{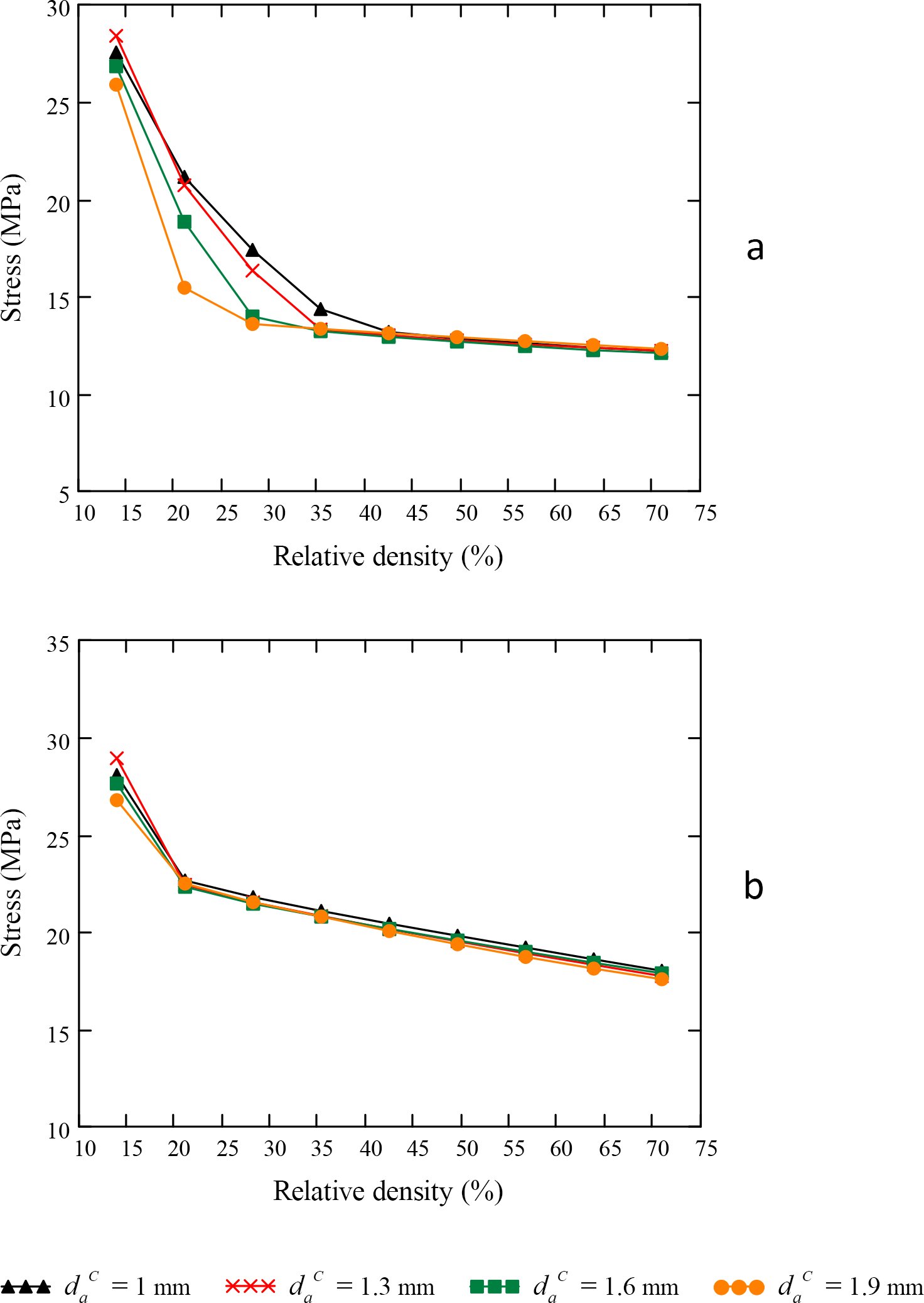}
\caption{Диаграммы зависимости \(\sigma_{\max}\) в трехслойном композите от \(\rho_{rel}\) при \(F_{y} = 30\)~Н в первой постановке численных экспериментов при жестком защемлении (a) \dualcite{ref77}{{[}\bibnum{ref77}{]}}{{[}\bibnum{ref77}{]}} и опирании с упругим поворотом (b)}\label{fig:24}
\end{figure}

На рисунке \ref{fig:24} представлены диаграммы зависимости \(\sigma_{\max}\) в трехслойном композите от \(\rho_{rel}\) при \(F_{y} = 30\) Н в первой постановке численных экспериментов при жестком защемлении (рисунок \ref{fig:24}a) и опирании с упругим поворотом (рисунок \ref{fig:24}b).

\vspace{-0.575\baselineskip}

\clearpage
\iflatexml\else\phantomsection\fi
\subsection{3.2. Результаты численного моделирования композитных пластин во второй постановке численных экспериментов при постоянном объеме твердого тела сот и варьируемой толщине заполнителя}

\begin{figure}[H]
\centering
\includegraphics[keepaspectratio,width=5.47244in,height=7.94303in]{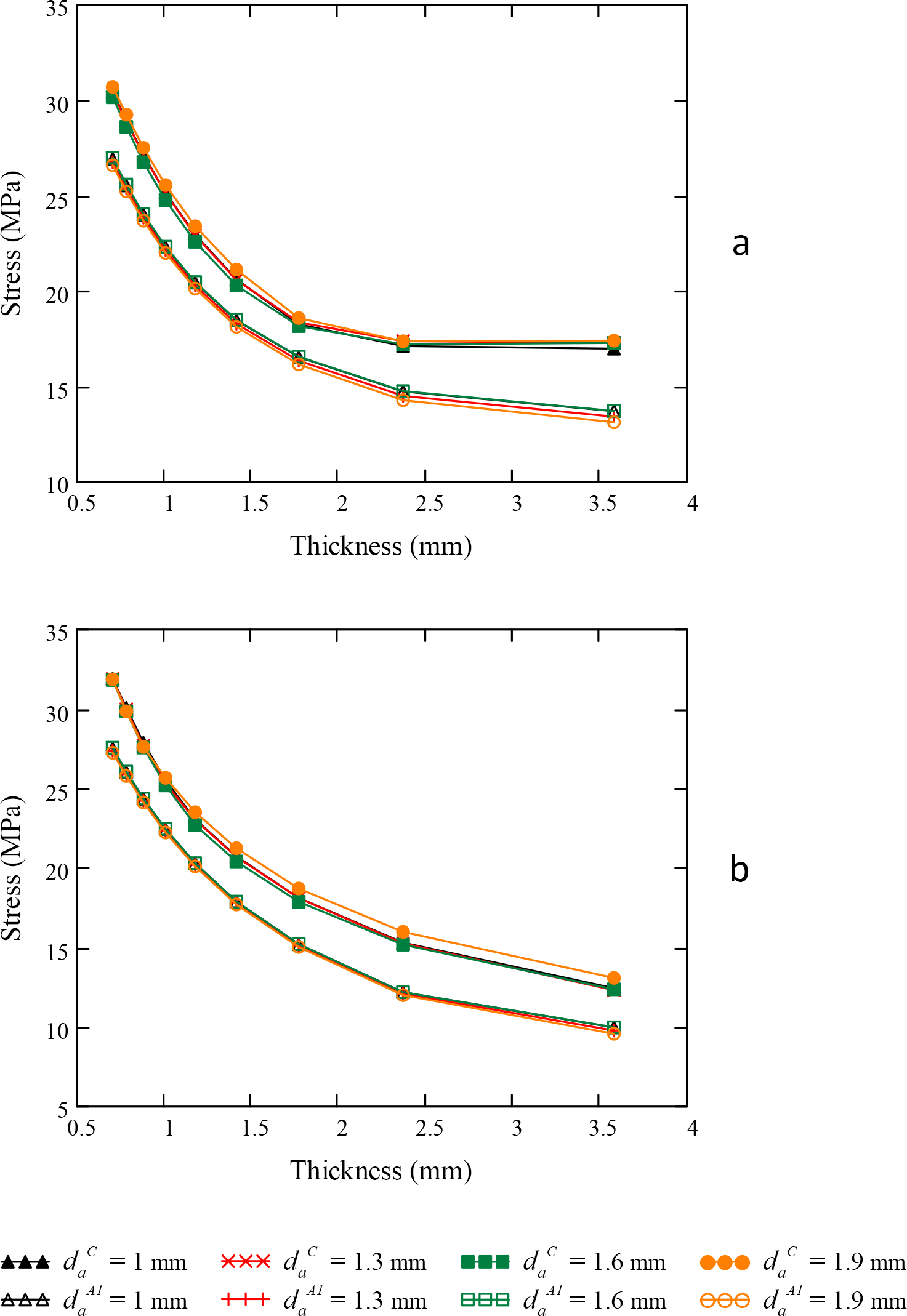}
\caption{Диаграммы зависимости \(\sigma_{\max}\) в верхнем (a) и нижнем (b) сплошном слое композитных пластин от \(t_{cl}\) при \(F_{y} = 60\) Н во второй постановке численных экспериментов при жестком защемлении, полученные с применением Comsol Multiphysics \dualcite{ref77,ref78}{{[}\bibnum{ref77}, \bibnum{ref78}{]}}{{[}\bibnum{ref77}, \bibnum{ref78}{]}} и Алгоритма 1}\label{fig:25}
\end{figure}

\begin{figure}[H]
\centering
\includegraphics[keepaspectratio,width=5.47244in,height=7.93931in]{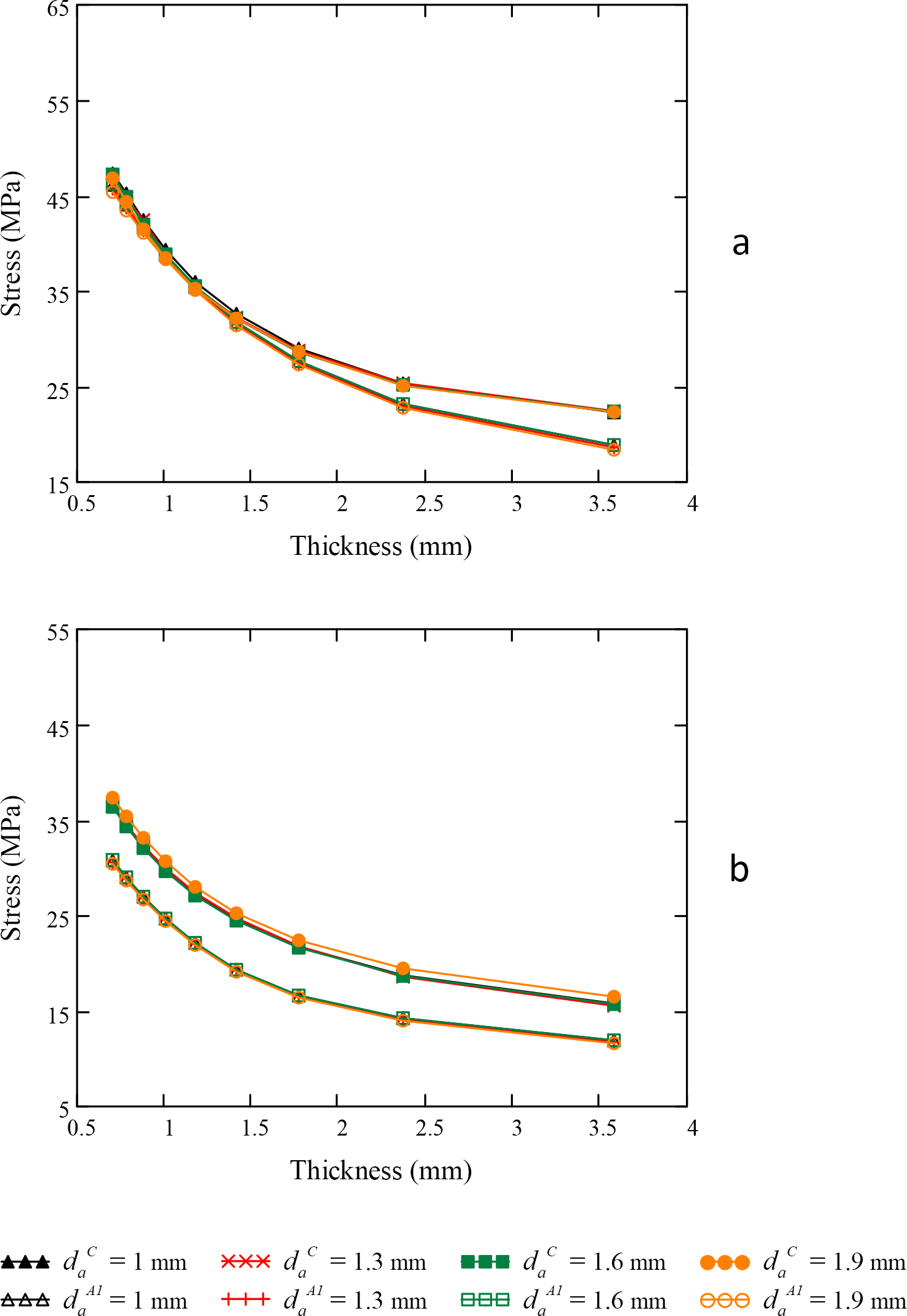}
\caption{Диаграммы зависимости \(\sigma_{\max}\) в верхнем (a) и нижнем (b) сплошном слое композитных пластин от \(t_{cl}\) при \(F_{y} = 60\) Н во второй постановке численных экспериментов при опирании с упругим поворотом, полученные с применением Comsol Multiphysics \dualcite{ref78}{{[}\bibnum{ref78}{]}}{{[}\bibnum{ref78}{]}} и Алгоритма 1}\label{fig:26}
\end{figure}

\begin{figure}[H]
\centering
\includegraphics[keepaspectratio,width=5.46643in,height=7.9343in]{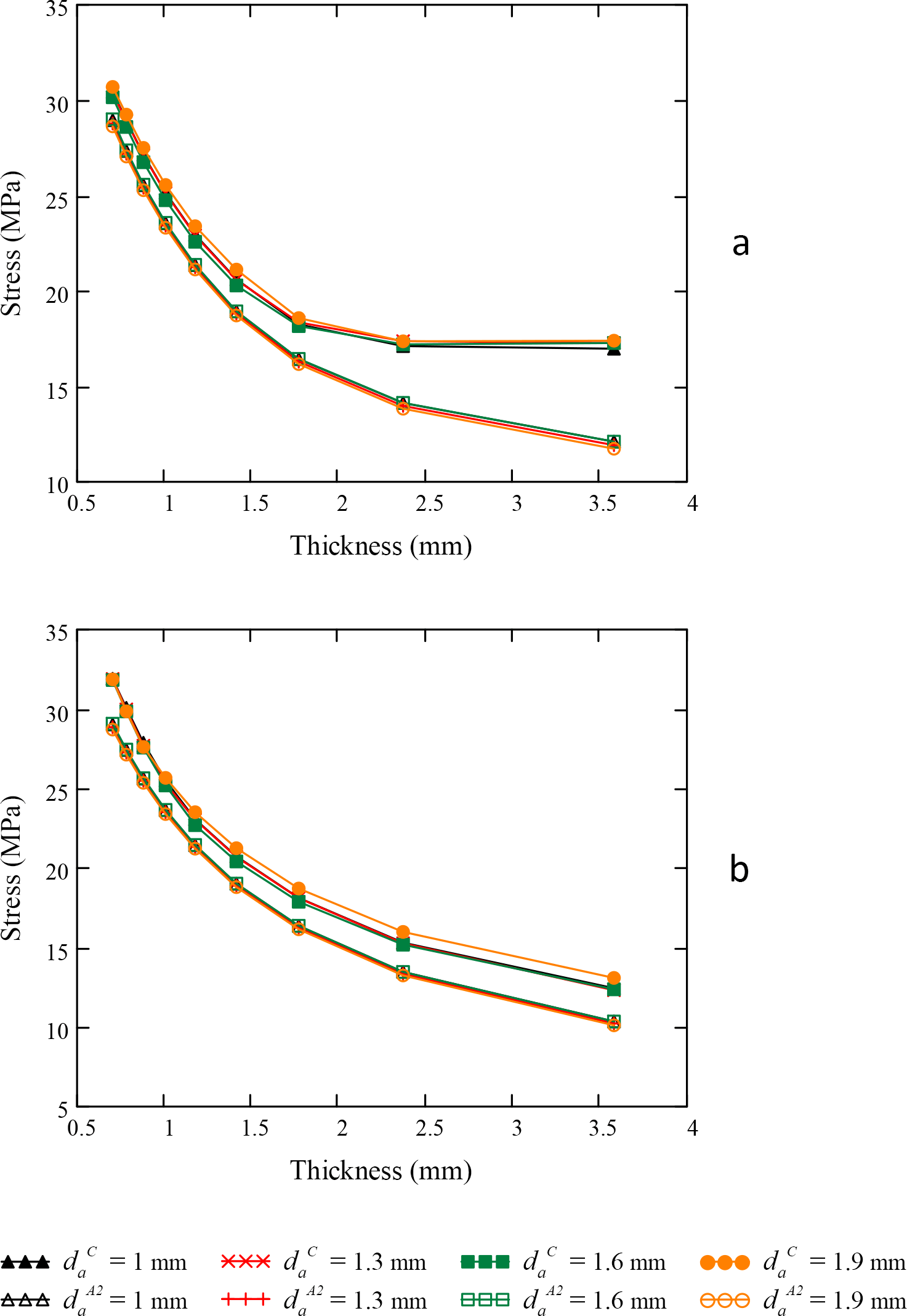}
\caption{Диаграммы зависимости \(\sigma_{\max}\) в верхнем (a) и нижнем (b) сплошном слое композитных пластин от \(t_{cl}\) при \(F_{y} = 60\) Н во второй постановке численных экспериментов при жестком защемлении, полученные с применением Comsol Multiphysics \dualcite{ref77,ref78}{{[}\bibnum{ref77}, \bibnum{ref78}{]}}{{[}\bibnum{ref77}, \bibnum{ref78}{]}} и Алгоритма 2}\label{fig:27}
\end{figure}

\begin{figure}[H]
\centering
\includegraphics[keepaspectratio,width=5.46899in,height=7.9343in]{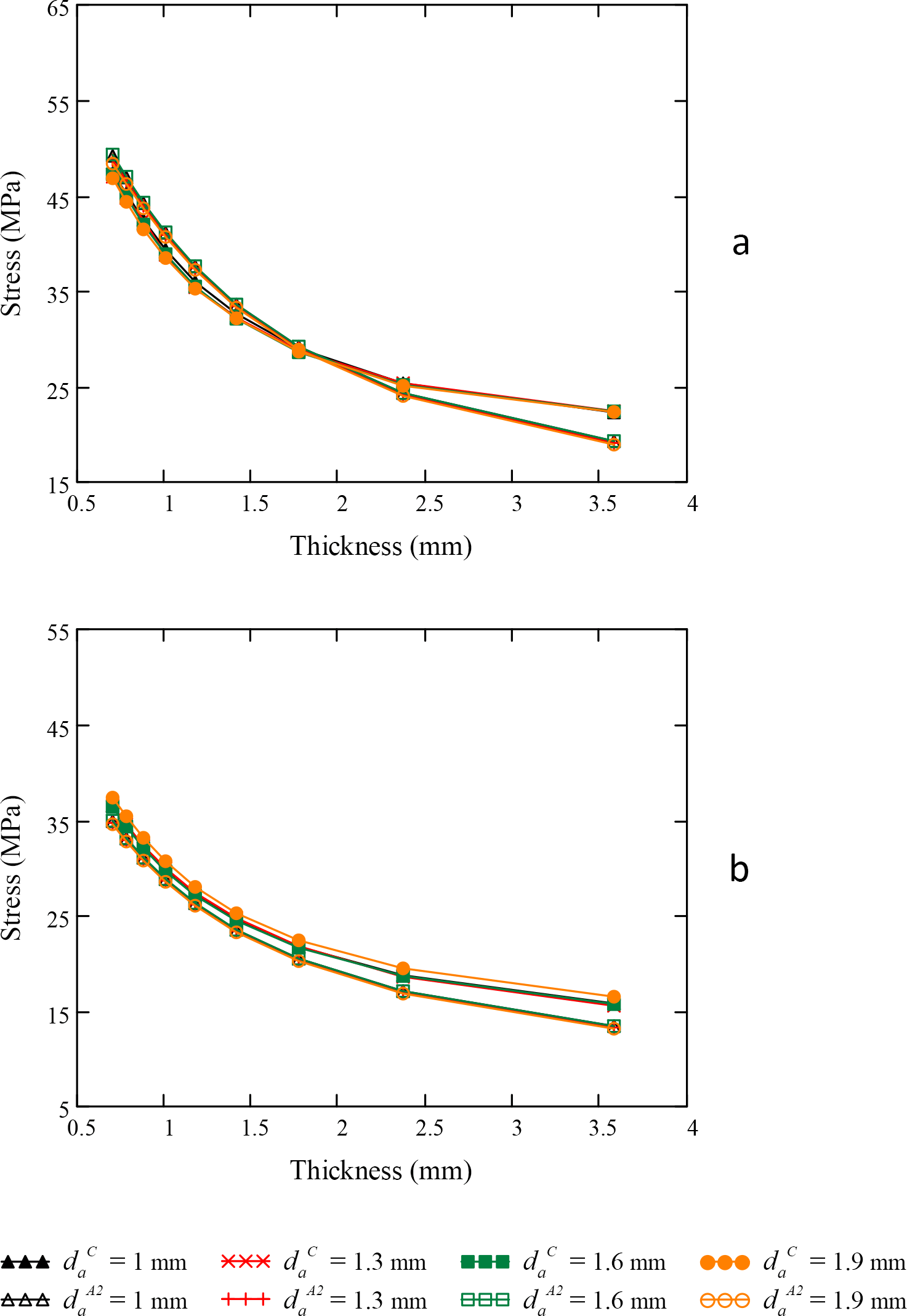}
\caption{Диаграммы зависимости \(\sigma_{\max}\) в верхнем (a) и нижнем (b) сплошном слое композитных пластин от \(t_{cl}\) при \(F_{y} = 60\) Н во второй постановке численных экспериментов при опирании с упругим поворотом, полученные с применением Comsol Multiphysics \dualcite{ref78}{{[}\bibnum{ref78}{]}}{{[}\bibnum{ref78}{]}} и Алгоритма 2}\label{fig:28}
\end{figure}

Диаграммы зависимости \(\sigma_{\max}\) во внешних слоях композитных пластин от \(t_{cl}\) при \(F_{y} = 60\) Н во второй постановке численных экспериментов при жестком защемлении и опирании с упругим поворотом представлены на рисунках \ref{fig:25} и \ref{fig:26} соответственно. При этом \(\rho_{rel}\) заполнителя изменяется в аналогичном диапазоне с первой постановкой экспериментов, однако соты имеют разную толщину из-за равного объема их твердого тела. На рисунках \ref{fig:27} и \ref{fig:28} представлены диаграммы с аналогичной зависимостью в тех же граничных условиях, однако с использованием дополнительного алгоритма решения плоской задачи, в котором сплошные слои моделируются с применением несовместных конечных элементов. Как и в случае первой постановки, графики, полученные с применением системы Comsol Multiphysics и Алгоритма 2, демонстрируют хорошее соответствие, а Алгоритм 1 приводит к более заметной погрешности в большинстве случаев.

\begin{figure}[H]
\centering
\includegraphics[keepaspectratio,width=5.47244in,height=7.71469in]{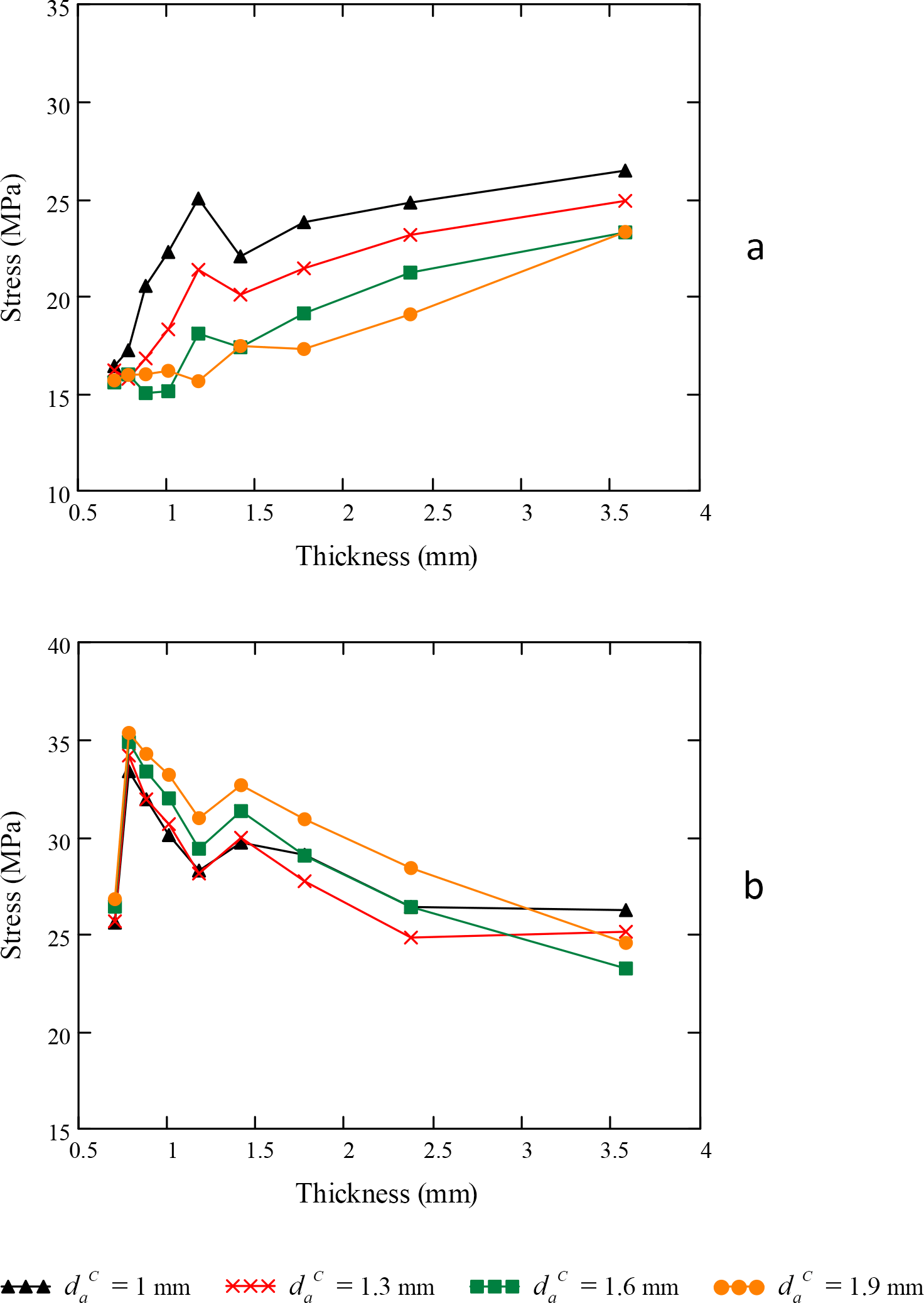}
\caption{Диаграммы зависимости \(\sigma_{\max}\) в сотовом заполнителе композитных пластин от \(t_{cl}\) при \(F_{y} = 60\)~Н во второй постановке численных экспериментов при жестком защемлении (a) \dualcite{ref77,ref78}{{[}\bibnum{ref77}, \bibnum{ref78}{]}}{{[}\bibnum{ref77}, \bibnum{ref78}{]}} и опирании с упругим поворотом~(b) \dualcite{ref78}{{[}\bibnum{ref78}{]}}{{[}\bibnum{ref78}{]}}}\label{fig:29}
\end{figure}

На рисунке \ref{fig:29} представлены диаграммы зависимости \(\sigma_{\max}\) в сотовом заполнителе композитных пластин от \(t_{cl}\) при \(F_{y} = 60\) Н во второй постановке численных экспериментов. При жестком защемлении пластин с увеличением \(t_{cl}\) (уменьшением \(\rho_{rel}\)) немонотонно возрастают \(\sigma_{\max}\) в тетракиральных сотах с различной дискретизацией (рисунок \ref{fig:29}a), однако по сравнению с рисунком \ref{fig:23}a интенсивность возрастания ниже. При опирании с упругим поворотом в тетракиральных сотах с различной дискретизацией наблюдается более узкий разброс значений \(\sigma_{\max}\) со слабо выраженной закономерностью (рисунок \ref{fig:29}b). При этом дискретизация сот также не влияет на значения \(\sigma_{\max}\) во внешних слоях композитных пластин (рисунки \ref{fig:25}\iflatexml--\else--\fi\ref{fig:28}).

\vspace{-0.575\baselineskip}
\vspace{-0.575\baselineskip}
\vspace{-0.575\baselineskip}
\vspace{-0.575\baselineskip}

\begin{figure}[H]
\centering
\includegraphics[keepaspectratio,width=6.69291in,height=3.12739in]{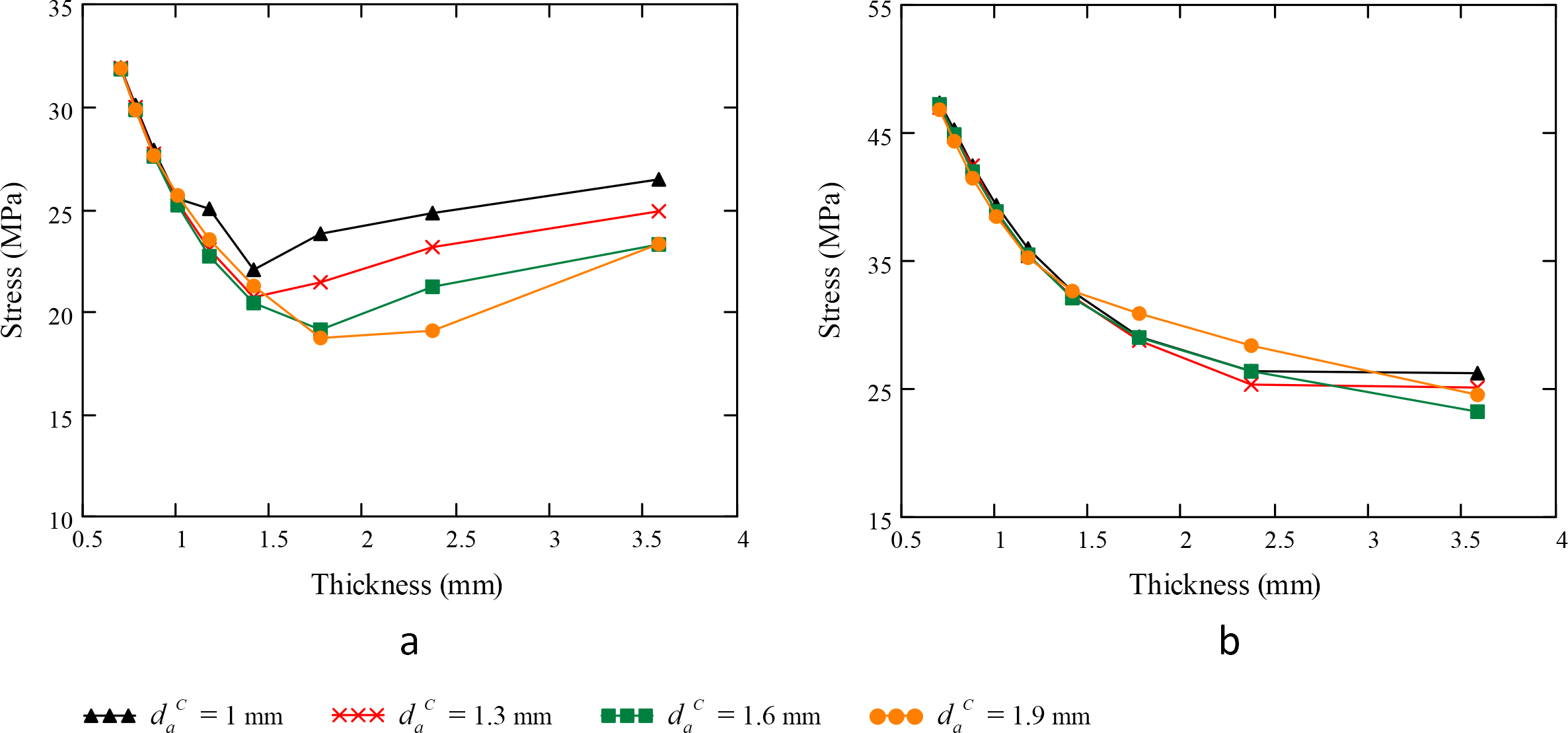}
\caption{Диаграммы зависимости \(\sigma_{\max}\) в трехслойном композите от \(t_{cl}\) при \(F_{y} = 60\)~Н во второй постановке численных экспериментов при жестком защемлении (a) \dualcite{ref77}{{[}\bibnum{ref77}{]}}{{[}\bibnum{ref77}{]}} и опирании с упругим поворотом (b)}\label{fig:30}
\end{figure}

На рисунке \ref{fig:30} представлены диаграммы зависимости \(\sigma_{\max}\) в трехслойном композите от \(t_{cl}\) при \(F_{y} = 60\) Н во второй постановке численных экспериментов при жестком защемлении (рисунок \ref{fig:30}a) и опирании с упругим поворотом (рисунок \ref{fig:30}b). В условиях жесткого защемления (рисунок \ref{fig:30}a) при варьировании толщины пластин посредством изменения относительной плотности заполнителя можно получить минимальные значения максимального напряжения в композитах. При этом минимальное значение максимального напряжения у большинства композитных пластин находится в точке перехода от сотового заполнителя к сплошному слою при определенной толщине сот.

\iflatexml\else\phantomsection\fi
\subsection{3.3. Напряженное состояние трехслойных пластин в первой и во второй постановках численных экспериментов при жестком защемлении и опирании с упругим поворотом}

\begin{figure}[H]
\centering
\includegraphics[keepaspectratio,width=3.93559in,height=8.20787in]{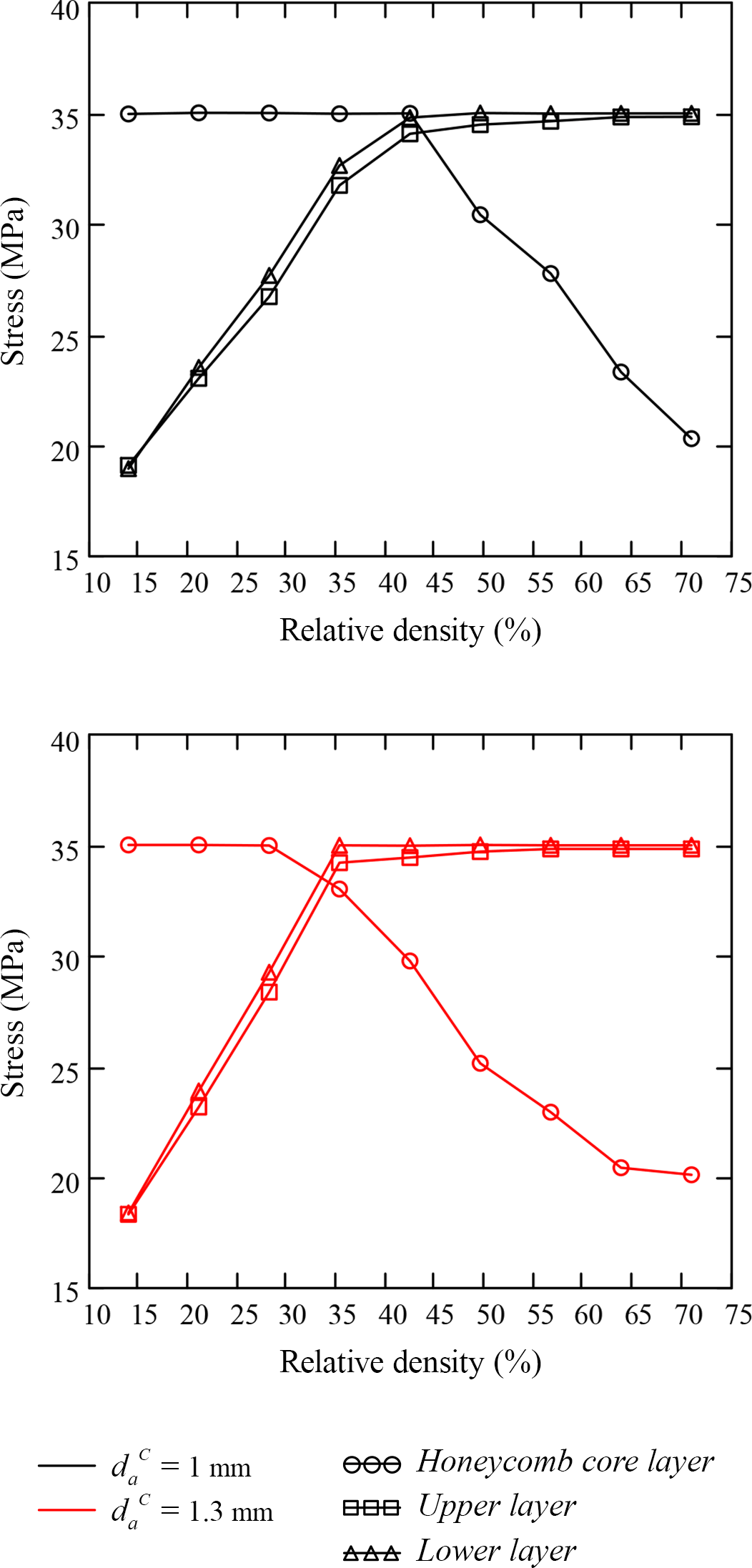}
\caption{Диаграммы зависимости \(\sigma_{\max}\) в слоях композитных пластин от \(\rho_{rel}\) при \(\sigma_{\max} = \sigma_{pl}\) в первой постановке численных экспериментов при \(d_{a} = 1\) и \(d_{a} = 1.3\) мм в условиях жесткого защемления \dualcite{ref77}{{[}\bibnum{ref77}{]}}{{[}\bibnum{ref77}{]}}}\label{fig:31}
\end{figure}

\begin{figure}[H]
\centering
\includegraphics[keepaspectratio,width=3.93559in,height=8.20787in]{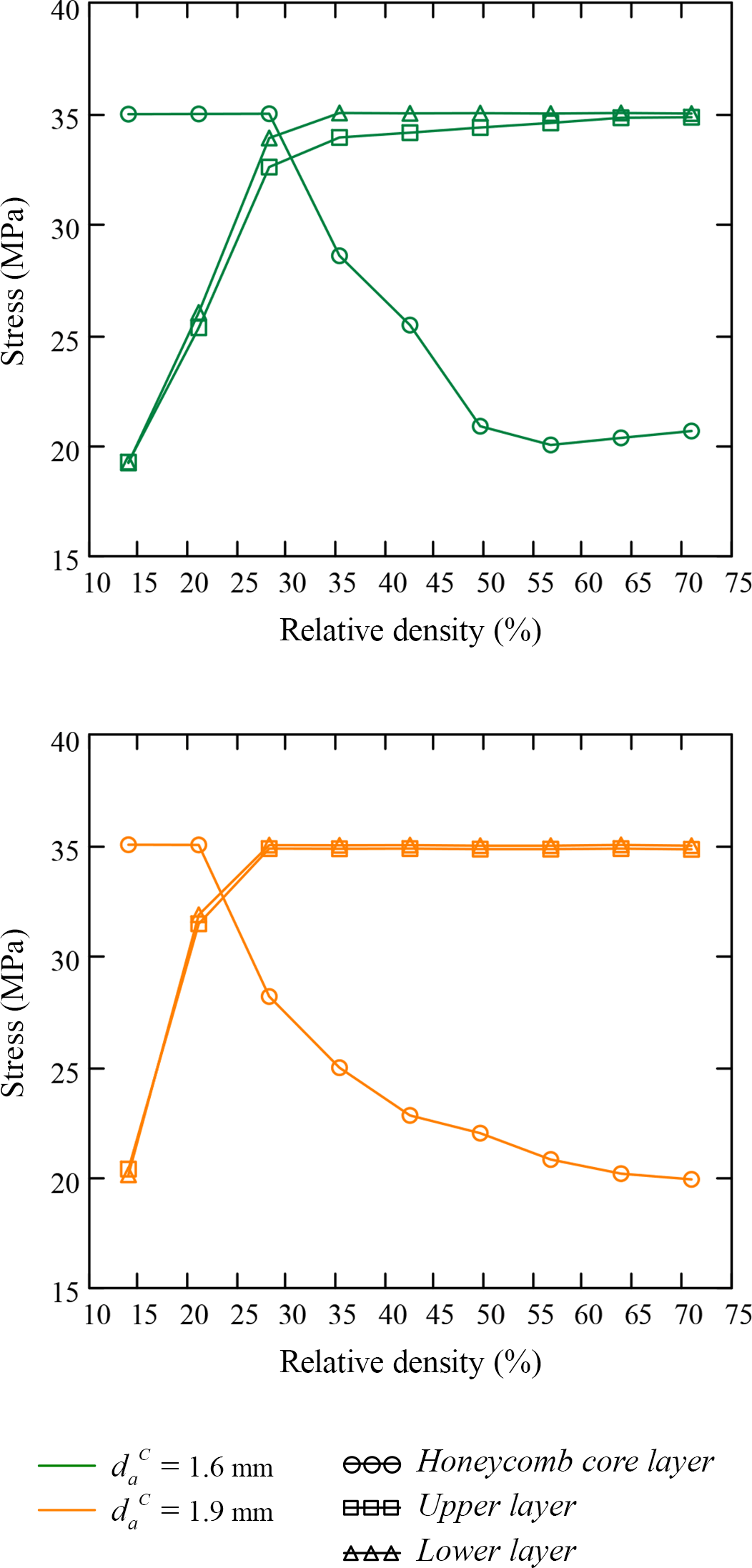}
\caption{Диаграммы зависимости \(\sigma_{\max}\) в слоях композитных пластин от \(\rho_{rel}\) при \(\sigma_{\max} = \sigma_{pl}\) в первой постановке численных экспериментов при \(d_{a} = 1.6\) и \(d_{a} = 1.9\) мм в условиях жесткого защемления \dualcite{ref77}{{[}\bibnum{ref77}{]}}{{[}\bibnum{ref77}{]}}}\label{fig:32}
\end{figure}

\begin{figure}[H]
\centering
\includegraphics[keepaspectratio,width=3.93559in,height=8.20787in]{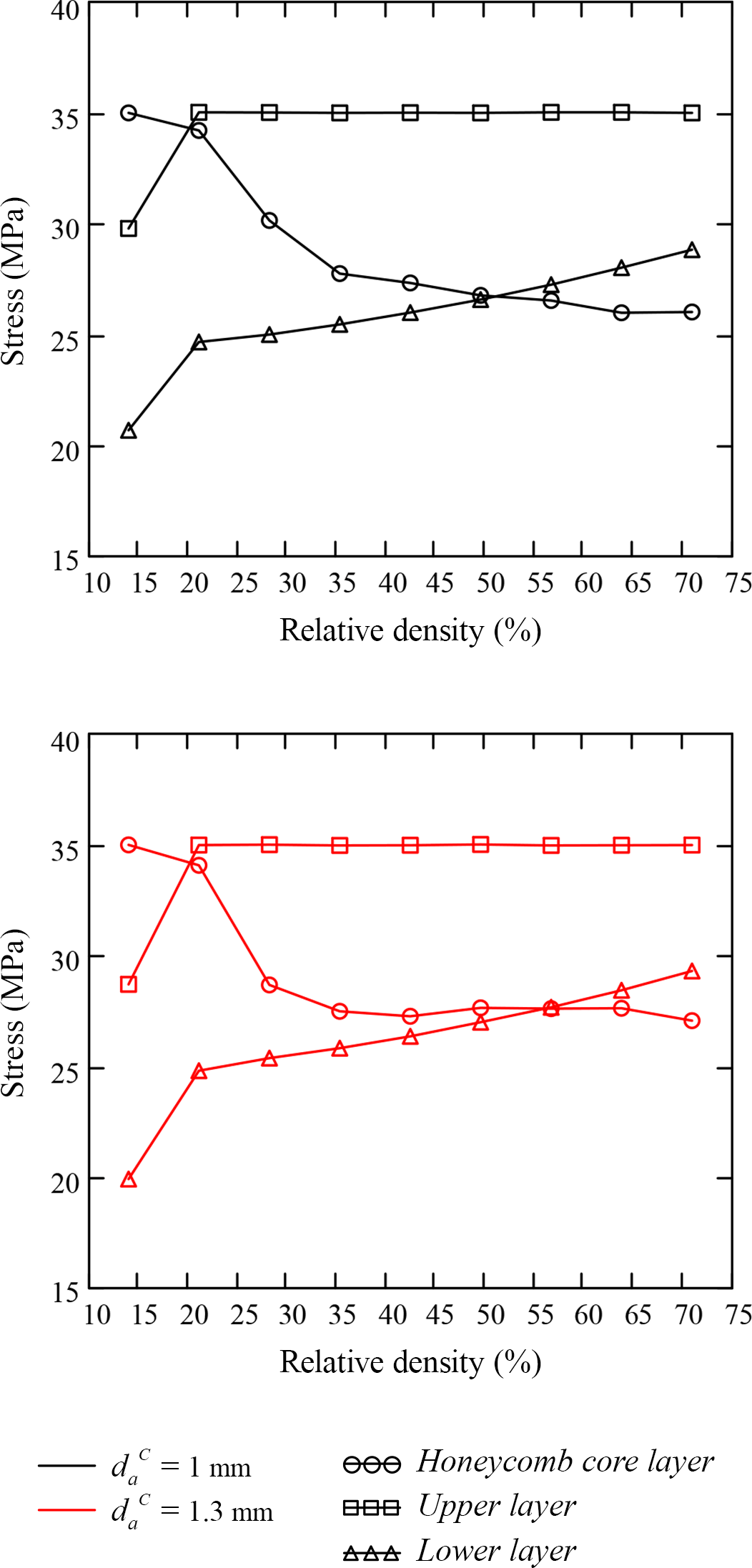}
\caption{Диаграммы зависимости \(\sigma_{\max}\) в слоях композитных пластин от \(\rho_{rel}\) при \(\sigma_{\max} = \sigma_{pl}\) в первой постановке численных экспериментов при \(d_{a} = 1\) и \(d_{a} = 1.3\) мм в условиях опирания с упругим поворотом}\label{fig:33}
\end{figure}

\begin{figure}[H]
\centering
\includegraphics[keepaspectratio,width=3.93559in,height=8.20787in]{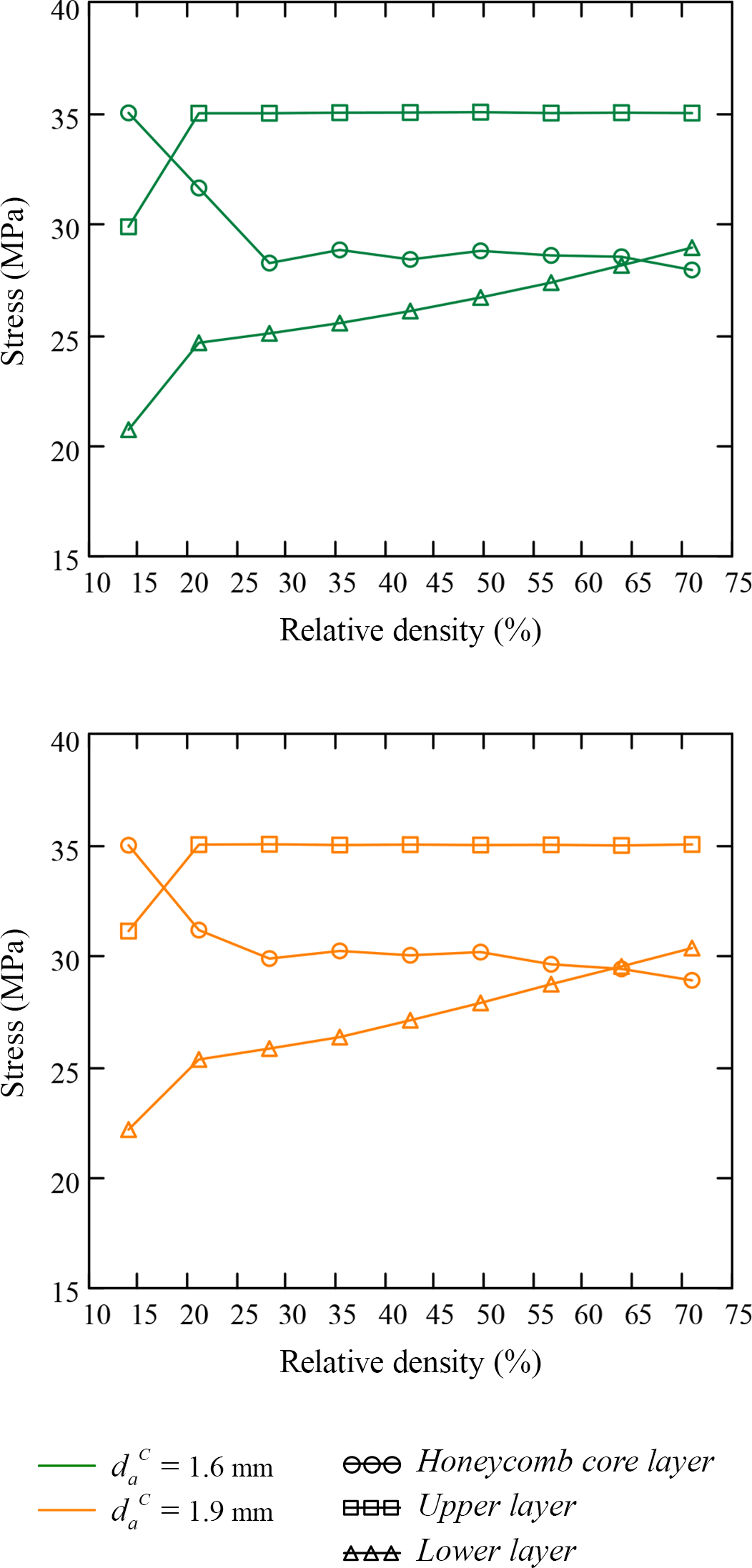}
\caption{Диаграммы зависимости \(\sigma_{\max}\) в слоях композитных пластин от \(\rho_{rel}\) при \(\sigma_{\max} = \sigma_{pl}\) в первой постановке численных экспериментов при \(d_{a} = 1.6\) и \(d_{a} = 1.9\) мм в условиях опирания с упругим поворотом}\label{fig:34}
\end{figure}

\begin{figure}[H]
\centering
\includegraphics[keepaspectratio,width=3.93559in,height=8.20787in]{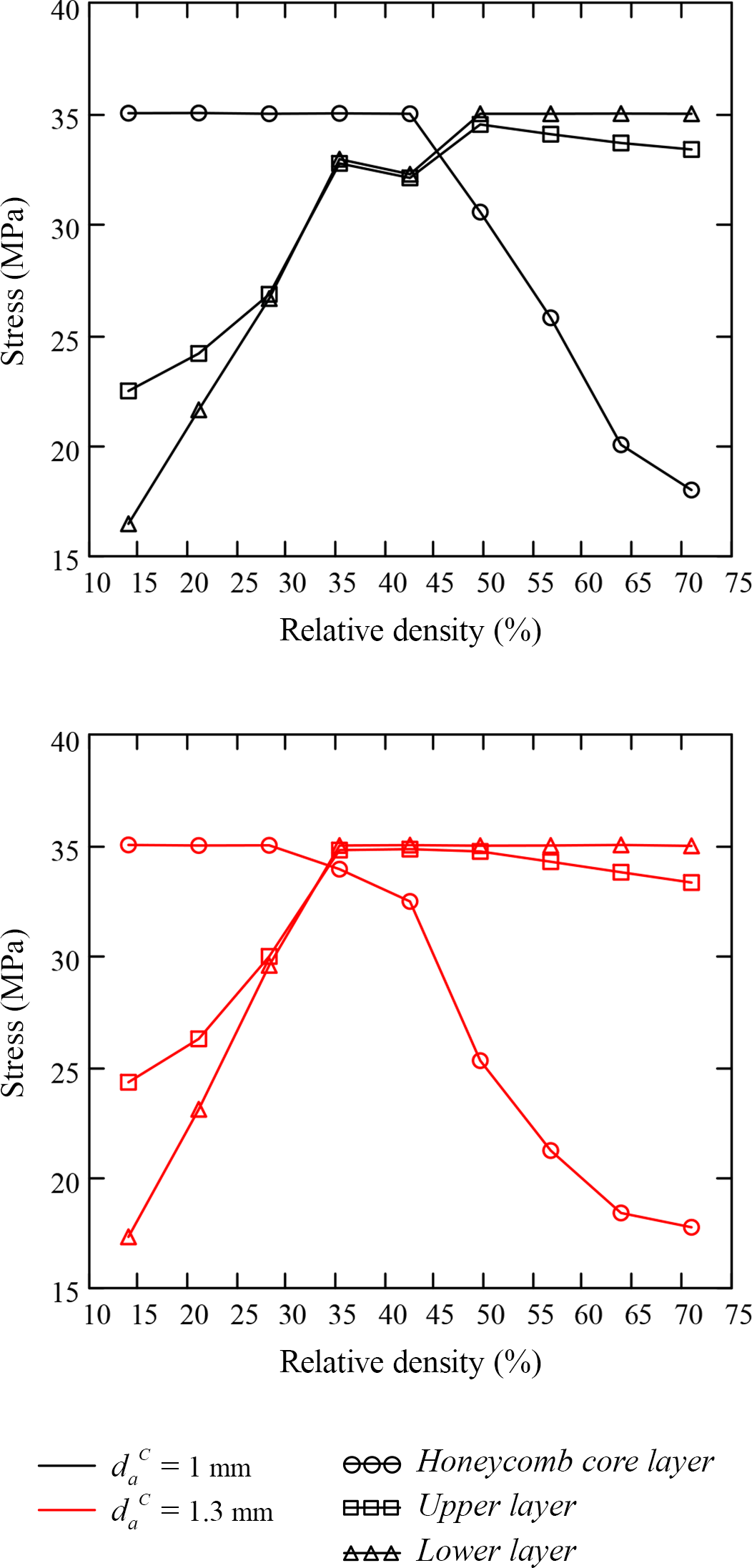}
\caption{Диаграммы зависимости \(\sigma_{\max}\) в слоях композитных пластин от \(\rho_{rel}\) при \(\sigma_{\max} = \sigma_{pl}\) во второй постановке численных экспериментов при \(d_{a} = 1\) и \(d_{a} = 1.3\) мм в условиях жесткого защемления \dualcite{ref77}{{[}\bibnum{ref77}{]}}{{[}\bibnum{ref77}{]}}}\label{fig:35}
\end{figure}

\begin{figure}[H]
\centering
\includegraphics[keepaspectratio,width=3.93559in,height=8.20787in]{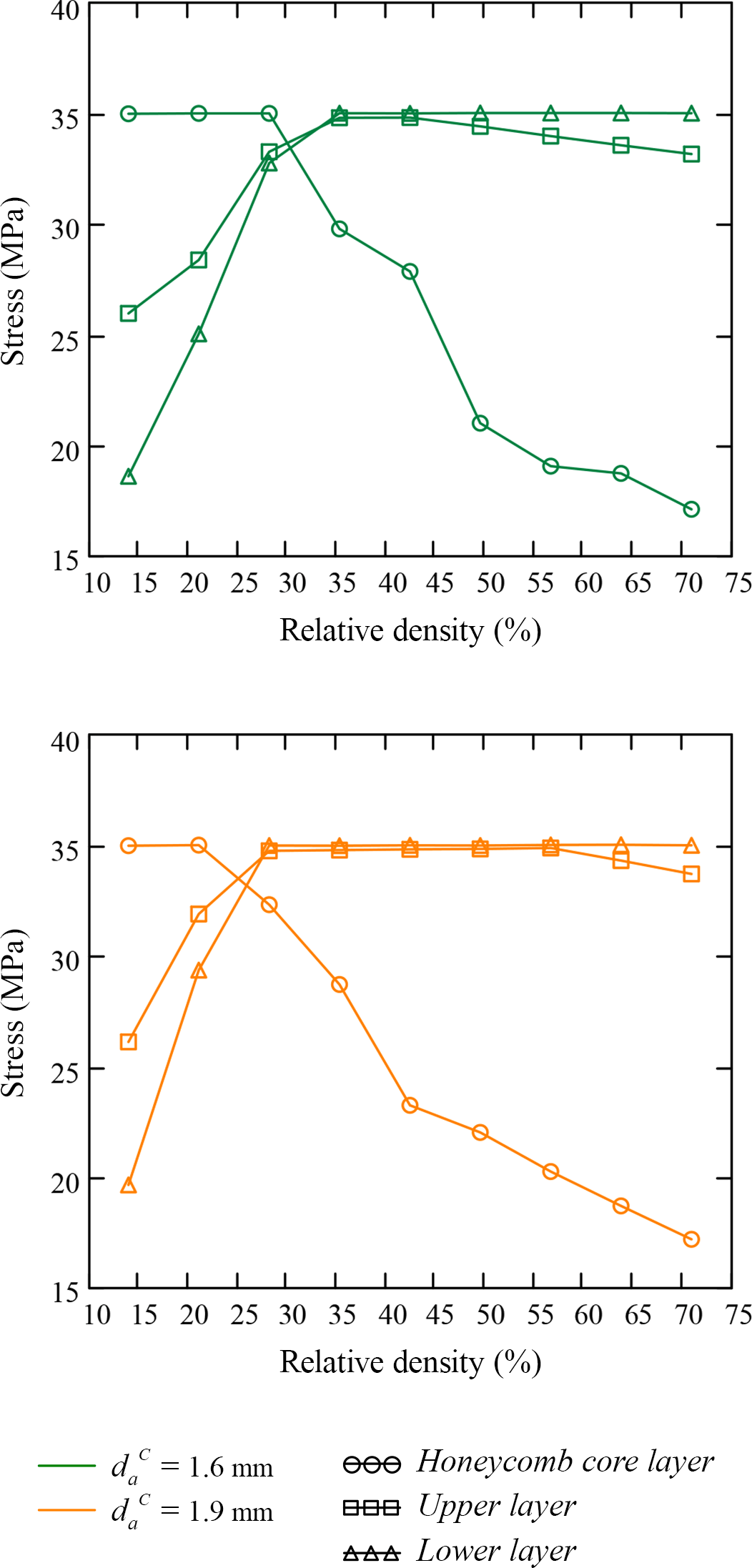}
\caption{Диаграммы зависимости \(\sigma_{\max}\) в слоях композитных пластин от \(\rho_{rel}\) при \(\sigma_{\max} = \sigma_{pl}\) во второй постановке численных экспериментов при \(d_{a} = 1.6\) и \(d_{a} = 1.9\) мм в условиях жесткого защемления \dualcite{ref77}{{[}\bibnum{ref77}{]}}{{[}\bibnum{ref77}{]}}}\label{fig:36}
\end{figure}

\begin{figure}[H]
\centering
\includegraphics[keepaspectratio,width=3.93559in,height=8.20787in]{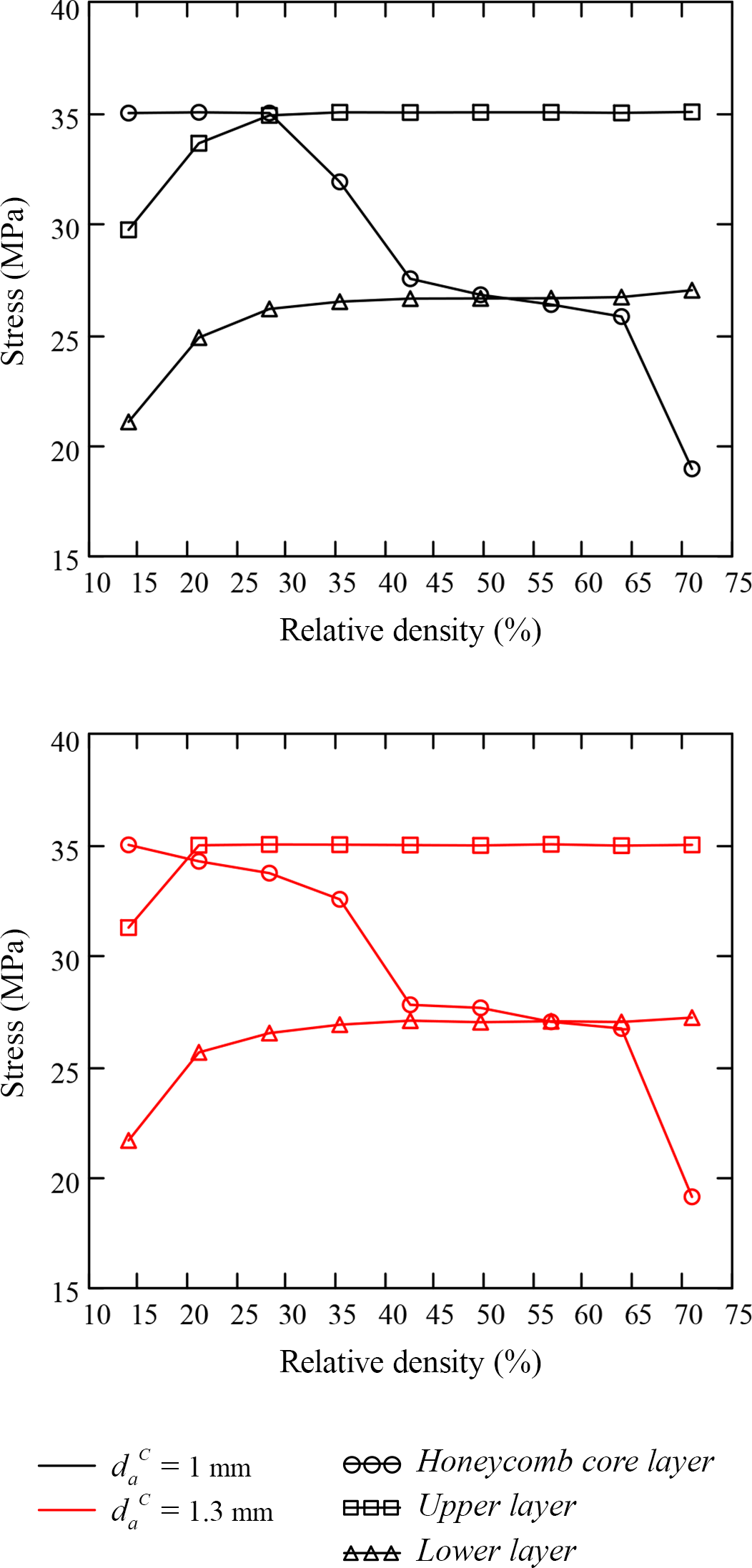}
\caption{Диаграммы зависимости \(\sigma_{\max}\) в слоях композитных пластин от \(\rho_{rel}\) при \(\sigma_{\max} = \sigma_{pl}\) во второй постановке численных экспериментов при \(d_{a} = 1\) и \(d_{a} = 1.3\) мм в условиях опирания с упругим поворотом}\label{fig:37}
\end{figure}

\begin{figure}[H]
\centering
\includegraphics[keepaspectratio,width=3.93559in,height=8.20787in]{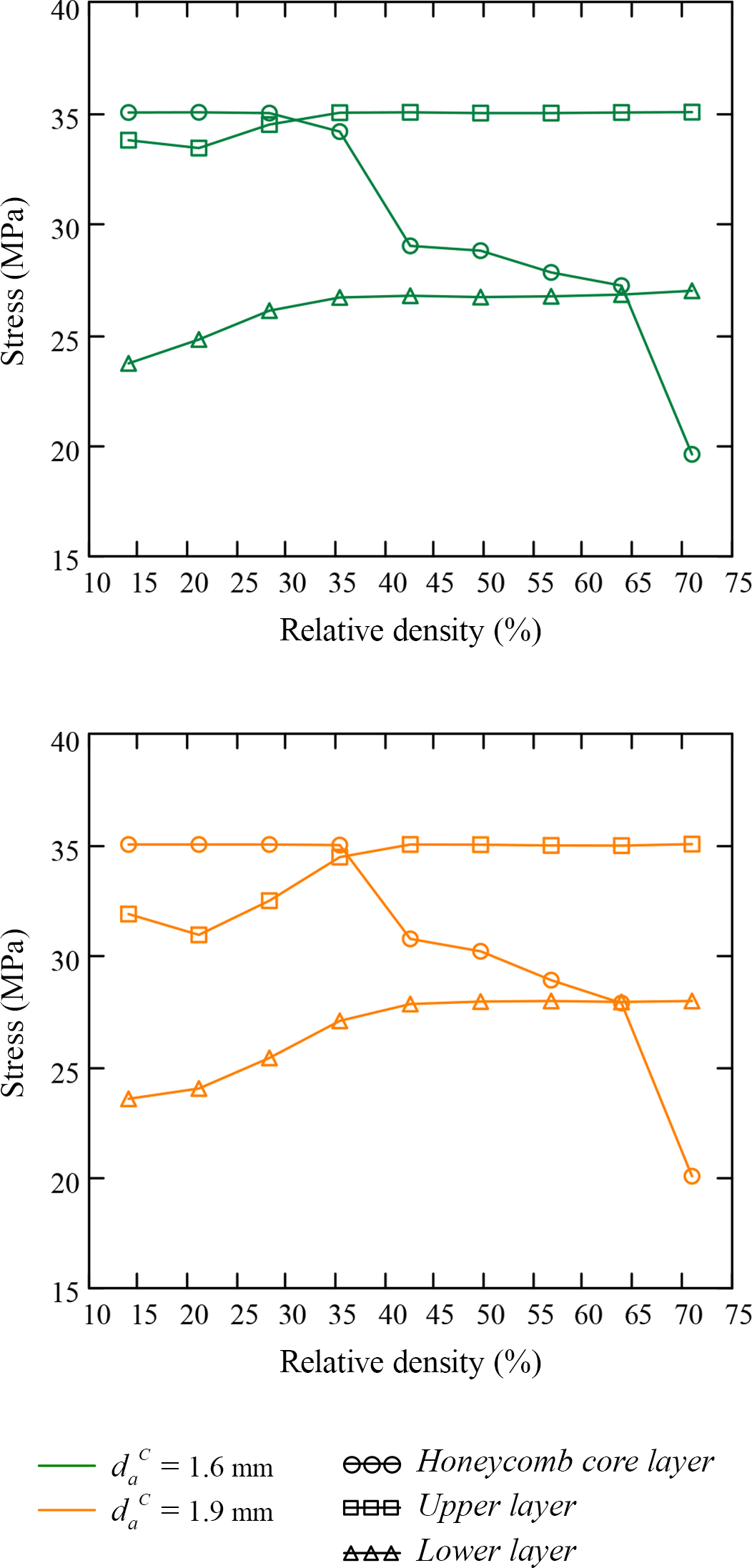}
\caption{Диаграммы зависимости \(\sigma_{\max}\) в слоях композитных пластин от \(\rho_{rel}\) при \(\sigma_{\max} = \sigma_{pl}\) во второй постановке численных экспериментов при \(d_{a} = 1.6\) и \(d_{a} = 1.9\) мм в условиях опирания с упругим поворотом}\label{fig:38}
\end{figure}

Результаты анализа зависимости пиковых напряжений в слоях композитных пластин при \(\sigma_{\max} = \sigma_{pl}\) от \(\rho_{rel}\) заполнителя с использованием трехмерного конечно-элементного моделирования в первой постановке численных экспериментов представлены при жестком защемлении (таблицы \ref{tab:1}, \ref{tab:2}, рисунки \ref{fig:31}, \ref{fig:32}) и опирании с упругим поворотом (таблицы \ref{tab:3}, \ref{tab:4}, рисунки \ref{fig:33}, \ref{fig:34}). Также результаты аналогичной зависимости представлены во второй постановке численных экспериментов при жестком защемлении (таблицы \ref{tab:5}, \ref{tab:6}, рисунки \ref{fig:35}, \ref{fig:36}) и опирании с упругим поворотом (таблицы \ref{tab:7}, \ref{tab:8}, рисунки \ref{fig:37}, \ref{fig:38}).

С увеличением \(\rho_{rel}\) заполнителя у композитных пластин при жестком защемлении переход критических напряжений от сотовых прослоек к сплошным слоям в обеих постановках осуществляется с одинаковой очередностью в равных диапазонах (рисунки \ref{fig:31}, \ref{fig:32}, \ref{fig:35}, \ref{fig:36}). При этом из диаграмм на рисунках \ref{fig:23}a и \ref{fig:29}a следует, что при постоянном значении \(\rho_{rel}\) и снижении дискретизации сотовых структур снижается \(\sigma_{\max}\) в заполнителе для большинства значений \(\rho_{rel}\) из рассматриваемого диапазона. Аналогичный вывод можно сделать на основе диаграмм распределения напряжений для большинства композитов, рассмотренных в работе \dualcite{ref79}{{[}\bibnum{ref79}{]}}{{[}\bibnum{ref79}{]}}.

\begin{figure}[H]
\centering
\includegraphics[keepaspectratio,width=6.69375in,height=4.40486in]{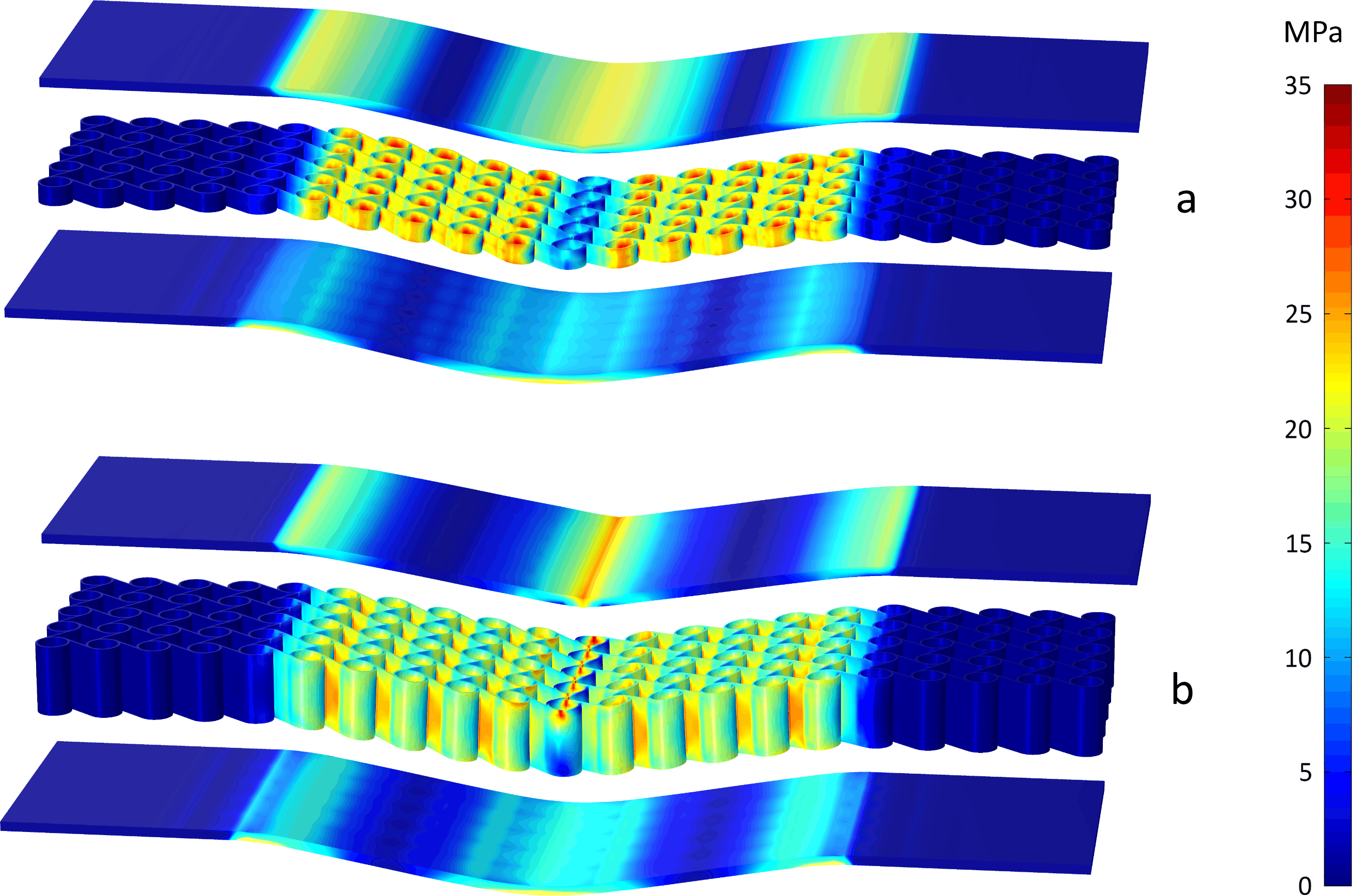}
\caption{Диаграммы распределения напряжений в тонкой (a) и относительно толстой (b) композитной пластине при \(\sigma_{\max} = \sigma_{pl}\) в условиях жесткого защемления \dualcite{ref77}{{[}\bibnum{ref77}{]}}{{[}\bibnum{ref77}{]}}}\label{fig:39}
\end{figure}

На рисунке \ref{fig:39} представлены диаграммы распределения напряжений в композитных пластинах при \(d_{a} = 1.6\) мм, \(\rho_{rel} = 14\) \%, \(\sigma_{\max} = \sigma_{pl}\) в двух постановках численных экспериментов при жестком защемлении. Критические напряжения в относительно толстых композитных пластинах локализуются преимущественно в области приложения силы, а в тонких композитных пластинах имеют более равномерное распределение.

\clearpage
\begingroup
\setstretch{1.2}
\iflatexml
\begin{table}
\centering
\caption{Результаты анализа зависимости \(\sigma_{\max}\) в слоях композитных пластин при \(\sigma_{\max} = \sigma_{pl}\) от \(\rho_{rel}\) заполнителя в первой постановке численных экспериментов при \(d_{a} = 1\) и \(d_{a} = 1.3\) мм в условиях жесткого защемления}\label{tab:1}
\begin{tabular}{ccccccc}
\lxBeginTableHead
\lxTableColumnHead \(d_{a}\) (мм) & \lxTableColumnHead \(t_{sw}\) (мм) & \lxTableColumnHead \(F_{y}\) (Н) & \multicolumn{3}{c}{\lxTableColumnHead \(\sigma_{\max}\) (МПа)} & \lxTableColumnHead \(\rho_{rel}\) (\%) \\
& & & \lxTableColumnHead в сотах & \lxTableColumnHead в верхнем сплошном слое & \lxTableColumnHead в нижнем сплошном слое & \\
\lxEndTableHead
\toprule
\tabledataA{*}{\midrule}
\bottomrule
\end{tabular}
\end{table}
\else
\begin{longtable}[]{@{}
  >{\centering\arraybackslash}m{(\columnwidth - 12\tabcolsep) * \real{0.0923}}
  >{\centering\arraybackslash}m{(\columnwidth - 12\tabcolsep) * \real{0.1383}}
  >{\centering\arraybackslash}m{(\columnwidth - 12\tabcolsep) * \real{0.1383}}
  >{\centering\arraybackslash}m{(\columnwidth - 12\tabcolsep) * \real{0.1383}}
  >{\centering\arraybackslash}m{(\columnwidth - 12\tabcolsep) * \real{0.1772}}
  >{\centering\arraybackslash}m{(\columnwidth - 12\tabcolsep) * \real{0.1772}}
  >{\centering\arraybackslash}m{(\columnwidth - 12\tabcolsep) * \real{0.1383}}@{}}
\caption{Результаты анализа зависимости \(\sigma_{\max}\) в слоях композитных пластин при \(\sigma_{\max} = \sigma_{pl}\) от \(\rho_{rel}\) заполнителя в первой постановке численных экспериментов при \(d_{a} = 1\) и \(d_{a} = 1.3\) мм в условиях жесткого защемления}\label{tab:1}\\
\toprule\noalign{}
\multirow{2}{=}{\begin{minipage}[c]{\linewidth}\centering
\vspace{0.35\baselineskip}
\(d_{a}\)\\[-0.05\baselineskip](мм)
\end{minipage}} & \multirow{2}{=}{\begin{minipage}[c]{\linewidth}\centering
\vspace{0.35\baselineskip}
\(t_{sw}\)\\[-0.05\baselineskip](мм)
\end{minipage}} & \multirow{2}{=}{\begin{minipage}[c]{\linewidth}\centering
\vspace{0.35\baselineskip}
\(F_{y}\)\\[-0.05\baselineskip](Н)
\end{minipage}} & \multicolumn{3}{>{\centering\arraybackslash}m{(\columnwidth - 12\tabcolsep) * \real{0.4927} + 4\tabcolsep}}{
\begin{minipage}[c]{\linewidth}\centering
\(\sigma_{\max}\) (МПа)
\end{minipage}} & \multirow{2}{=}{\begin{minipage}[c]{\linewidth}\centering
\vspace{0.35\baselineskip}
\(\rho_{rel}\)\\[-0.05\baselineskip](\%)
\end{minipage}} \\
& & & \begin{minipage}[c]{\linewidth}\centering
\vspace{1.04\baselineskip}в сотах
\end{minipage} & \begin{minipage}[c]{\linewidth}\centering\setstretch{0.90}
\vspace{0.76\baselineskip}в верхнем сплошном слое
\end{minipage} & \begin{minipage}[c]{\linewidth}\centering\setstretch{0.90}
\vspace{0.76\baselineskip}в нижнем сплошном слое
\end{minipage} \\
\midrule\noalign{}
\endhead
\bottomrule\noalign{}
\endlastfoot
\tabledataA{=}{\midrule\noalign{}}
\end{longtable}
\fi
\endgroup
\clearpage

\begingroup
\setstretch{1.2}
\iflatexml
\begin{table}
\centering
\caption{Результаты анализа зависимости \(\sigma_{\max}\) в слоях композитных пластин при \(\sigma_{\max} = \sigma_{pl}\) от \(\rho_{rel}\) заполнителя в первой постановке численных экспериментов при \(d_{a} = 1.6\) и \(d_{a} = 1.9\) мм в условиях жесткого защемления}\label{tab:2}
\begin{tabular}{ccccccc}
\lxBeginTableHead
\lxTableColumnHead \(d_{a}\) (мм) & \lxTableColumnHead \(t_{sw}\) (мм) & \lxTableColumnHead \(F_{y}\) (Н) & \multicolumn{3}{c}{\lxTableColumnHead \(\sigma_{\max}\) (МПа)} & \lxTableColumnHead \(\rho_{rel}\) (\%) \\
& & & \lxTableColumnHead в сотах & \lxTableColumnHead в верхнем сплошном слое & \lxTableColumnHead в нижнем сплошном слое & \\
\lxEndTableHead
\toprule
\tabledataB{*}{\midrule}
\bottomrule
\end{tabular}
\end{table}
\else
\begin{longtable}[]{@{}
  >{\centering\arraybackslash}m{(\columnwidth - 12\tabcolsep) * \real{0.0923}}
  >{\centering\arraybackslash}m{(\columnwidth - 12\tabcolsep) * \real{0.1383}}
  >{\centering\arraybackslash}m{(\columnwidth - 12\tabcolsep) * \real{0.1383}}
  >{\centering\arraybackslash}m{(\columnwidth - 12\tabcolsep) * \real{0.1383}}
  >{\centering\arraybackslash}m{(\columnwidth - 12\tabcolsep) * \real{0.1773}}
  >{\centering\arraybackslash}m{(\columnwidth - 12\tabcolsep) * \real{0.1773}}
  >{\centering\arraybackslash}m{(\columnwidth - 12\tabcolsep) * \real{0.1383}}@{}}
\caption{Результаты анализа зависимости \(\sigma_{\max}\) в слоях композитных пластин при \(\sigma_{\max} = \sigma_{pl}\) от \(\rho_{rel}\) заполнителя в первой постановке численных экспериментов при \(d_{a} = 1.6\) и \(d_{a} = 1.9\) мм в условиях жесткого защемления}\label{tab:2}\\
\toprule\noalign{}
\multirow{2}{=}{\begin{minipage}[c]{\linewidth}\centering
\vspace{0.35\baselineskip}
\(d_{a}\)\\[-0.05\baselineskip](мм)
\end{minipage}} & \multirow{2}{=}{\begin{minipage}[c]{\linewidth}\centering
\vspace{0.35\baselineskip}
\(t_{sw}\)\\[-0.05\baselineskip](мм)
\end{minipage}} & \multirow{2}{=}{\begin{minipage}[c]{\linewidth}\centering
\vspace{0.35\baselineskip}
\(F_{y}\)\\[-0.05\baselineskip](Н)
\end{minipage}} & \multicolumn{3}{>{\centering\arraybackslash}m{(\columnwidth - 12\tabcolsep) * \real{0.4929} + 4\tabcolsep}}{
\begin{minipage}[c]{\linewidth}\centering
\(\sigma_{\max}\) (МПа)
\end{minipage}} & \multirow{2}{=}{\begin{minipage}[c]{\linewidth}\centering
\vspace{0.35\baselineskip}
\(\rho_{rel}\)\\[-0.05\baselineskip](\%)
\end{minipage}} \\
& & & \begin{minipage}[c]{\linewidth}\centering
\vspace{1.04\baselineskip}в сотах
\end{minipage} & \begin{minipage}[c]{\linewidth}\centering\setstretch{0.90}
\vspace{0.76\baselineskip}в верхнем сплошном слое
\end{minipage} & \begin{minipage}[c]{\linewidth}\centering\setstretch{0.90}
\vspace{0.76\baselineskip}в нижнем сплошном слое
\end{minipage} \\
\midrule\noalign{}
\endhead
\bottomrule\noalign{}
\endlastfoot
\tabledataB{=}{\midrule\noalign{}}
\end{longtable}
\fi
\endgroup
\clearpage

\begingroup
\setstretch{1.2}
\iflatexml
\begin{table}
\centering
\caption{Результаты анализа зависимости \(\sigma_{\max}\) в слоях композитных пластин при \(\sigma_{\max} = \sigma_{pl}\) от \(\rho_{rel}\) заполнителя в первой постановке численных экспериментов при \(d_{a} = 1\) и \(d_{a} = 1.3\) мм в условиях опирания с упругим поворотом}\label{tab:3}
\begin{tabular}{ccccccc}
\lxBeginTableHead
\lxTableColumnHead \(d_{a}\) (мм) & \lxTableColumnHead \(t_{sw}\) (мм) & \lxTableColumnHead \(F_{y}\) (Н) & \multicolumn{3}{c}{\lxTableColumnHead \(\sigma_{\max}\) (МПа)} & \lxTableColumnHead \(\rho_{rel}\) (\%) \\
& & & \lxTableColumnHead в сотах & \lxTableColumnHead в верхнем сплошном слое & \lxTableColumnHead в нижнем сплошном слое & \\
\lxEndTableHead
\toprule
\tabledataC{*}{\midrule}
\bottomrule
\end{tabular}
\end{table}
\else
\begin{longtable}[]{@{}
  >{\centering\arraybackslash}m{(\columnwidth - 12\tabcolsep) * \real{0.0923}}
  >{\centering\arraybackslash}m{(\columnwidth - 12\tabcolsep) * \real{0.1383}}
  >{\centering\arraybackslash}m{(\columnwidth - 12\tabcolsep) * \real{0.1383}}
  >{\centering\arraybackslash}m{(\columnwidth - 12\tabcolsep) * \real{0.1383}}
  >{\centering\arraybackslash}m{(\columnwidth - 12\tabcolsep) * \real{0.1773}}
  >{\centering\arraybackslash}m{(\columnwidth - 12\tabcolsep) * \real{0.1773}}
  >{\centering\arraybackslash}m{(\columnwidth - 12\tabcolsep) * \real{0.1383}}@{}}
\caption{Результаты анализа зависимости \(\sigma_{\max}\) в слоях композитных пластин при \(\sigma_{\max} = \sigma_{pl}\) от \(\rho_{rel}\) заполнителя в первой постановке численных экспериментов при \(d_{a} = 1\) и \(d_{a} = 1.3\) мм в условиях опирания с упругим поворотом}\label{tab:3}\\
\toprule\noalign{}
\multirow{2}{=}{\begin{minipage}[c]{\linewidth}\centering
\vspace{0.35\baselineskip}
\(d_{a}\)\\[-0.05\baselineskip](мм)
\end{minipage}} & \multirow{2}{=}{\begin{minipage}[c]{\linewidth}\centering
\vspace{0.35\baselineskip}
\(t_{sw}\)\\[-0.05\baselineskip](мм)
\end{minipage}} & \multirow{2}{=}{\begin{minipage}[c]{\linewidth}\centering
\vspace{0.35\baselineskip}
\(F_{y}\)\\[-0.05\baselineskip](Н)
\end{minipage}} & \multicolumn{3}{>{\centering\arraybackslash}m{(\columnwidth - 12\tabcolsep) * \real{0.4929} + 4\tabcolsep}}{
\begin{minipage}[c]{\linewidth}\centering
\(\sigma_{\max}\) (МПа)
\end{minipage}} & \multirow{2}{=}{\begin{minipage}[c]{\linewidth}\centering
\vspace{0.35\baselineskip}
\(\rho_{rel}\)\\[-0.05\baselineskip](\%)
\end{minipage}} \\
& & & \begin{minipage}[c]{\linewidth}\centering
\vspace{1.04\baselineskip}в сотах
\end{minipage} & \begin{minipage}[c]{\linewidth}\centering\setstretch{0.90}
\vspace{0.76\baselineskip}в верхнем сплошном слое
\end{minipage} & \begin{minipage}[c]{\linewidth}\centering\setstretch{0.90}
\vspace{0.76\baselineskip}в нижнем сплошном слое
\end{minipage} \\
\midrule\noalign{}
\endhead
\bottomrule\noalign{}
\endlastfoot
\tabledataC{=}{\midrule\noalign{}}
\end{longtable}
\fi
\endgroup
\clearpage

\begingroup
\setstretch{1.2}
\iflatexml
\begin{table}
\centering
\caption{Результаты анализа зависимости \(\sigma_{\max}\) в слоях композитных пластин при \(\sigma_{\max} = \sigma_{pl}\) от \(\rho_{rel}\) заполнителя в первой постановке численных экспериментов при \(d_{a} = 1.6\) и \(d_{a} = 1.9\) мм в условиях опирания с упругим поворотом}\label{tab:4}
\begin{tabular}{ccccccc}
\lxBeginTableHead
\lxTableColumnHead \(d_{a}\) (мм) & \lxTableColumnHead \(t_{sw}\) (мм) & \lxTableColumnHead \(F_{y}\) (Н) & \multicolumn{3}{c}{\lxTableColumnHead \(\sigma_{\max}\) (МПа)} & \lxTableColumnHead \(\rho_{rel}\) (\%) \\
& & & \lxTableColumnHead в сотах & \lxTableColumnHead в верхнем сплошном слое & \lxTableColumnHead в нижнем сплошном слое & \\
\lxEndTableHead
\toprule
\tabledataD{*}{\midrule}
\bottomrule
\end{tabular}
\end{table}
\else
\begin{longtable}[]{@{}
  >{\centering\arraybackslash}m{(\columnwidth - 12\tabcolsep) * \real{0.0923}}
  >{\centering\arraybackslash}m{(\columnwidth - 12\tabcolsep) * \real{0.1383}}
  >{\centering\arraybackslash}m{(\columnwidth - 12\tabcolsep) * \real{0.1383}}
  >{\centering\arraybackslash}m{(\columnwidth - 12\tabcolsep) * \real{0.1383}}
  >{\centering\arraybackslash}m{(\columnwidth - 12\tabcolsep) * \real{0.1773}}
  >{\centering\arraybackslash}m{(\columnwidth - 12\tabcolsep) * \real{0.1773}}
  >{\centering\arraybackslash}m{(\columnwidth - 12\tabcolsep) * \real{0.1383}}@{}}
\caption{Результаты анализа зависимости \(\sigma_{\max}\) в слоях композитных пластин при \(\sigma_{\max} = \sigma_{pl}\) от \(\rho_{rel}\) заполнителя в первой постановке численных экспериментов при \(d_{a} = 1.6\) и \(d_{a} = 1.9\) мм в условиях опирания с упругим поворотом}\label{tab:4}\\
\toprule\noalign{}
\multirow{2}{=}{\begin{minipage}[c]{\linewidth}\centering
\vspace{0.35\baselineskip}
\(d_{a}\)\\[-0.05\baselineskip](мм)
\end{minipage}} & \multirow{2}{=}{\begin{minipage}[c]{\linewidth}\centering
\vspace{0.35\baselineskip}
\(t_{sw}\)\\[-0.05\baselineskip](мм)
\end{minipage}} & \multirow{2}{=}{\begin{minipage}[c]{\linewidth}\centering
\vspace{0.35\baselineskip}
\(F_{y}\)\\[-0.05\baselineskip](Н)
\end{minipage}} & \multicolumn{3}{>{\centering\arraybackslash}m{(\columnwidth - 12\tabcolsep) * \real{0.4929} + 4\tabcolsep}}{
\begin{minipage}[c]{\linewidth}\centering
\(\sigma_{\max}\) (МПа)
\end{minipage}} & \multirow{2}{=}{\begin{minipage}[c]{\linewidth}\centering
\vspace{0.35\baselineskip}
\(\rho_{rel}\)\\[-0.05\baselineskip](\%)
\end{minipage}} \\
& & & \begin{minipage}[c]{\linewidth}\centering
\vspace{1.04\baselineskip}в сотах
\end{minipage} & \begin{minipage}[c]{\linewidth}\centering\setstretch{0.90}
\vspace{0.76\baselineskip}в верхнем сплошном слое
\end{minipage} & \begin{minipage}[c]{\linewidth}\centering\setstretch{0.90}
\vspace{0.76\baselineskip}в нижнем сплошном слое
\end{minipage} \\
\midrule\noalign{}
\endhead
\bottomrule\noalign{}
\endlastfoot
\tabledataD{=}{\midrule\noalign{}}
\end{longtable}
\fi
\endgroup
\clearpage

\begingroup
\setstretch{1.2}
\iflatexml
\begin{table}
\centering
\caption{Результаты анализа зависимости \(\sigma_{\max}\) в слоях композитных пластин при \(\sigma_{\max} = \sigma_{pl}\) от \(\rho_{rel}\) заполнителя во второй постановке численных экспериментов при \(d_{a} = 1\) и \(d_{a} = 1.3\) мм в условиях жесткого защемления}\label{tab:5}
\begin{tabular}{ccccccc}
\lxBeginTableHead
\lxTableColumnHead \(d_{a}\) (мм) & \lxTableColumnHead \(t_{sw}\) (мм) & \lxTableColumnHead \(F_{y}\) (Н) & \multicolumn{3}{c}{\lxTableColumnHead \(\sigma_{\max}\) (МПа)} & \lxTableColumnHead \(\rho_{rel}\) (\%) \\
& & & \lxTableColumnHead в сотах & \lxTableColumnHead в верхнем сплошном слое & \lxTableColumnHead в нижнем сплошном слое & \\
\lxEndTableHead
\toprule
\tabledataE{*}{\midrule}
\bottomrule
\end{tabular}
\end{table}
\else
\begin{longtable}[]{@{}
  >{\centering\arraybackslash}m{(\columnwidth - 12\tabcolsep) * \real{0.0923}}
  >{\centering\arraybackslash}m{(\columnwidth - 12\tabcolsep) * \real{0.1383}}
  >{\centering\arraybackslash}m{(\columnwidth - 12\tabcolsep) * \real{0.1383}}
  >{\centering\arraybackslash}m{(\columnwidth - 12\tabcolsep) * \real{0.1383}}
  >{\centering\arraybackslash}m{(\columnwidth - 12\tabcolsep) * \real{0.1773}}
  >{\centering\arraybackslash}m{(\columnwidth - 12\tabcolsep) * \real{0.1773}}
  >{\centering\arraybackslash}m{(\columnwidth - 12\tabcolsep) * \real{0.1383}}@{}}
\caption{Результаты анализа зависимости \(\sigma_{\max}\) в слоях композитных пластин при \(\sigma_{\max} = \sigma_{pl}\) от \(\rho_{rel}\) заполнителя во второй постановке численных экспериментов при \(d_{a} = 1\) и \(d_{a} = 1.3\) мм в условиях жесткого защемления}\label{tab:5}\\
\toprule\noalign{}
\multirow{2}{=}{\begin{minipage}[c]{\linewidth}\centering
\vspace{0.35\baselineskip}
\(d_{a}\)\\[-0.05\baselineskip](мм)
\end{minipage}} & \multirow{2}{=}{\begin{minipage}[c]{\linewidth}\centering
\vspace{0.35\baselineskip}
\(t_{sw}\)\\[-0.05\baselineskip](мм)
\end{minipage}} & \multirow{2}{=}{\begin{minipage}[c]{\linewidth}\centering
\vspace{0.35\baselineskip}
\(F_{y}\)\\[-0.05\baselineskip](Н)
\end{minipage}} & \multicolumn{3}{>{\centering\arraybackslash}m{(\columnwidth - 12\tabcolsep) * \real{0.4929} + 4\tabcolsep}}{
\begin{minipage}[c]{\linewidth}\centering
\(\sigma_{\max}\) (МПа)
\end{minipage}} & \multirow{2}{=}{\begin{minipage}[c]{\linewidth}\centering
\vspace{0.35\baselineskip}
\(\rho_{rel}\)\\[-0.05\baselineskip](\%)
\end{minipage}} \\
& & & \begin{minipage}[c]{\linewidth}\centering
\vspace{1.04\baselineskip}в сотах
\end{minipage} & \begin{minipage}[c]{\linewidth}\centering\setstretch{0.90}
\vspace{0.76\baselineskip}в верхнем сплошном слое
\end{minipage} & \begin{minipage}[c]{\linewidth}\centering\setstretch{0.90}
\vspace{0.76\baselineskip}в нижнем сплошном слое
\end{minipage} \\
\midrule\noalign{}
\endhead
\bottomrule\noalign{}
\endlastfoot
\tabledataE{=}{\midrule\noalign{}}
\end{longtable}
\fi
\endgroup
\clearpage

\begingroup
\setstretch{1.2}
\iflatexml
\begin{table}
\centering
\caption{Результаты анализа зависимости \(\sigma_{\max}\) в слоях композитных пластин при \(\sigma_{\max} = \sigma_{pl}\) от \(\rho_{rel}\) заполнителя во второй постановке численных экспериментов при \(d_{a} = 1.6\) и \(d_{a} = 1.9\) мм в условиях жесткого защемления}\label{tab:6}
\begin{tabular}{ccccccc}
\lxBeginTableHead
\lxTableColumnHead \(d_{a}\) (мм) & \lxTableColumnHead \(t_{sw}\) (мм) & \lxTableColumnHead \(F_{y}\) (Н) & \multicolumn{3}{c}{\lxTableColumnHead \(\sigma_{\max}\) (МПа)} & \lxTableColumnHead \(\rho_{rel}\) (\%) \\
& & & \lxTableColumnHead в сотах & \lxTableColumnHead в верхнем сплошном слое & \lxTableColumnHead в нижнем сплошном слое & \\
\lxEndTableHead
\toprule
\tabledataF{*}{\midrule}
\bottomrule
\end{tabular}
\end{table}
\else
\begin{longtable}[]{@{}
  >{\centering\arraybackslash}m{(\columnwidth - 12\tabcolsep) * \real{0.0923}}
  >{\centering\arraybackslash}m{(\columnwidth - 12\tabcolsep) * \real{0.1383}}
  >{\centering\arraybackslash}m{(\columnwidth - 12\tabcolsep) * \real{0.1383}}
  >{\centering\arraybackslash}m{(\columnwidth - 12\tabcolsep) * \real{0.1383}}
  >{\centering\arraybackslash}m{(\columnwidth - 12\tabcolsep) * \real{0.1773}}
  >{\centering\arraybackslash}m{(\columnwidth - 12\tabcolsep) * \real{0.1773}}
  >{\centering\arraybackslash}m{(\columnwidth - 12\tabcolsep) * \real{0.1383}}@{}}
\caption{Результаты анализа зависимости \(\sigma_{\max}\) в слоях композитных пластин при \(\sigma_{\max} = \sigma_{pl}\) от \(\rho_{rel}\) заполнителя во второй постановке численных экспериментов при \(d_{a} = 1.6\) и \(d_{a} = 1.9\) мм в условиях жесткого защемления}\label{tab:6}\\
\toprule\noalign{}
\multirow{2}{=}{\begin{minipage}[c]{\linewidth}\centering
\vspace{0.35\baselineskip}
\(d_{a}\)\\[-0.05\baselineskip](мм)
\end{minipage}} & \multirow{2}{=}{\begin{minipage}[c]{\linewidth}\centering
\vspace{0.35\baselineskip}
\(t_{sw}\)\\[-0.05\baselineskip](мм)
\end{minipage}} & \multirow{2}{=}{\begin{minipage}[c]{\linewidth}\centering
\vspace{0.35\baselineskip}
\(F_{y}\)\\[-0.05\baselineskip](Н)
\end{minipage}} & \multicolumn{3}{>{\centering\arraybackslash}m{(\columnwidth - 12\tabcolsep) * \real{0.4929} + 4\tabcolsep}}{
\begin{minipage}[c]{\linewidth}\centering
\(\sigma_{\max}\) (МПа)
\end{minipage}} & \multirow{2}{=}{\begin{minipage}[c]{\linewidth}\centering
\vspace{0.35\baselineskip}
\(\rho_{rel}\)\\[-0.05\baselineskip](\%)
\end{minipage}} \\
& & & \begin{minipage}[c]{\linewidth}\centering
\vspace{1.04\baselineskip}в сотах
\end{minipage} & \begin{minipage}[c]{\linewidth}\centering\setstretch{0.90}
\vspace{0.76\baselineskip}в верхнем сплошном слое
\end{minipage} & \begin{minipage}[c]{\linewidth}\centering\setstretch{0.90}
\vspace{0.76\baselineskip}в нижнем сплошном слое
\end{minipage} \\
\midrule\noalign{}
\endhead
\bottomrule\noalign{}
\endlastfoot
\tabledataF{=}{\midrule\noalign{}}
\end{longtable}
\fi
\endgroup
\clearpage

\begingroup
\setstretch{1.2}
\iflatexml
\begin{table}
\centering
\caption{Результаты анализа зависимости \(\sigma_{\max}\) в слоях композитных пластин при \(\sigma_{\max} = \sigma_{pl}\) от \(\rho_{rel}\) заполнителя во второй постановке численных экспериментов при \(d_{a} = 1\) и \(d_{a} = 1.3\) мм в условиях опирания с упругим поворотом}\label{tab:7}
\begin{tabular}{ccccccc}
\lxBeginTableHead
\lxTableColumnHead \(d_{a}\) (мм) & \lxTableColumnHead \(t_{sw}\) (мм) & \lxTableColumnHead \(F_{y}\) (Н) & \multicolumn{3}{c}{\lxTableColumnHead \(\sigma_{\max}\) (МПа)} & \lxTableColumnHead \(\rho_{rel}\) (\%) \\
& & & \lxTableColumnHead в сотах & \lxTableColumnHead в верхнем сплошном слое & \lxTableColumnHead в нижнем сплошном слое & \\
\lxEndTableHead
\toprule
\tabledataG{*}{\midrule}
\bottomrule
\end{tabular}
\end{table}
\else
\begin{longtable}[]{@{}
  >{\centering\arraybackslash}m{(\columnwidth - 12\tabcolsep) * \real{0.0923}}
  >{\centering\arraybackslash}m{(\columnwidth - 12\tabcolsep) * \real{0.1383}}
  >{\centering\arraybackslash}m{(\columnwidth - 12\tabcolsep) * \real{0.1383}}
  >{\centering\arraybackslash}m{(\columnwidth - 12\tabcolsep) * \real{0.1383}}
  >{\centering\arraybackslash}m{(\columnwidth - 12\tabcolsep) * \real{0.1773}}
  >{\centering\arraybackslash}m{(\columnwidth - 12\tabcolsep) * \real{0.1773}}
  >{\centering\arraybackslash}m{(\columnwidth - 12\tabcolsep) * \real{0.1383}}@{}}
\caption{Результаты анализа зависимости \(\sigma_{\max}\) в слоях композитных пластин при \(\sigma_{\max} = \sigma_{pl}\) от \(\rho_{rel}\) заполнителя во второй постановке численных экспериментов при \(d_{a} = 1\) и \(d_{a} = 1.3\) мм в условиях опирания с упругим поворотом}\label{tab:7}\\
\toprule\noalign{}
\multirow{2}{=}{\begin{minipage}[c]{\linewidth}\centering
\vspace{0.35\baselineskip}
\(d_{a}\)\\[-0.05\baselineskip](мм)
\end{minipage}} & \multirow{2}{=}{\begin{minipage}[c]{\linewidth}\centering
\vspace{0.35\baselineskip}
\(t_{sw}\)\\[-0.05\baselineskip](мм)
\end{minipage}} & \multirow{2}{=}{\begin{minipage}[c]{\linewidth}\centering
\vspace{0.35\baselineskip}
\(F_{y}\)\\[-0.05\baselineskip](Н)
\end{minipage}} & \multicolumn{3}{>{\centering\arraybackslash}m{(\columnwidth - 12\tabcolsep) * \real{0.4929} + 4\tabcolsep}}{
\begin{minipage}[c]{\linewidth}\centering
\(\sigma_{\max}\) (МПа)
\end{minipage}} & \multirow{2}{=}{\begin{minipage}[c]{\linewidth}\centering
\vspace{0.35\baselineskip}
\(\rho_{rel}\)\\[-0.05\baselineskip](\%)
\end{minipage}} \\
& & & \begin{minipage}[c]{\linewidth}\centering
\vspace{1.04\baselineskip}в сотах
\end{minipage} & \begin{minipage}[c]{\linewidth}\centering\setstretch{0.90}
\vspace{0.76\baselineskip}в верхнем сплошном слое
\end{minipage} & \begin{minipage}[c]{\linewidth}\centering\setstretch{0.90}
\vspace{0.76\baselineskip}в нижнем сплошном слое
\end{minipage} \\
\midrule\noalign{}
\endhead
\bottomrule\noalign{}
\endlastfoot
\tabledataG{=}{\midrule\noalign{}}
\end{longtable}
\fi
\endgroup
\clearpage

\begingroup
\setstretch{1.2}
\iflatexml
\begin{table}
\centering
\caption{Результаты анализа зависимости \(\sigma_{\max}\) в слоях композитных пластин при \(\sigma_{\max} = \sigma_{pl}\) от \(\rho_{rel}\) заполнителя во второй постановке численных экспериментов при \(d_{a} = 1.6\) и \(d_{a} = 1.9\) мм в условиях опирания с упругим поворотом}\label{tab:8}
\begin{tabular}{ccccccc}
\lxBeginTableHead
\lxTableColumnHead \(d_{a}\) (мм) & \lxTableColumnHead \(t_{sw}\) (мм) & \lxTableColumnHead \(F_{y}\) (Н) & \multicolumn{3}{c}{\lxTableColumnHead \(\sigma_{\max}\) (МПа)} & \lxTableColumnHead \(\rho_{rel}\) (\%) \\
& & & \lxTableColumnHead в сотах & \lxTableColumnHead в верхнем сплошном слое & \lxTableColumnHead в нижнем сплошном слое & \\
\lxEndTableHead
\toprule
\tabledataH{*}{\midrule}
\bottomrule
\end{tabular}
\end{table}
\else
\begin{longtable}[]{@{}
  >{\centering\arraybackslash}m{(\columnwidth - 12\tabcolsep) * \real{0.0923}}
  >{\centering\arraybackslash}m{(\columnwidth - 12\tabcolsep) * \real{0.1383}}
  >{\centering\arraybackslash}m{(\columnwidth - 12\tabcolsep) * \real{0.1383}}
  >{\centering\arraybackslash}m{(\columnwidth - 12\tabcolsep) * \real{0.1383}}
  >{\centering\arraybackslash}m{(\columnwidth - 12\tabcolsep) * \real{0.1773}}
  >{\centering\arraybackslash}m{(\columnwidth - 12\tabcolsep) * \real{0.1773}}
  >{\centering\arraybackslash}m{(\columnwidth - 12\tabcolsep) * \real{0.1383}}@{}}
\caption{Результаты анализа зависимости \(\sigma_{\max}\) в слоях композитных пластин при \(\sigma_{\max} = \sigma_{pl}\) от \(\rho_{rel}\) заполнителя во второй постановке численных экспериментов при \(d_{a} = 1.6\) и \(d_{a} = 1.9\) мм в условиях опирания с упругим поворотом}\label{tab:8}\\
\toprule\noalign{}
\multirow{2}{=}{\begin{minipage}[c]{\linewidth}\centering
\vspace{0.35\baselineskip}
\(d_{a}\)\\[-0.05\baselineskip](мм)
\end{minipage}} & \multirow{2}{=}{\begin{minipage}[c]{\linewidth}\centering
\vspace{0.35\baselineskip}
\(t_{sw}\)\\[-0.05\baselineskip](мм)
\end{minipage}} & \multirow{2}{=}{\begin{minipage}[c]{\linewidth}\centering
\vspace{0.35\baselineskip}
\(F_{y}\)\\[-0.05\baselineskip](Н)
\end{minipage}} & \multicolumn{3}{>{\centering\arraybackslash}m{(\columnwidth - 12\tabcolsep) * \real{0.4929} + 4\tabcolsep}}{
\begin{minipage}[c]{\linewidth}\centering
\(\sigma_{\max}\) (МПа)
\end{minipage}} & \multirow{2}{=}{\begin{minipage}[c]{\linewidth}\centering
\vspace{0.35\baselineskip}
\(\rho_{rel}\)\\[-0.05\baselineskip](\%)
\end{minipage}} \\
& & & \begin{minipage}[c]{\linewidth}\centering
\vspace{1.04\baselineskip}в сотах
\end{minipage} & \begin{minipage}[c]{\linewidth}\centering\setstretch{0.90}
\vspace{0.76\baselineskip}в верхнем сплошном слое
\end{minipage} & \begin{minipage}[c]{\linewidth}\centering\setstretch{0.90}
\vspace{0.76\baselineskip}в нижнем сплошном слое
\end{minipage} \\
\midrule\noalign{}
\endhead
\bottomrule\noalign{}
\endlastfoot
\tabledataH{=}{\midrule\noalign{}}
\end{longtable}
\fi
\endgroup
\clearpage

\vspace{-0.575\baselineskip}
\vspace{-0.575\baselineskip}
\vspace{-0.575\baselineskip}
\vspace{-0.575\baselineskip}
\vspace{-0.575\baselineskip}
\vspace{-0.575\baselineskip}
\vspace{-0.575\baselineskip}
\vspace{-0.575\baselineskip}

\iflatexml\else\phantomsection\fi
\subsection{3.4. Результаты лабораторных испытаний трехслойных композитов с тетракиральным сотовым заполнителем, полученных с применением аддитивных технологий}

\begin{figure}[H]
\centering
\includegraphics[keepaspectratio,width=5.47244in,height=3.64432in]{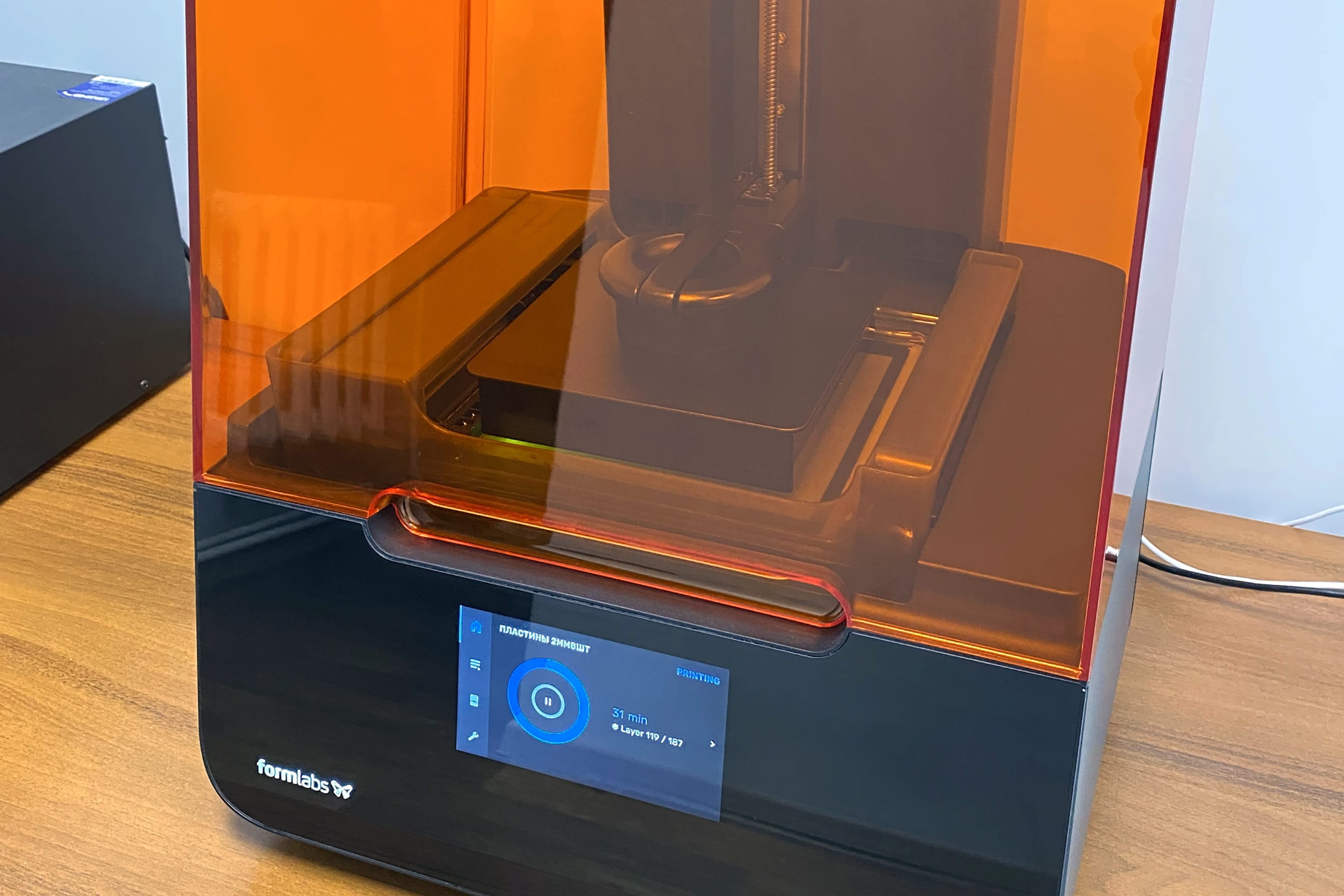}
\caption{Печать сотовых заполнителей и внешних слоев на 3D-принтере \dualcite{ref78}{{[}\bibnum{ref78}{]}}{{[}\bibnum{ref78}{]}}}\label{fig:40}
\end{figure}

\begin{figure}[H]
\centering
\vspace{0.38\baselineskip}
\includegraphics[keepaspectratio,width=5.47244in,height=3.64432in]{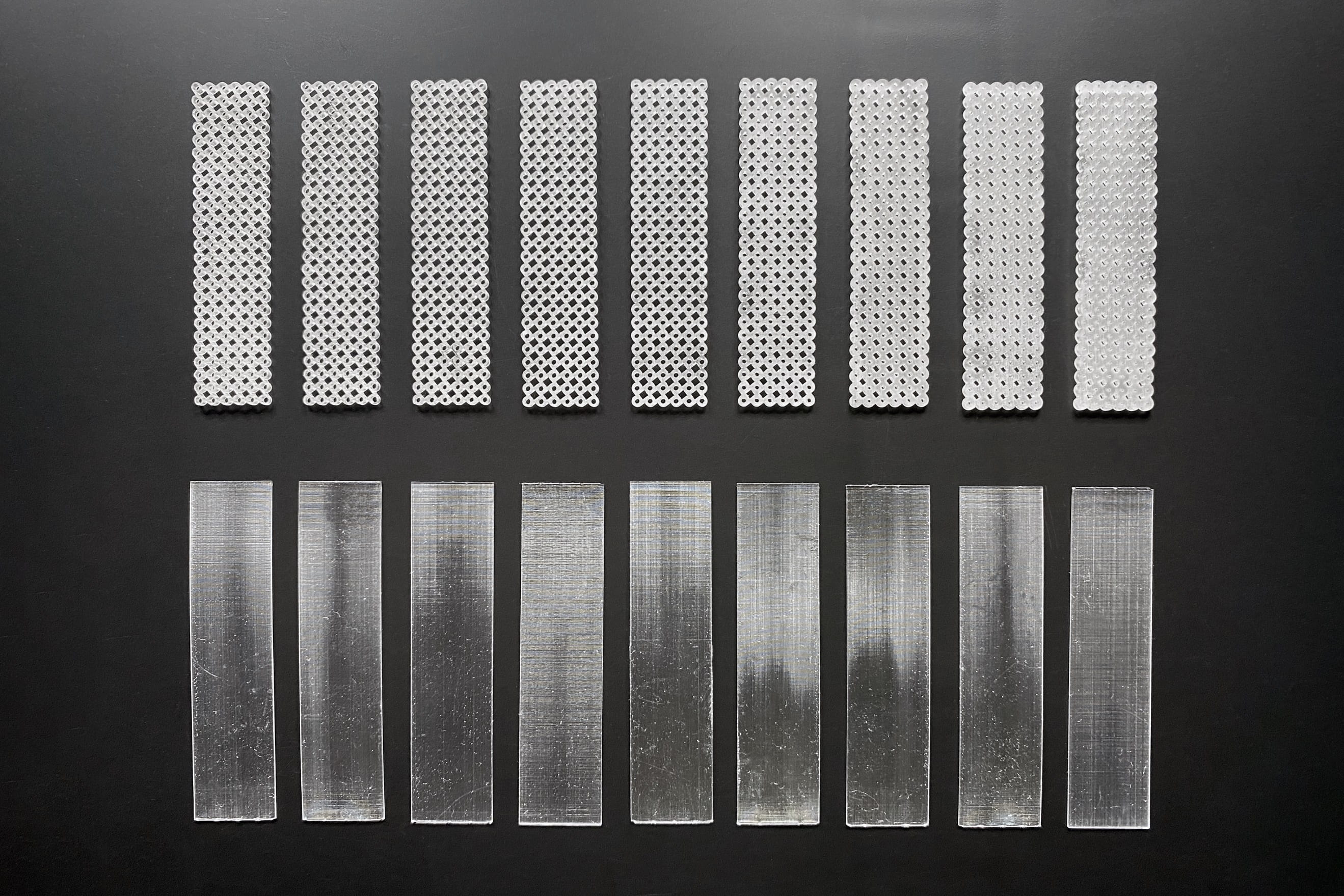}
\caption{Тетракиральные соты и сплошные слои для сборки композитов \dualcite{ref78}{{[}\bibnum{ref78}{]}}{{[}\bibnum{ref78}{]}}}\label{fig:41}
\end{figure}

\begin{figure}[H]
\centering
\includegraphics[keepaspectratio,width=5.47244in,height=3.64432in]{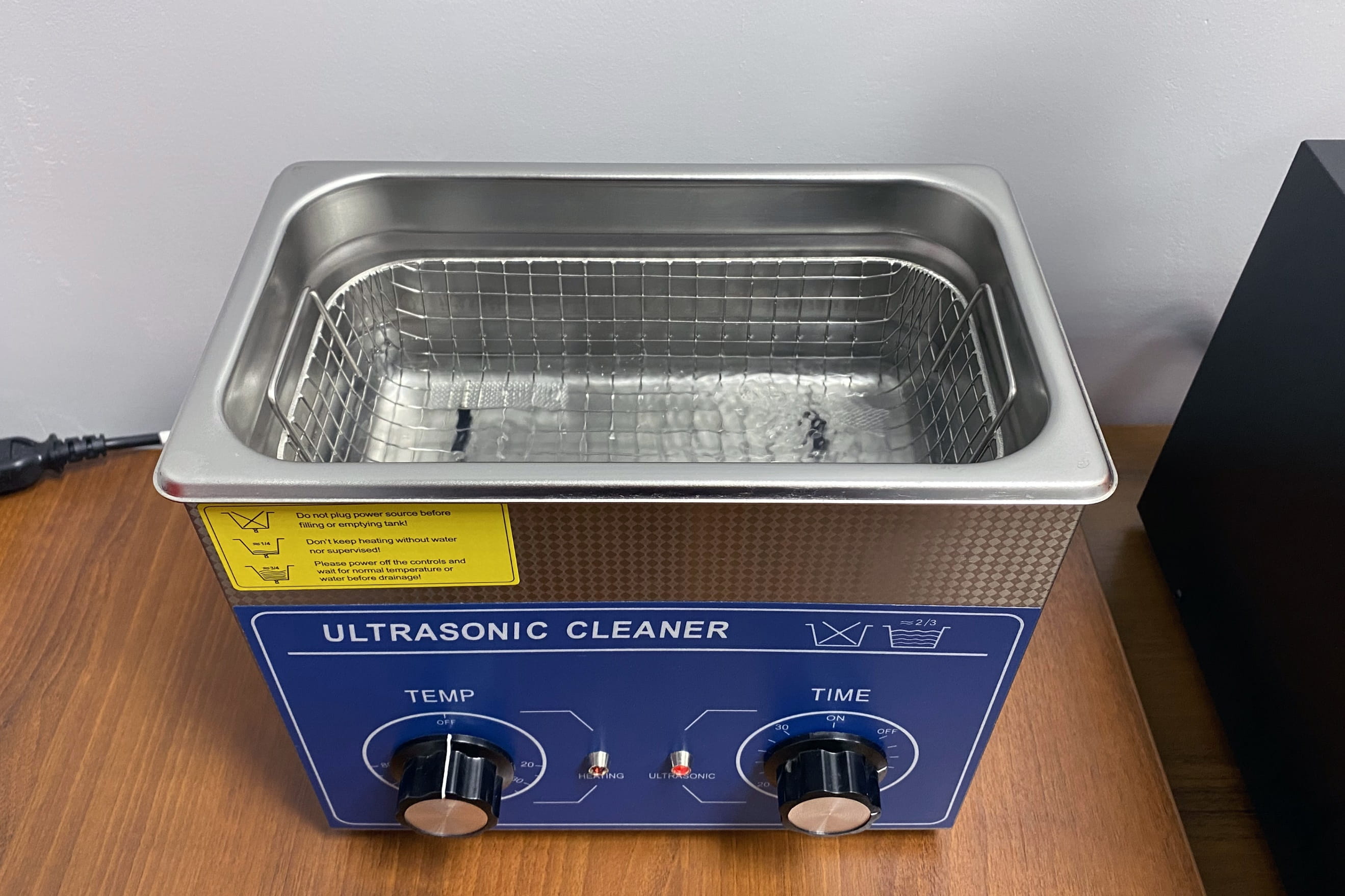}
\caption{Промывка сот и внешних слоев в ультразвуковой ванне \dualcite{ref78}{{[}\bibnum{ref78}{]}}{{[}\bibnum{ref78}{]}}}\label{fig:42}
\end{figure}

\begin{figure}[H]
\centering
\vspace{0.38\baselineskip}
\includegraphics[keepaspectratio,width=5.47244in,height=3.65in]{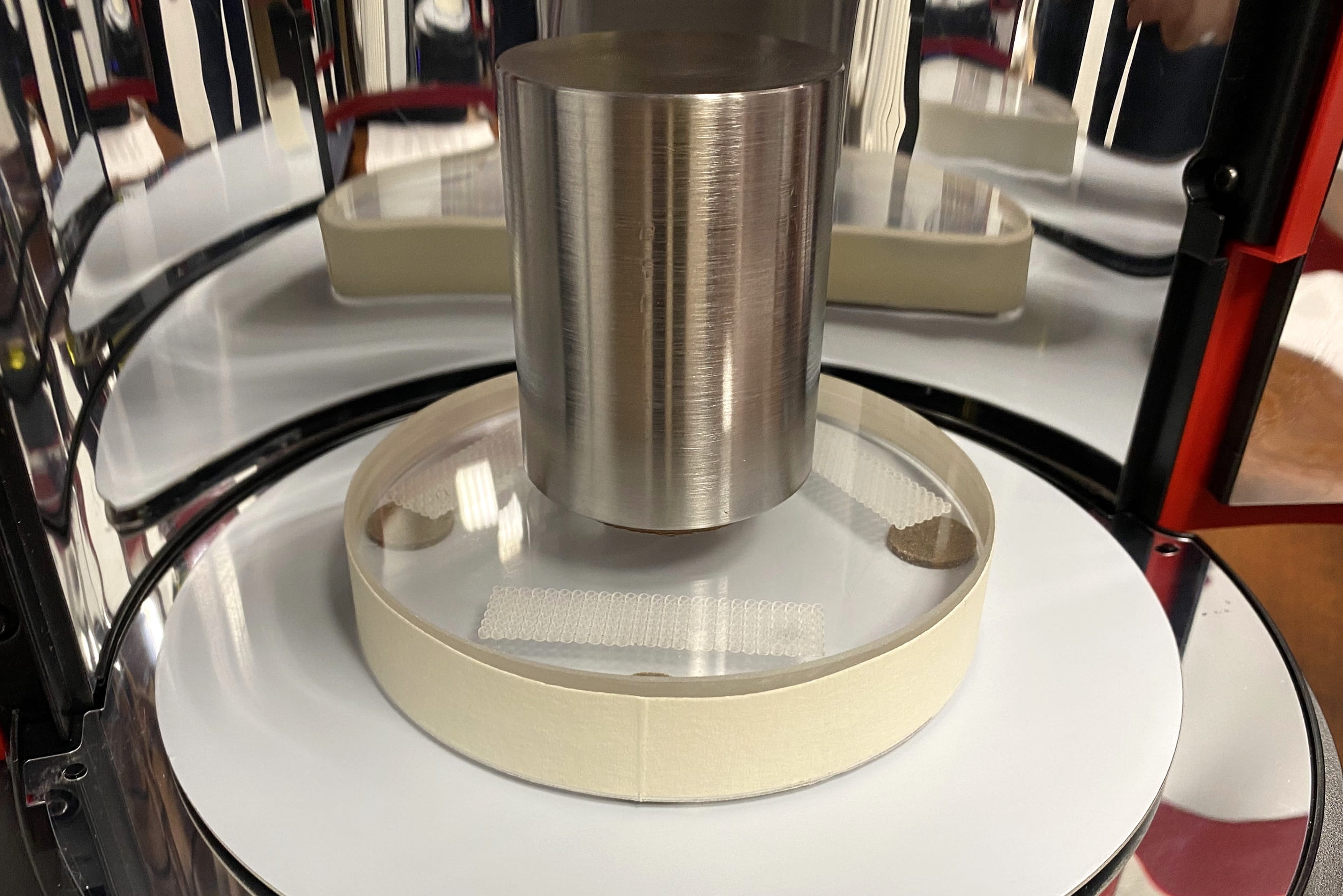}
\caption{Полимеризация сотовых заполнителей в ультрафиолетовой камере \dualcite{ref78}{{[}\bibnum{ref78}{]}}{{[}\bibnum{ref78}{]}}}\label{fig:43}
\end{figure}

\begin{figure}[H]
\centering
\includegraphics[keepaspectratio,width=5.47244in,height=3.64432in]{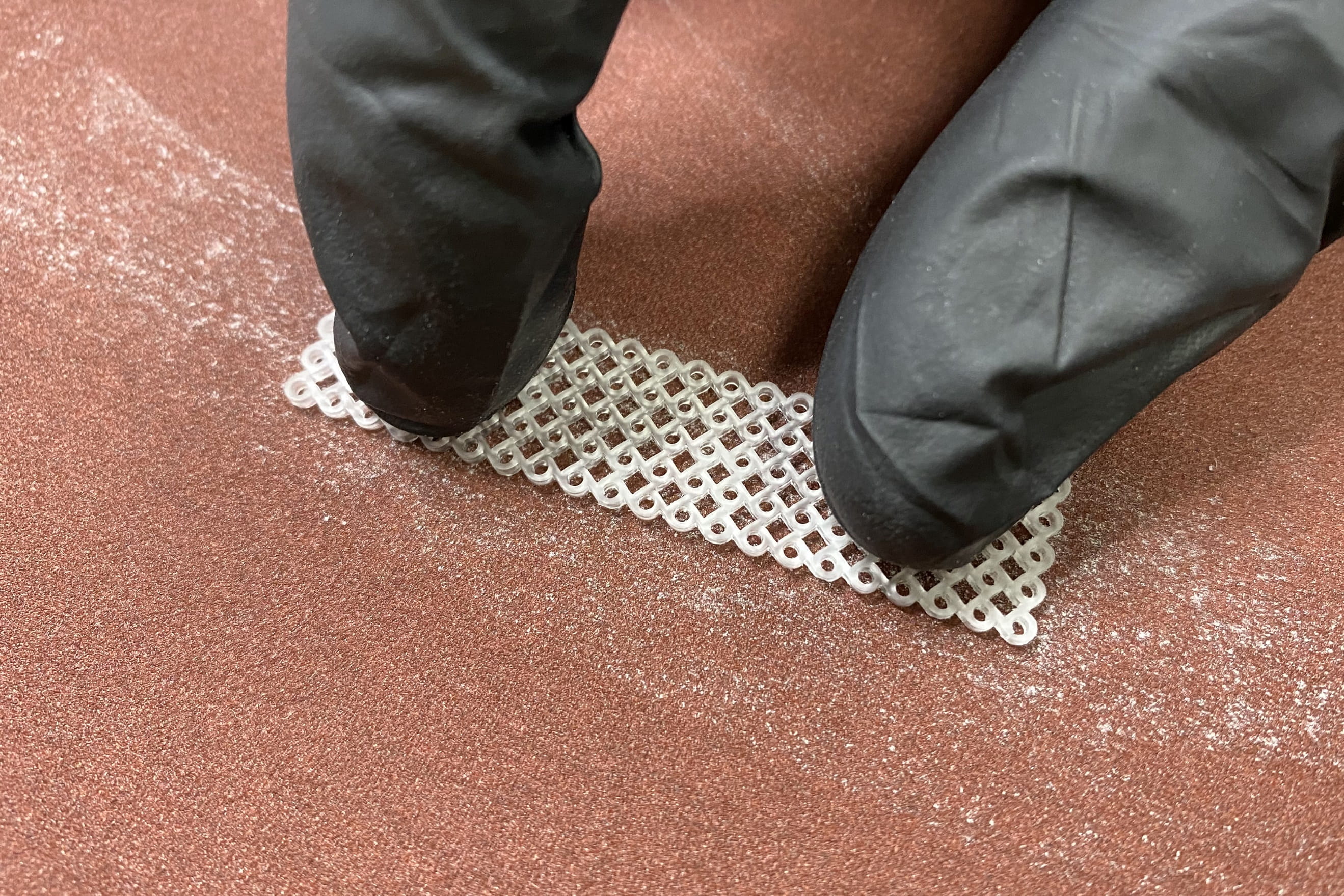}
\caption{Шлифовка сотовых заполнителей \dualcite{ref78}{{[}\bibnum{ref78}{]}}{{[}\bibnum{ref78}{]}}}\label{fig:44}
\end{figure}

\begin{figure}[H]
\centering
\vspace{0.38\baselineskip}
\includegraphics[keepaspectratio,width=5.47244in,height=3.64432in]{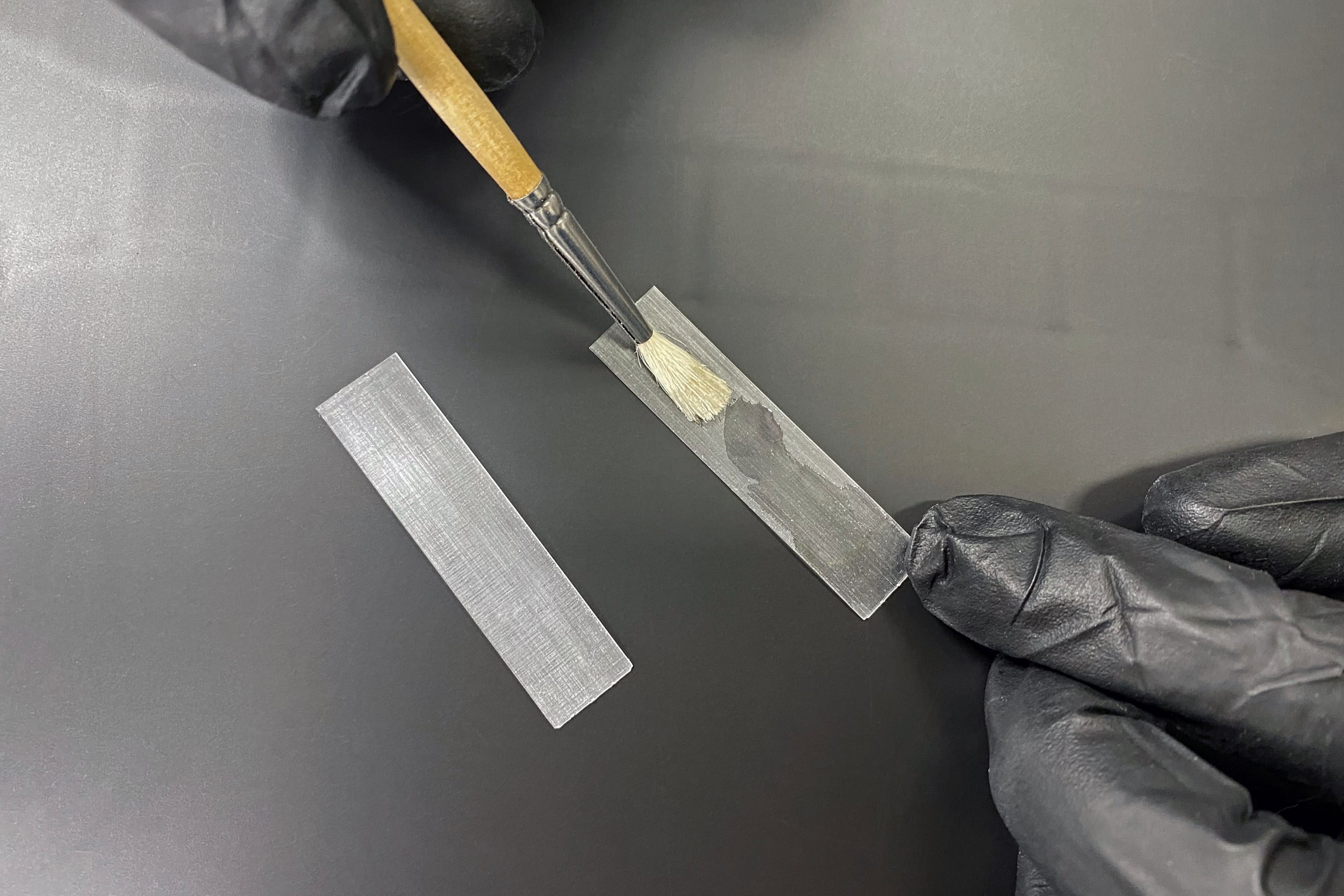}
\caption{Нанесение кистью клеевого слоя \dualcite{ref78}{{[}\bibnum{ref78}{]}}{{[}\bibnum{ref78}{]}}}\label{fig:45}
\end{figure}

Серия композитных пластин согласно первой постановке экспериментов была изготовлена с применением аддитивных технологий в лаборатории. В качестве материала использовалась фотополимерная смола Formlabs Clear Resin \dualcite{ref81}{{[}\bibnum{ref81}{]}}{{[}\bibnum{ref81}{]}}, а печать слоев композита производилась на 3D-принтере Formlabs Form 3. Технологический процесс изготовления экспериментальных образцов состоял из следующих этапов \dualcite{ref78}{{[}\bibnum{ref78}{]}}{{[}\bibnum{ref78}{]}}:

\iflatexml
\renewcommand{\labelenumi}{\arabic{enumi}.}
\begin{enumerate}
\fi
\fullwidthnumitem{1}{Печать сотовых заполнителей и внешних слоев на 3D-принтере Formlabs Form 3 (рисунки \ref{fig:40}, \ref{fig:41}) при толщине слоя печати 100 мкм;}
\fullwidthnumitem{2}{Промывка напечатанных слоев изопропиловым спиртом в ультразвуковой ванне Uniz UC-4120 (рисунок \ref{fig:42}) в течение 5 минут;}
\fullwidthnumitem{3}{Предварительное отверждение сотовых заполнителей и внешних слоев между двумя кварцевыми стеклами диаметром 110 мм и толщиной 10 мм в ультрафиолетовой камере XYZPrinting (рисунок \ref{fig:43}) при длине волны \(\lambda = 405\) нм в течение 5 минут. Одновременно отверждались три слоя одинаковой толщины (либо сот, либо внешних слоев) с использованием дополнительного груза массой 1 кг (масса стекла -- 200 граммов);}
\fullwidthnumitem{4}{Шлифовка сотовых заполнителей с запасом толщины наждачной бумагой (рисунок \ref{fig:44}) с зернистостью P240 (получение плоской поверхности) и P1000 (финишная шлифовка);}
\fullwidthnumitem{5}{Промывка сот изопропиловым спиртом в ультразвуковой ванне Uniz UC-4120 в течение 5 минут;}
\fullwidthnumitem{6}{Нанесение кистью клеевого слоя Formlabs Clear Resin на внешние слои (рисунок \ref{fig:45}) с последующей сборкой композита. При этом толщина сот -- 1.2 мм, толщина внешних слоев -- 0.4 мм, а толщина клеевого слоя \(\leq\) 0.1 мм;}
\fullwidthnumitem{7}{Окончательное отверждение собранного композита между кварцевыми стеклами с дополнительным грузом в полимеризационной камере XYZPrinting при длине волны \(\lambda = 405\) нм в течение 5 минут.}
\iflatexml
\end{enumerate}
\fi

Традиционные способы изготовления слоистых композитов предполагают сборку слоев с применением различных клеев или припоев \dualcite{ref9,ref11,ref22}{{[}\bibnum{ref9}, \bibnum{ref11}, \bibnum{ref22}{]}}{{[}\bibnum{ref9}, \bibnum{ref11}, \bibnum{ref22}{]}}. При этом композиты со значительным отличием коэффициента Пуассона у слоев подвержены раннему расслоению вследствие внешних воздействий \dualcite{ref3}{{[}\bibnum{ref3}{]}}{{[}\bibnum{ref3}{]}}. В приведенном выше техпроцессе предусмотрено изготовление слоистой пластины из одного материала, что приводит к образованию практически монолитной конструкции и снижает риск расслоения.

\begin{figure}[H]
\centering
\includegraphics[keepaspectratio,width=5.51181in,height=3.78255in]{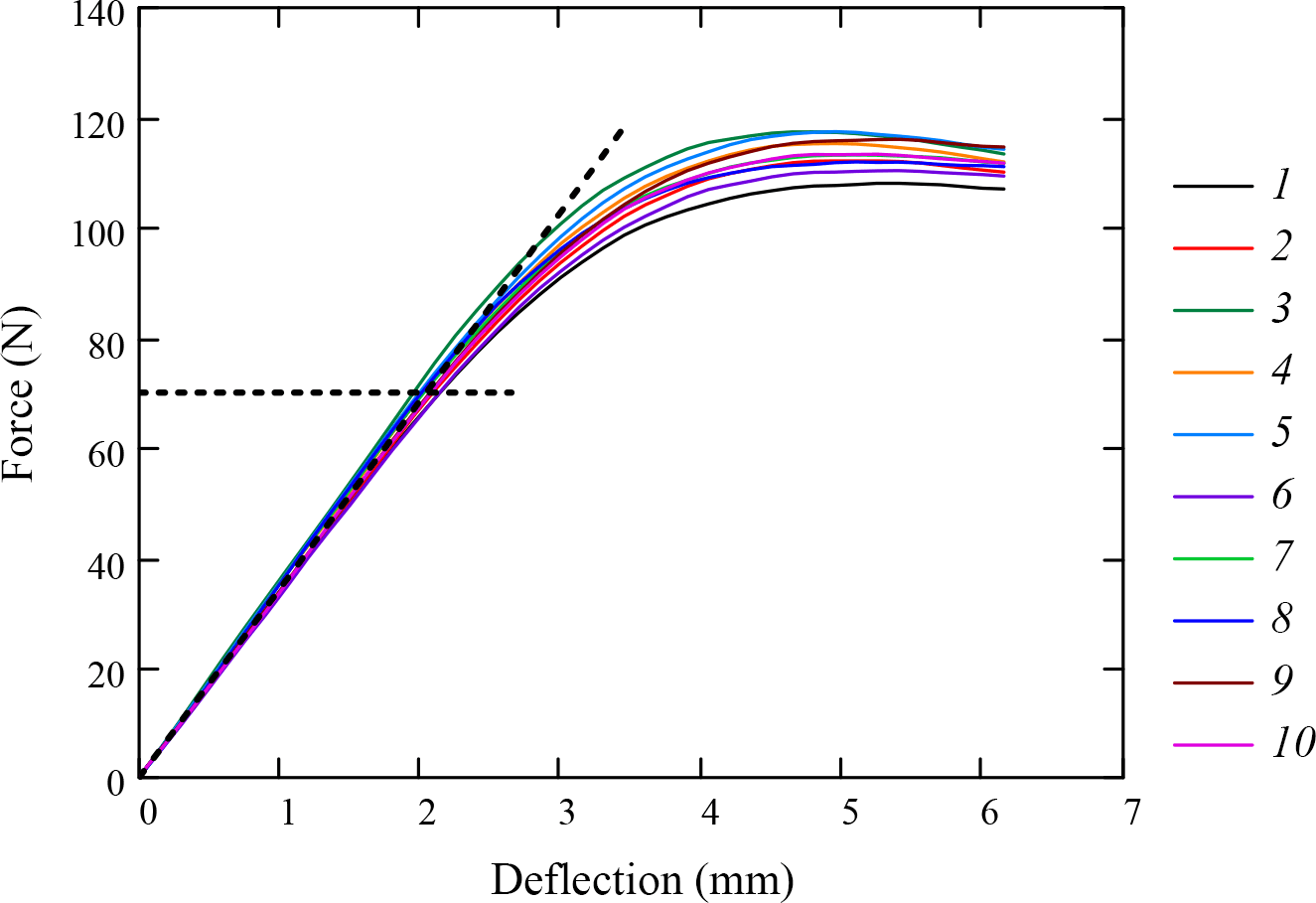}
\caption{Диаграмма зависимости \(F_{y}\) от прогиба сплошных пластин при трехточечном изгибе \dualcite{ref78}{{[}\bibnum{ref78}{]}}{{[}\bibnum{ref78}{]}}}\label{fig:46}
\end{figure}

Поскольку в численном моделировании использовались параметры материала Formlabs Clear из спецификации \dualcite{ref81}{{[}\bibnum{ref81}{]}}{{[}\bibnum{ref81}{]}}, выбор суммарного времени отверждения композитов осуществлялся посредством сопоставления модуля упругости при изгибе пластин со значением, приведенным в \dualcite{ref81}{{[}\bibnum{ref81}{]}}{{[}\bibnum{ref81}{]}}. Испытания сплошных пластин с размерами \(54 \times 13 \times 2\) мм производились с применением электромеханической системы Instron 5982 при частичном соответствии стандарту ASTM D790-10 \dualcite{ref44}{{[}\bibnum{ref44}{]}}{{[}\bibnum{ref44}{]}}. При этом радиус пуансона и опор -- 5 и 2.5 мм соответственно, расстояние между опорами -- 30 мм, а скорость опускания пуансона -- 4 мм/мин. На рисунке \ref{fig:46} представлена диаграмма зависимости нагрузки от прогиба для 10 образцов сплошной пластины, которые отверждались между кварцевыми стеклами с дополнительным грузом (согласно техпроцессу) в течение 10 мин (кривые для каждого образца показаны разным цветом). В случае полимеризации пластин в течение 10 минут среднее значение модуля упругости материала при изгибе составляет 2.3 ГПа (рисунок \ref{fig:46}), что хорошо соответствует значению из спецификации -- 2.2 ГПа \dualcite{ref81}{{[}\bibnum{ref81}{]}}{{[}\bibnum{ref81}{]}}. Также осуществлялось сопоставление значений силы, при которых максимальные напряжения в сплошных пластинах достигали предела пропорциональности. По результатам лабораторных испытаний среднее значение силы хорошо согласуется с расчетным значением, полученным с применением Comsol Multiphysics.

\newpage
\vspace{-0.575\baselineskip}
\vspace{-0.575\baselineskip}
\vspace{-0.575\baselineskip}
\vspace{-0.575\baselineskip}
\vspace{-0.575\baselineskip}
\vspace{-0.575\baselineskip}

\begin{figure}[H]
\centering
\includegraphics[keepaspectratio,width=6.29921in,height=3.76142in]{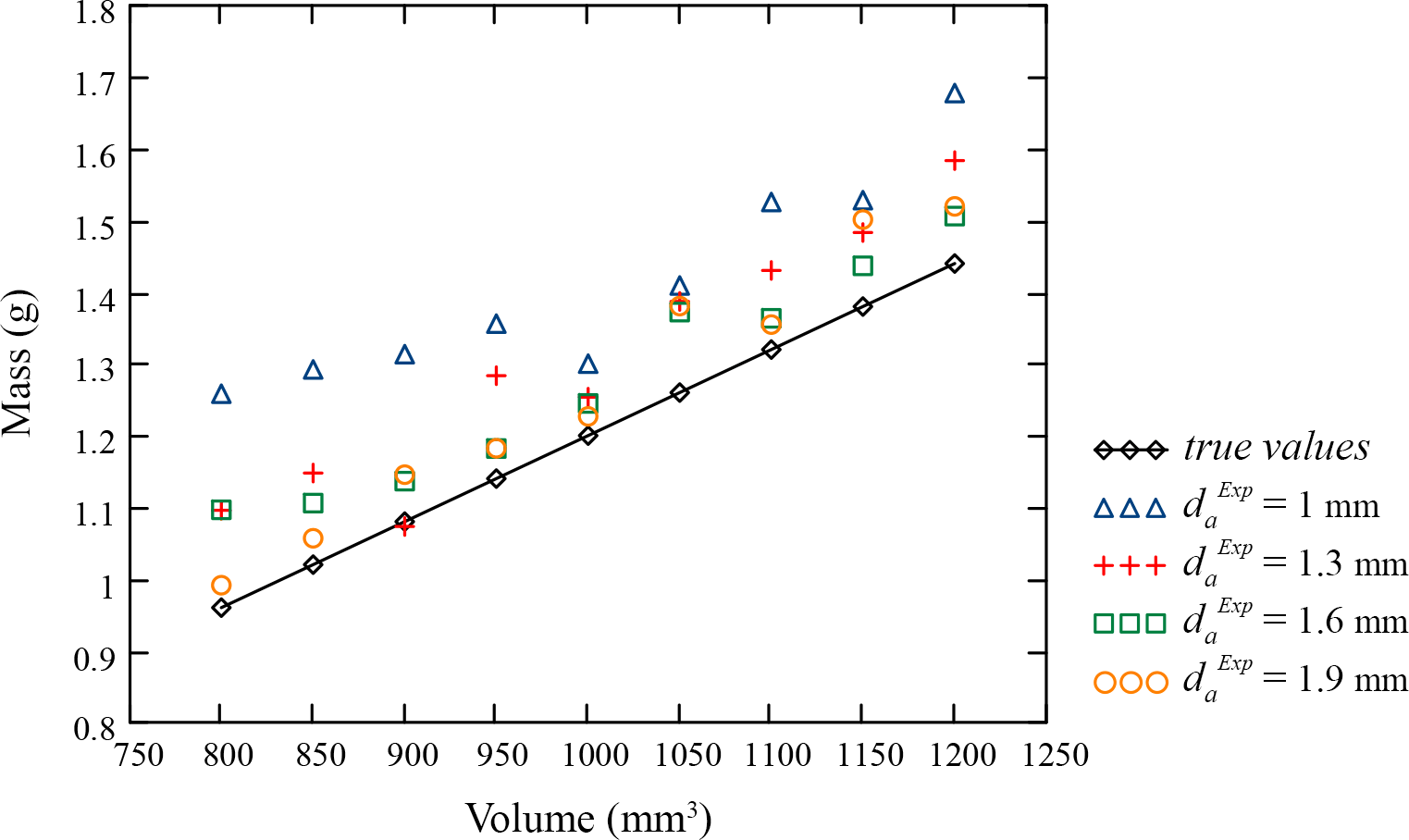}
\caption{Отклонение фактической массы изготовленных композитов от номинальных значений по CAD-модели \dualcite{ref78}{{[}\bibnum{ref78}{]}}{{[}\bibnum{ref78}{]}}}\label{fig:47}
\end{figure}

\begin{figure}[H]
\centering
\vspace{0.38\baselineskip}
\includegraphics[keepaspectratio,width=6.69291in,height=2.99193in]{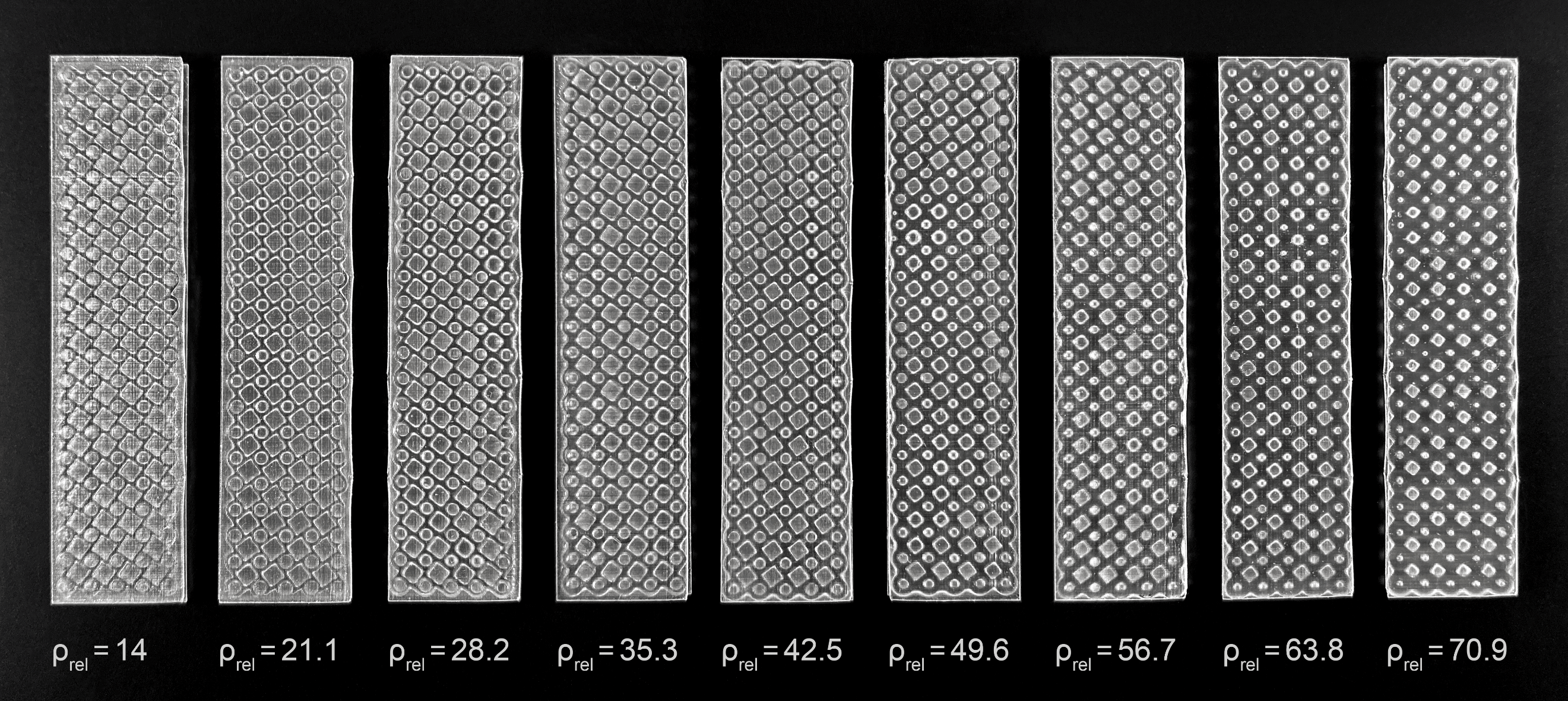}
\caption{Изготовленные композиты с тетракиральными сотами при \(d_{a} = 1.6\) мм и \(14 \leq \rho_{rel} \leq 70.9\) \% \dualcite{ref78}{{[}\bibnum{ref78}{]}}{{[}\bibnum{ref78}{]}}}\label{fig:48}
\end{figure}

\begin{figure}[H]
\centering
\includegraphics[keepaspectratio,width=5.51181in,height=3.67434in]{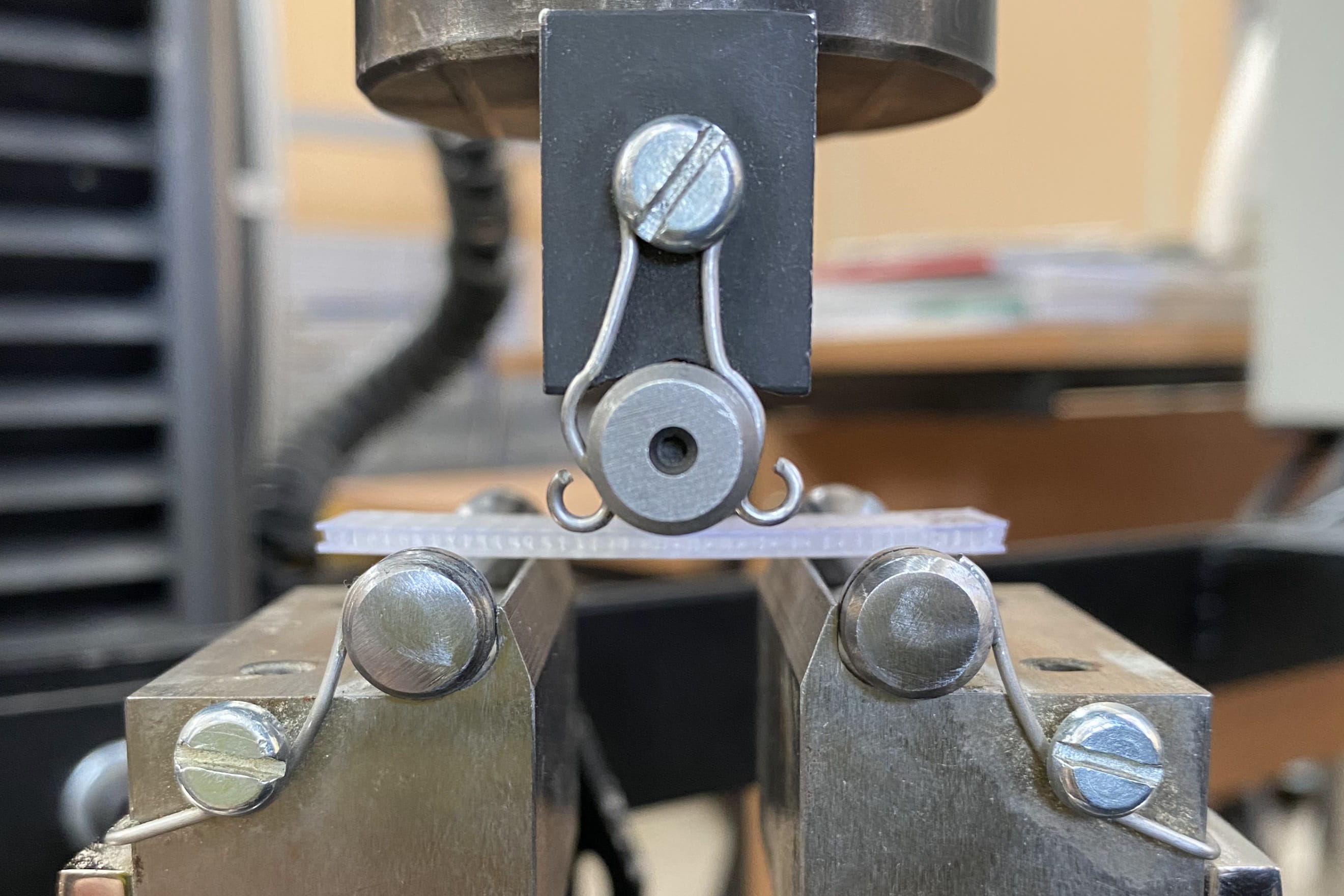}
\caption{Испытание композита в условиях трехточечного изгиба \dualcite{ref78}{{[}\bibnum{ref78}{]}}{{[}\bibnum{ref78}{]}}}\label{fig:49}
\end{figure}

Наряду с визуальным контролем качества изготовленных композитов осуществлялось сравнение их фактической массы с расчетными значениями, полученными в системе проектирования 3D-моделей (рисунок \ref{fig:47}). Средние отклонения фактической массы композитов от номинальных значений по CAD-модели при \(d_{a}\) = 1, 1.3, 1.6 и 1.9 мм составляют 17.9, 8.9, 6.2 и 5 \% соответственно, что в совокупности с визуальным контролем (рисунок \ref{fig:48}) свидетельствует о приемлемом качестве изготовления образцов. Следует отметить, что с увеличением дискретизации сот (уменьшением \(d_{a}\) от 1.9 до 1 мм) увеличивается общая длина контуров структуры сот в плоскости, что приводит к увеличению объема твердого тела напечатанных сот вследствие погрешности печати (рисунок \ref{fig:47}). Испытания композитных пластин на трехточечный изгиб проводились также с использованием оборудования Instron 5982 (рисунок \ref{fig:49}).

\vspace{-0.575\baselineskip}
\vspace{-0.575\baselineskip}

\begin{figure}[H]
\centering
\includegraphics[keepaspectratio,width=3.93973in,height=7.83464in]{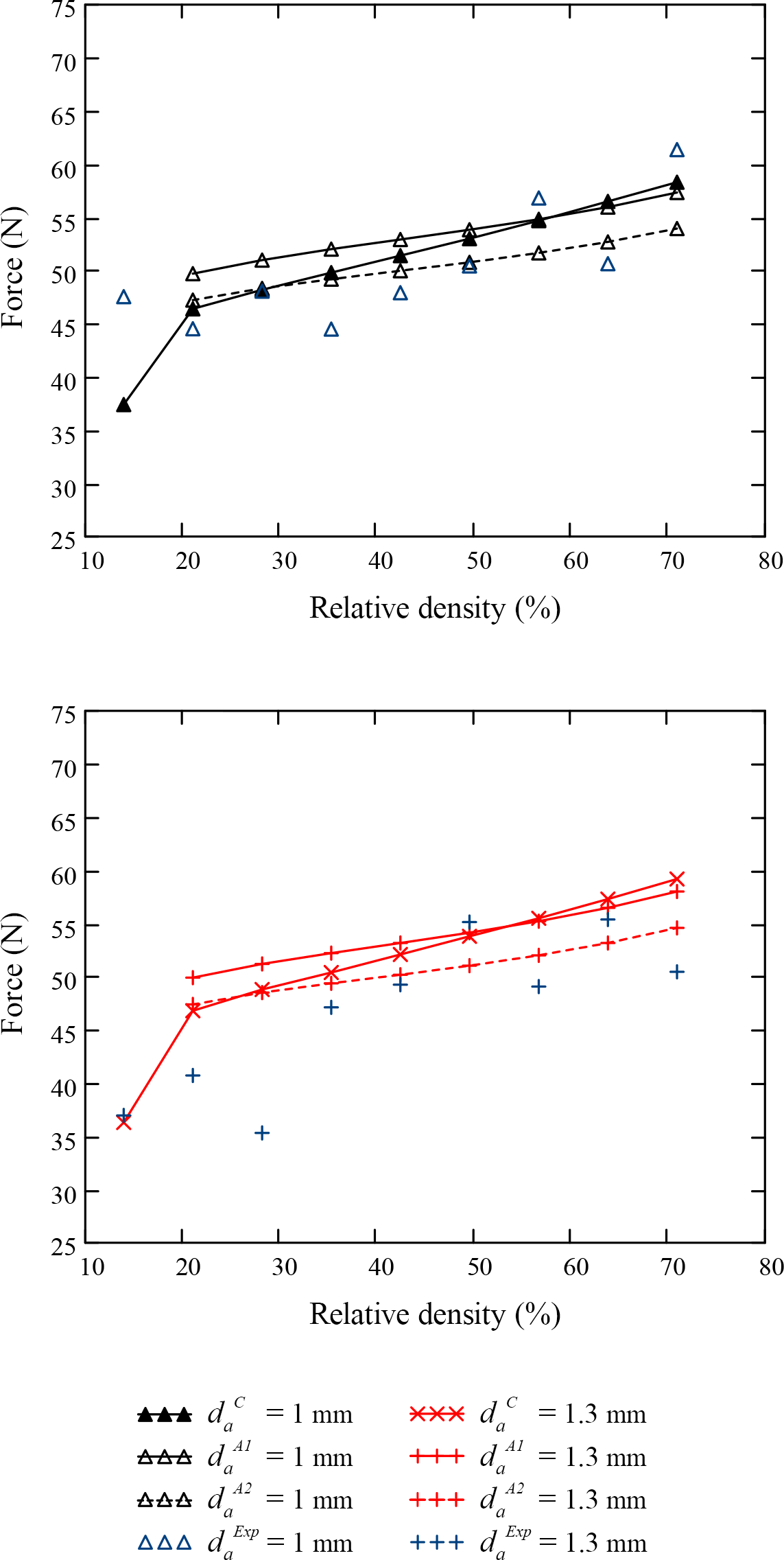}
\caption{Диаграммы зависимости \(F_{y}\) от \(\rho_{rel}\) при \(\sigma_{\max} = \sigma_{pl}\) в первой постановке численных и лабораторных экспериментов при \(d_{a} = 1\) и \(d_{a} = 1.3\) мм в условиях трехточечного изгиба, полученные с применением Comsol Multiphysics \dualcite{ref78}{{[}\bibnum{ref78}{]}}{{[}\bibnum{ref78}{]}} и Алгоритмов 1 и 2}\label{fig:50}
\end{figure}

\begin{figure}[H]
\centering
\includegraphics[keepaspectratio,width=3.93973in,height=7.83464in]{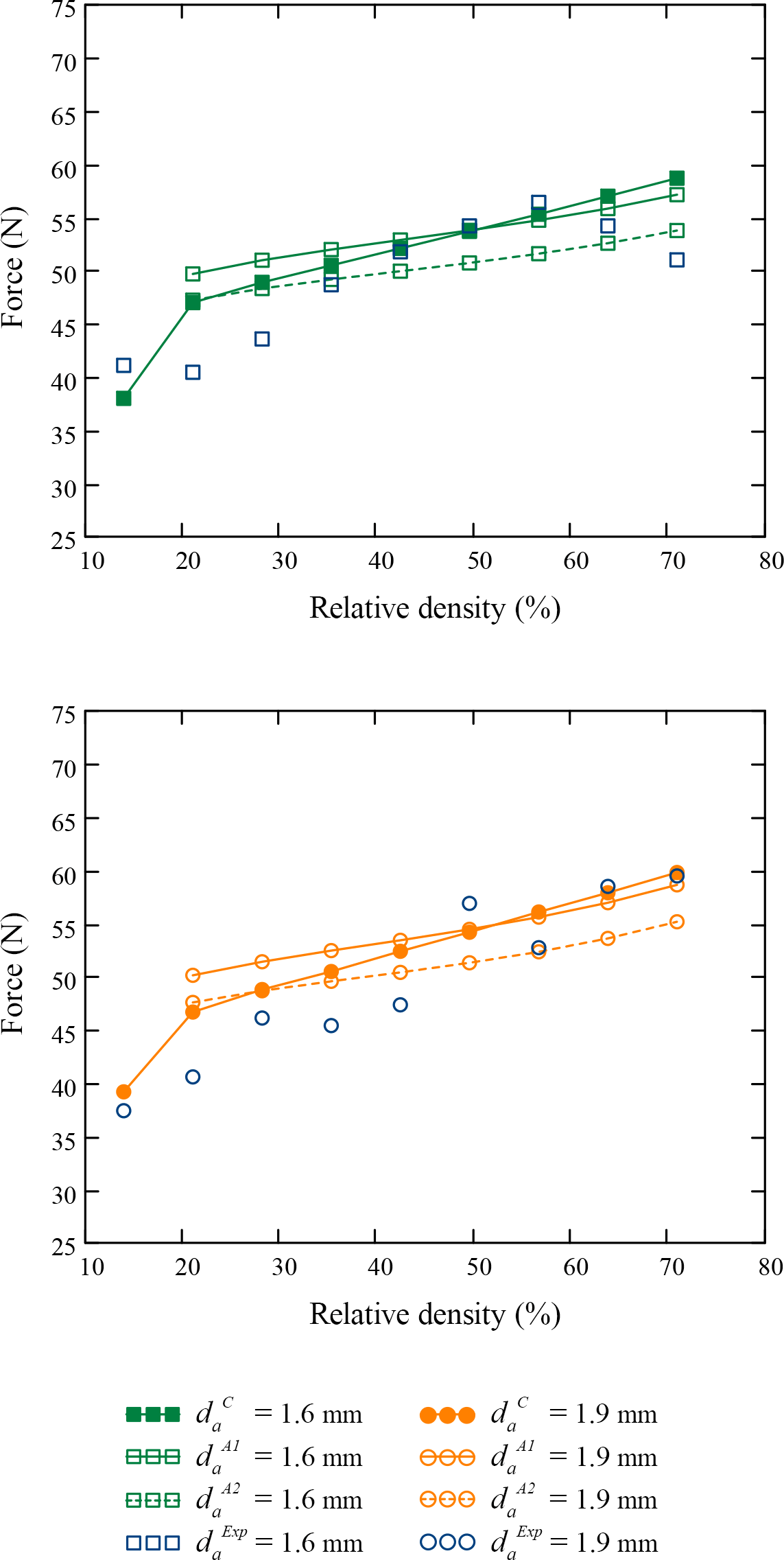}
\caption{Диаграммы зависимости \(F_{y}\) от \(\rho_{rel}\) при \(\sigma_{\max} = \sigma_{pl}\) в первой постановке численных и лабораторных экспериментов при \(d_{a} = 1.6\) и \(d_{a} = 1.9\) мм в условиях трехточечного изгиба, полученные с применением Comsol Multiphysics \dualcite{ref78}{{[}\bibnum{ref78}{]}}{{[}\bibnum{ref78}{]}} и Алгоритмов 1 и 2}\label{fig:51}
\end{figure}

\begin{figure}[H]
\centering
\includegraphics[keepaspectratio,width=6.29921in,height=4.25546in]{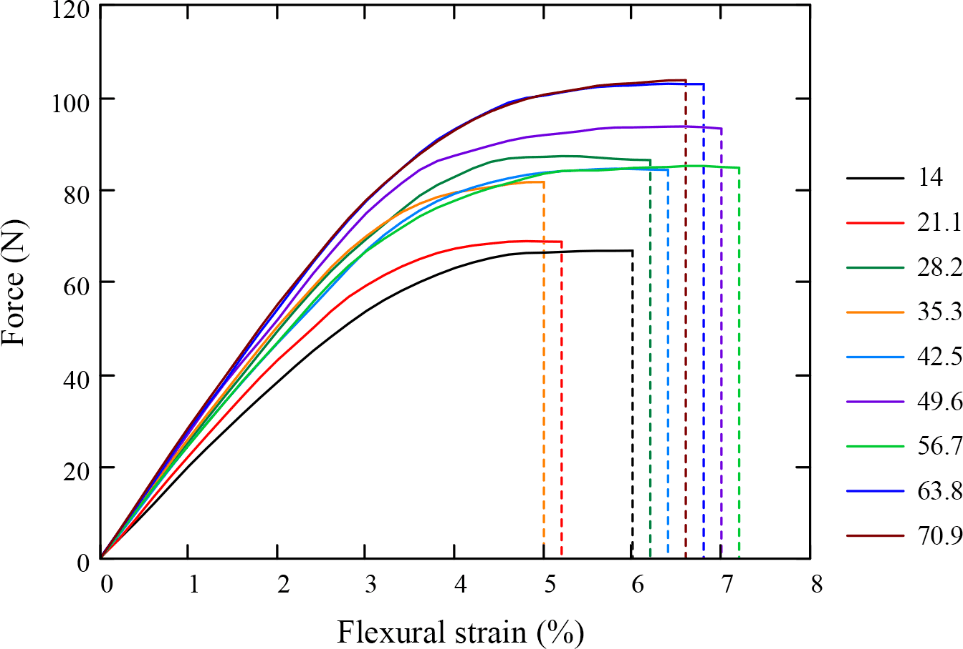}
\caption{Диаграмма зависимости \(F_{y}\) от \(\varepsilon_{f}\) композитов при \(d_{a} = 1.9\) мм и \(14 \leq \rho_{rel} \leq 70.9\) \% в первой постановке лабораторных экспериментов при трехточечном изгибе \dualcite{ref78}{{[}\bibnum{ref78}{]}}{{[}\bibnum{ref78}{]}}}\label{fig:52}
\end{figure}

Представлены диаграммы зависимости \(F_{y}\) от \(\rho_{rel}\) при \(\sigma_{\max} = \sigma_{pl}\) в первой постановке численных и лабораторных экспериментов при трехточечном изгибе композитных пластин (рисунки \ref{fig:50}, \ref{fig:51}). В рамках лабораторных испытаний \(\left( {d_{a}}^{Exp} \right)\) значения силы, при которых максимальные напряжения в композитах достигали предела пропорциональности, определялись по диаграмме зависимости \(F_{y}\) от относительной деформации при изгибе \(\varepsilon_{f}\) (рисунок \ref{fig:52}), которая рассчитывалась по стандарту ASTM D790-10 \dualcite{ref44}{{[}\bibnum{ref44}{]}}{{[}\bibnum{ref44}{]}}. При этом радиус пуансона и опор -- 5 и 2.5 мм соответственно, а скорость опускания пуансона -- 4 мм/мин.

\vspace{-0.575\baselineskip}
\vspace{-0.575\baselineskip}

\clearpage
\iflatexml\else\phantomsection\fi
\section{ЗАКЛЮЧЕНИЕ}

Основные результаты работы можно сформулировать следующим образом:

\iflatexml
\renewcommand{\labelenumi}{\arabic{enumi}.}
\begin{enumerate}
\fi
\fullwidthnumitem{1}{Построены эффективные математические конечно-элементные модели слоистой пластины с тетракиральными сотовыми прослойками с использованием совместных и несовместных конечных элементов для описания напряженно-деформированного состояния в условиях статического воздействия в рамках теории упругости.}
\fullwidthnumitem{2}{Разработаны алгоритмы численной реализации математического моделирования напряженно-деформированного состояния многослойных пластин с тетракиральными сотами в условиях статического изгиба посредством решения плоской задачи теории упругости методом конечных элементов с применением совместных и несовместных элементов.}
\fullwidthnumitem{3}{Выполнено численное моделирование напряженного состояния трехслойных пластин на основе первой постановки численных экспериментов при постоянной толщине слоев и варьируемой относительной плотности заполнителя методом конечных элементов посредством алгоритмов решения плоской задачи и путем трехмерного моделирования в Comsol Multiphysics.}
\fullwidthnumitem{4}{Выполнено численное моделирование напряженного состояния трехслойных пластин на основе второй постановки численных экспериментов при постоянном объеме твердого тела сот и варьируемой толщине заполнителя методом конечных элементов посредством алгоритмов решения плоской задачи и путем трехмерного моделирования в Comsol Multiphysics.}
\fullwidthnumitem{5}{В первой постановке численных экспериментов анализ результатов показал, что в условиях жесткого защемления пластин при постоянной нагрузке с увеличением относительной плотности сот от 14 до 71 \% разница между максимальными напряжениями в сотах с различной дискретизацией сначала интенсивно увеличивается, а затем уменьшается до минимального значения. Однако в условиях опирания с упругим поворотом на всем диапазоне изменения относительной плотности максимальные напряжения в сотах с различной дискретизацией имеют одинаковый разброс в пределах погрешности метода конечных элементов.}
\fullwidthnumitem{6}{Во второй постановке численных экспериментов анализ результатов показал, что в условиях жесткого защемления пластин при постоянной нагрузке с увеличением толщины заполнителя (уменьшением относительной плотности) немонотонно возрастают максимальные напряжения в тетракиральных сотах с различной дискретизацией. Однако в условиях опирания с упругим поворотом в сотах с различной дискретизацией наблюдается более узкий разброс максимальных напряжений со слабо выраженной закономерностью. При варьировании толщины пластин посредством изменения относительной плотности сот в условиях жесткого защемления можно получить минимальные значения максимального напряжения в композитах. При этом минимальное значение максимального напряжения у большинства композитных пластин находится в точке перехода от сотового заполнителя к сплошному слою при определенной толщине сот.}
\fullwidthnumitem{7}{В обеих постановках численных экспериментов при жестком защемлении снижение дискретизации сотовых структур (уменьшение числа элементарных ячеек) при постоянной относительной плотности приводит к снижению максимального напряжения в заполнителе для большинства значений относительной плотности в рассматриваемом диапазоне.}
\fullwidthnumitem{8}{Дискретизация сотового заполнителя не влияет на значения максимальных напряжений во внешних слоях композитных пластин в обеих постановках численных экспериментов с применением обоих граничных условий. Результаты моделирования в системе Comsol Multiphysics и алгоритма решения плоской задачи с использованием несовместных элементов демонстрируют хорошее соответствие, а применение совместных элементов приводит к более высокой погрешности в большинстве случаев.}
\fullwidthnumitem{9}{Результаты лабораторных испытаний трехслойных пластин с тетракиральным сотовым заполнителем в первой постановке экспериментов при трехточечном изгибе качественно соответствуют результатам конечно-элементного моделирования.}
\fullwidthnumitem{10}{Использование разработанных алгоритмов численного моделирования механического поведения многослойных композитов с тетракиральными сотами путем решения плоской задачи не требует проектирования детальных 3D-моделей и построения соответствующих конечно-элементных сеток. Это особенно актуально в отношении сотовых заполнителей, поскольку получение качественной сетки конечных элементов отдельной конфигурации сотовой структуры в трехмерной задаче, как правило, занимает дополнительное время, а дальнейшее использование объемной сетки требует повышенных вычислительных затрат. В предложенных алгоритмах решения плоской задачи тетракиральные соты представляются как непрерывная среда с усредненными свойствами, а упругие константы эквивалентной среды определяются аналитически через геометрические параметры сот. Этот подход позволяет эффективно анализировать слоистые композиты с тетракиральными сотами различной конфигурации без частого перестроения сетки плоских конечных элементов, при этом получаемые численные результаты сопоставимы с результатами трехмерного конечно-элементного моделирования.}
\iflatexml
\end{enumerate}
\fi

\clearpage
\renewcommand{\refname}{СПИСОК ЛИТЕРАТУРЫ}
\iflatexml
\else
\phantomsection\addcontentsline{toc}{section}{\refname}
\fi

\end{document}